\documentclass[12pt]{amsart}
\usepackage{amsfonts,amssymb,amscd,amsthm,amsmath,euscript}
\usepackage{graphpap}
\usepackage{epic}
\usepackage[matrix,arrow,curve]{xy}

 

\newtheorem{theorem}[subsection]{Theorem}
\newtheorem*{theorem*}{Theorem}
\newtheorem{lemma}[subsection]{Lemma}
\newtheorem{proposition}[subsection]{Proposition}
\newtheorem{corollary}[subsection]{Corollary}

\theoremstyle{definition}
\newtheorem{definition}[subsection]{Definition}
\newtheorem*{problem}{Problem}
\newtheorem{example}[subsection]{Example}
\theoremstyle{remark}
\newtheorem{remark}[subsection]{Remark}

\theoremstyle{definition}

\newcommand{\mt}[1]{\operatorname{#1}}
\newcommand{\EEE}{{\mathbb E}}
\newcommand{\DDD}{{\mathbb D}}
\newcommand{\AAA}{{\mathbb A}}
\newcommand{\QQ}{{\mathbb Q}}
\newcommand{\ZZ}{{\mathbb Z}}
\newcommand{\CC}{{\mathbb C}}
\newcommand{\OO}{{\mathcal O}}
\newcommand{\RR}{{\mathbb R}}
\newcommand{\PP}{{\mathbb P}}

\newcommand{\NN}{{\mathbb N}}
\newcommand{\FFF}{{\mathbb F}}
\newcommand{\PPP}{{\EuScript{P}}}

\newcommand{\Supp}{\mt{Supp}}

\newcommand{\Sing}{\mt{Sing}}
\newcommand{\Diff}{\mt{Diff}}

\newcommand{\D}{{\Delta}}
\newcommand{\Hom}{\mt{Hom}}
\newcommand{\Exc}{\mt{Exc}}
\newcommand{\codim}{\mt{codim}}
\newcommand{\mult}{\mt{mult}}

\newcommand{\NE}{\overline{\mt{NE}}}
\newcommand{\down}[1]{\llcorner #1 \lrcorner}

\newcommand{\fr}[1]{\{ #1\}}
\newcommand{\G}{\Gamma}

\newcommand{\LCS}{\mt{LCS}}

\title{Blow-ups of three-dimensional toric singularities}
\author{S.~A.~Kudryavtsev}

\date{}

\address{Scientific Research Institute of Precise Instruments, 453 Department,
Dekabristov str. 51, 127490 Moscow, Russia}

\email{kudryavtsev@myrambler.ru}

\dedicatory{In memory of Vasily Alexeevich Iskovskikh}

\begin{document}
\begin{abstract}
The purely log terminal blow-ups of three-dimensional terminal
toric singularities are described.
The three-dimensional divisorial contractions
$f\colon (Y,E)\to (X\ni P)$ are described provided that
$\Exc f=E$ is an irreducible divisor, 
$(X\ni P)$ is a toric terminal singularity,
$f(E)$ is a toric subvariety and
$Y$ has canonical singularities.
\end{abstract}

\maketitle

\section*{\bf Introduction}
Let $(X \ni P)$ be a log canonical singularity and let $f\colon Y \to X$ be its blow-up.
Suppose that the exceptional locus of $f$ consists of only one irreducible divisor:
$\Exc f=E$. Then $f\colon (Y,E) \to (X \ni P)$
is called a {\it purely log terminal blow-up}, {\it canonical blow-up} or {\it terminal blow-up},
if 1), 2) or 3) are satisfied respectively:
1) $K_Y+E$ is plt and $-E$ is $f$-ample;
2) $-K_Y$ is $f$-ample and $Y$ has canonical singularities;
3) $-K_Y$ is $f$-ample and $Y$ has terminal singularities.

The definition of plt blow-up
implicitly requires that the divisor
$E$ be $\QQ$-Cartier. Hence $Y$ is a
$\QQ$-gorenstein variety.
By the inversion of adjunction (see \cite[Theorem 17.6]{Koetal})
$K_E+\Diff_E(0)=(K_Y+E)|_E$ is klt.

The importance of study of purely log terminal blow-ups is that: some very important questions of
birational geometry for $n$-dimensional varieties, contractions can be reduced to the smaller dimension $n-1$,
using purely log terminal blow-ups (for instance, see
the papers \cite{PrLect}, \cite{Sh2}, \cite{PrSh} and \cite{PrEll}).
In dimension two, purely log terminal blow-ups are completely classified and the classification of two-dimensional 
non-divisorial log terminal extremal contractions of local type
is obtained using them \cite{PrLect}.
For three-dimensional varieties the first similar problem is to get the same explicit geometric classification of three-dimensional
Mori contraction of local type as in two-dimensional case.
The next problem is the first difficulty to realize this approach.
\begin{problem}
Describe the class of all log del Pezzo
surfaces, generic $\PP^1$-fibrations which can be the exceptional
divisors of some purely log terminal blow-ups of
three-dimensional terminal singularities.
\end{problem} 

Suppose that $f(E)=P$ is a point. 
Then we solve this problem in the case of terminal toric singularities 
(Theorem \ref{dim3plt}). Moreover we obtain the description of plt blow-ups of
$\QQ$-factorial three-dimensional toric singularities (Theorem \ref{dim32plt}).
Purely log terminal and canonical blow-ups are divided into toric and non-toric blow-ups up to analytic isomorphism.
The study of non-toric plt blow-ups is reduced to the description of plt triples
$(S,D,\G)$ in dimension two
(Definition \ref{triple}). 

Also we obtain the description of canonical blow-ups of three-dimensional terminal toric singularities
(Theorem \ref{dim3can}).
The study of non-toric canonical blow-ups is reduced to the description of the following two interrelated objects: a) toric canonical blow-ups of $(X\ni P)$ and
b) some triples $(S,D,\G)$ in dimension two.

Immediate corollary of Theorem \ref{dim3can} is that the terminal blow-ups of three-dimensional terminal toric singularities are toric up to analytic isomorphism.
This corollary was proved in the papers \cite{Kaw1}, 
\cite{Kawakita1} and \cite{Corti} by another methods.

Suppose that $f(E)$ is a one-dimensional toric subvariety (curve) of the toric singularity $(X\ni P)$.
Then the description of plt and canonical blow-ups is given in Theorems 
\ref{dim31_plt}, \ref{dim32_plt}, \ref{dim31_can} and in Corollary \ref{dim32_can}.

I am grateful to Professors  Yu.G.~Prokhorov and I.A.~Cheltsov for valuable advices.

\section{\bf Preliminary results and facts}
All varieties are algebraic and are assumed to be defined over
$\CC$, the complex number field.
The main definitions, notations and notions used in the paper are
given in
\cite{Koetal}, \cite{KMM}, \cite{PrLect}.
See \cite[Section 3.10]{BCHM} on minimal model program with scaling.
The definition of $\Diff$ and its main properties are given in the papers \cite[\S 3]{Sh1}, \cite[Chapter 16]{Koetal}.
By $(X\ni P)$ denote the algebraic germ of the variety $X$ at the point $P$.

A smooth point is a special case of {\it singularity} by our definition. For example, Du Val singularity of type $\AAA_0$ is a smooth point.

Let $f\colon Y\dashrightarrow X$ be a birational map and let $D$ be a divisor on the variety $X$. 
By $D_Y$ denote the proper transform of $D$ on the variety $Y$.
If $Y=\widetilde X$, $Y=X'$ or $Y=\overline X$, then for notational convenience we use the notation
$\widetilde D=D_{\widetilde X}$,
$D'=D_{X'}$ or $\overline D=D_{\overline X}$ respectively. 
The similar notation is used for subvarieties of $X$.

{\it The contraction} $f\colon Y\to X$ is a projective morphism of the normal variety such that $f_*\OO_Y=\OO_X$. {\it A blow-up} is a birational divisorial contraction. 
A $\QQ$-{\it factoriality} means analytical $\QQ$-factoriality in this paper.

The proper irreducible subvariety
$\G$ of $X$ is said to be {\it a center of canonical singularities} of $(X,D)$, if there exist
the birational morphism $f\colon Y\to X$ and the exceptional divisor $E\subset Y$ such that $\G=f(E)$ and $a(E,D)\leq 0$. 
The set of canonical singularity centers of
$(X,D)$ and $X$ is denoted by $\mt{CS}(X,D)$ and $\mt{CS}(X)$ respectively.

By our definition the {\it toric varieties, toric morphisms} are considered up to analytic isomorphism (analytical identification), if
they are not explicitly defined by fans. Shokurov's (hypothetical) criterion on the characterization
of toric varieties is formulated in \cite[Chapter 6]{Sh2}.
By definition of {\it weighted blow-up}, its center is a point always, that is, 
its every weight is positive.

We write all singularities of surface in brackets. For example, the notation
$S(\AAA_1+\frac15(1,2))$ means that the surface $S$ has two singular points 
of types $\AAA_1$ and $\frac15(1,2)$ exactly.

We actively use a structure of the local toric conic bundle $f\colon S\to (C\ni P)$, where $\dim S=2$ and $\rho(S/C)=1$.
By \cite[Lemma 7.1.11]{PrLect} the surface $S$ has two singularities of types $\frac1r(1,q)$ and $\frac1r(1,-q)$ over the point $P$ only,
where $r\geq 1$.

\begin{proposition}\label{lemma1}\cite[Lemma 6.2]{Koetal}
Let $f_i\colon Y_i\to X$ be two divisorial contractions of normal varieties, where
$\Exc f_i=E_i$ are irreducible divisors and $-E_i$ are $f_i$-ample divisors.
If $E_1$ and $E_2$ define the same discrete valuation of the function field
$\mathcal K(X)$, then
the contractions
$f_1$ and $f_2$ are isomorphic.
\end{proposition}

\begin{proposition}\label{lemma2}
Let $f_i\colon Y_i\to (X\ni P)$ be two divisorial contractions to a point $P$,
where $\Exc f_i=E_i$ are irreducible divisors. Suppose that the varieties $Y_i$, $X$
have log terminal singularities, $E_1$ and $E_2$
define the same discrete valuation of the function field
$\mathcal K(X)$, the divisor $-E_1$ is $f_1$-ample,
the divisor $-E_2$ is not $f_2$-ample.
Then there exists the small flopping contraction
$($with respect to $K_{Y_2})$ $g\colon Y_2\to Y_1$ such that $f_2$ and $f_1\circ g$ are isomorphic.
\begin{proof}
Let $K_{Y_2}=f_2^*K_X+aE_2$. If $a>0$, then we put $L=-K_{Y_2}$. If
$a\le 0$, then we put $L=-(K_{Y_2}+(-a+\varepsilon)E_2)$, where $\varepsilon$ is a sufficiently
small positive rational number. Since $-E_2$ is a $f_2$-nef divisor, then
the linear system $|nL|$ is free over $X$ for $n\gg 0$ and gives a contraction
$g\colon Y_2\to Y'_2$ over $X$ by the base point free theorem \cite[Remark 3.1.2]{KMM}.
A curve $C$ is exceptional for $g$ if and only if
$L\cdot C=E_2\cdot C=K_{Y_2}\cdot C=0$. Therefore
$g$ is a flopping contraction and $Y'_2\cong Y_1$ by Proposition \ref{lemma1}.
\end{proof}
\end{proposition}

The next example shows the idea of Proposition
\ref{lemma2}.
\begin{example}\label{example}
Let $(X\ni P)\cong (\{x_1x_2+x_3^2+x_4^4=0\}\subset (\CC^4_{x_1x_2x_3x_4},0))$. Consider the divisorial contraction
$f_1\colon Y_1\to (X\ni P)$ induced by the blow-up of the maximal ideal of the point
$(\CC^4 \ni 0)$. Then $\Exc f_1\cong \PP(1,1,2)$, the variety $Y_1$ has only one singular point denoted by $Q$, and 
$(Y_1\ni Q)\cong (\{y_1y_2+y_3^2+y_4^2=0\}\subset (\CC^4_{y_1y_2y_3y_4},0))$. 
This singularity is not $\QQ$-factorial and let 
$g\colon Y_2\to (Y_1\ni P)$ be its $\QQ$-factorialization.
We obtain the divisorial contraction
$f_2\colon Y_2\to (X\ni P)$, where $Y_2$ is a smooth 3-fold,
$\Exc f_2\cong \FFF_2$, and $-K_{Y_2}$ is not a $f_2$-ample divisor.
\end{example}

\begin{definition}\label{defcan}
Let $(X \ni P)$ be a log canonical singularity and let $f\colon Y \to X$ be its blow-up.
Suppose that the exceptional locus of $f$ consists of only one irreducible divisor:
$\Exc f=E$. Then $f\colon (Y,E) \to (X \ni P)$
is called a {\it canonical blow-up} if
$-K_Y$ is $f$-ample and $Y$ has canonical singularities.
Note that the definition of canonical blow-up implies that $(X\ni P)$ is a canonical singularity. 
The canonical blow-up is said to be {\it a terminal blow-up} if $Y$ has terminal singularities.
\end{definition}

\begin{remark}\label{canpropert}
Using the notation of Definition \ref{defcan}, we have the following properties
of canonical blow-ups.
\begin{enumerate}
\item[1)] The definition of canonical (resp. terminal) blow-up implies easily that $(X\ni P)$ is a canonical (resp. terminal) singularity.
\item[2)] The divisor $-E$ is $f$-ample and $a(E,0)>0$.
\item[3)] Let $f_i\colon (Y_i,E_i)\to (X \ni P)$ be two canonical blow-ups.
If $E_1$ and $E_2$ define the same discrete valuation of the function field
$\mathcal K(X)$ then
the blow-ups $f_1$ and $f_2$ are isomorphic by Proposition $\ref{lemma1}$.
\item[4)] Let $(X \ni P)$ be a
$\QQ$-factorial singularity. Then $Y$ is a $\QQ$-factorial variety also, $\rho(Y/X)=1$ and $\rho(E)=1$ \cite[\S 5]{IshiiPr}.
\end{enumerate}
\end{remark}

\begin{theorem}\label{caninductive}
Let $(X\ni P)$ be a canonical singularity and $(X\ni P,D)$ be a pair with canonical singularities, where $D$ is a boundary. 
Assume that $a(E,D)=0$ and $a(E,0)>0$ for some irreducible exceptional divisor $E$.
Then there exists a canonical blow-up such that its exceptional divisor and $E$
define the same discrete valuation of the function field $\mathcal K(X)$.
Moreover, if $E$ is a unique exceptional divisor with $a(E,D)=0$ then its canonical blow-up is a terminal blow-up.
\begin{proof}
By Proposition 21.6.1 of the paper \cite{Koetal} we consider the birational contraction
$\widetilde f\colon (\widetilde Y,\widetilde E)\to (X\ni P)$ with the following three properties:
\par 1) $\widetilde E$ is a unique irreducible exceptional divisor of $\Exc \widetilde f$;
\par 2) $\widetilde E$ and $E$ define the same discrete valuation of the function field $\mathcal K(X)$;
\par 3) if $(X\ni P)$ is $\QQ$-factorial then $\rho(\widetilde Y/X)=1$ and $\Exc \widetilde f=\widetilde E$.

The proof of Proposition 21.6.1 of \cite{Koetal} holds in any dimension since we can apply MMP with scaling to prove it.
Let $\widetilde f$ be not the required canonical blow-up.
If $\Exc \widetilde f=\widetilde E$ then by Proposition \ref{lemma2} we have $\widetilde f\cong f\circ g$, where $f$ is the required blow-up.
Consider the remaining case when
$\Exc \widetilde f=\widetilde E \cup \Delta$, where $\Delta \ne \emptyset$ and
$\codim_{\widetilde Y}\Delta \ge 2$. Let $H$ be a general Cartier divisor containing the set $\widetilde f(\Exc \widetilde f)$.
Then $K_{\widetilde Y}+D_{\widetilde Y}+\varepsilon H_{\widetilde Y}\equiv -\varepsilon a \widetilde E$ over $X$, where $a>0$.
For $0<\varepsilon\ll 1$ we apply $K_{\widetilde Y}+D_{\widetilde Y}$ -- MMP with scaling of $H_{\widetilde Y}$.
We obtain a birational map
$\varphi\colon \widetilde Y\dashrightarrow Y'$, which is a composition of log flips, and we also obtain
a divisorial contraction $f'\colon Y' \to X$ such that $\Exc f'=E'$, where
$E'$ is an irreducible divisor. Therefore, by Proposition \ref{lemma2} we have the required canonical blow-up.
\end{proof}
\end{theorem}

\begin{definition} \label{defplt}
Let $(X \ni P)$ be a log canonical singularity and let $f\colon Y \to X$ be its blow-up.
Suppose that the exceptional locus of $f$ consists of only one irreducible divisor:
$\Exc f=E$. Then $f\colon (Y,E) \to (X \ni P)$
is called a {\it purely log terminal blow-up} if the divisor
$K_Y+E$ is purely log terminal and $-E$ is $f$-ample.
\end{definition}

\begin{remark}
Definition
\ref{defplt} implicitly requires that the divisor
$E$ be $\QQ$-Cartier. Hence $Y$ is a
$\QQ$-gorenstein variety.
By the inversion of adjunction
$K_E+\Diff_E(0)=(K_Y+E)|_E$ is klt.
\end{remark}

\begin{remark}\label{pltpropert}
Using the notation of Definition \ref{defplt} we have the following properties
of purely log terminal blow-ups.
\begin{enumerate}

\item[1)] The variety $f(E)$ is normal \cite[Corollary 2.11]{Pr2}.

\item[2)]  If $(X \ni P)$ is a log terminal singularity then
$-(K_Y+E)$ is a $f$-ample divisor. 
A purely log terminal blow-up of
log terminal singularity always exists 
by Theorem 1.5 of \cite{Kud1} since we can apply MMP with scaling to prove it (see also Theorem \ref{inductive}).

\item[3)] If $(X \ni P)$ is a strictly log canonical singularity then
$a(E,0)=-1$. A purely log terminal blow-up of strictly
log canonical singularity exists if and only if
there is only one exceptional divisor with discrepancy
$-1$ \cite[Theorem 1.9]{Kud1}, since we can apply MMP with scaling to prove Theorem 1.9 of \cite{Kud1}.

\item[4)] If $(X \ni P)$ is a
$\QQ$-factorial singularity then $Y$ is a $\QQ$-factorial variety also, 
$\rho(Y/X)=1$ and $\rho(E)=1$ \cite[Remark 2.2]{Pr2}, \cite[\S 5]{IshiiPr}.
Hence, for $\QQ$-factorial singularity we can omit the requirement that
$-E$ be $f$-ample in
Definition \ref{defplt} because it holds automatically.

\item[5)] 
Let $f_i\colon (Y_i,E_i)\to (X \ni P)$ be two purely log terminal blow-ups.
If $E_1$ and $E_2$ define the same discrete valuation of the function field
$\mathcal K(X)$ then
the blow-ups $f_1$ and $f_2$ are isomorphic by Proposition $\ref{lemma1}$.

\item[6)] Let $-E$ be not a $f$-ample divisor in Definition $\ref{defplt}$.
Then such blow-up can differ from some plt blow-up by a small flopping contraction only
(with respect to the canonical divisor
$K_Y$) \cite[Corollary 1.13]{Kud1}. This statement is similar to Proposition $\ref{lemma2}$.

\item[7)] Let $f\colon (Y,E) \to (X \ni P)$ be a toric blow-up of a toric $\QQ$-gorenstein
singularity. Assume that $Y$ is a $\QQ$-gorenstein variety and $\Exc f=E$ is an irreducible divisor. 
It is obvious that 
$K_Y+E$ is a plt divisor. Therefore, if $(X\ni P)$ is $\QQ$-factorial singularity then $f$ is a plt blow-up.
\end{enumerate}
\end{remark}

\begin{theorem}\cite[Theorem 1.5]{Kud1}, \cite[Proposition 2.9]{Pr2}\label{inductive}
Let $X$ be  a kawamata log terminal variety and let $D\ne 0$ be a boundary on $X$ such that
$(X,D)$ is log canonical, but not purely log terminal.
Then there exists an inductive blow-up $f:Y\to X$
such that:
\begin{enumerate}
\item the exceptional locus of $f$ contains only one irreducible divisor
$E\ (\Exc(f)=E)$;
\item $K_Y+E+D_Y=f^*(K_X+D)$ is log canonical;
\item $K_Y+E+(1-\varepsilon)D_Y$ is purely log terminal and anti-ample over $X$
for any $\varepsilon>0$;
\item if $X$ is $\QQ$-factorial then $Y$ is also
$\QQ$-factorial and $\rho(Y/X)=1$.
\end{enumerate}
\begin{proof}
The proofs of \cite[Theorem 1.5]{Kud1}, \cite[Proposition 2.9]{Pr2} hold in any dimension since we can apply MMP with scaling to prove them.
\end{proof}
\end{theorem}

\begin{remark}
Inductive blow-up is a plt blow-up.
Conversely, for any plt blow-up 
$f\colon (Y,E) \to (X \ni P)$ there exists a pair $(X,D)$ such that $f$ is its inductive blow-up. Indeed,
put $D=f(\frac1nD_Y)$, where $D_Y\in |-n(K_Y+E)|$ is a general element for $n\gg 0$.
\end{remark}

\begin{definition}
Let $(X/Z,D)$ be a contraction of varieties, where $D$ is a subboundary. Then a
\textit{$\QQ$-complement} of $K_X+D$ is an effective $\QQ$-divisor $D'$
such that $D'\ge D$, $K_X+D'$ is log canonical and
$K_X+D'\sim_{\QQ} 0/Z$ for some $n\in\NN$.
\end{definition}

\begin{definition}
Let $(X/Z,D)$ be a contraction of varieties.
Let $D=S+B$ be a subboundary on
$X$ such that $B$ and $S$ have no common components, $S$ is an
effective integral divisor and $\down{B}\le 0$. Then we say that
$K_X+D$ is \textit{$n$-complementary} if there is a $\QQ$-divisor
$D^+$ (called an $n$-{\it complement}) such that
\begin{enumerate}
\item
$n(K_X+D^+)\sim 0/Z$ (in particular, $nD^+$ is an integral divisor);
\item
the divisor $K_X+D^+$ is log canonical;
\item
$nD^+\ge nS+\down{(n+1)B}$.
\end{enumerate}
The divisor $K_X+D^+$ is also called an $n$-{\it complement}.
\end{definition}

\begin{definition} For $n\in \NN$ put
$$
\PPP_n=\{ a \mid 0\le a \le 1,\ \down{(n+1)a}\ge n a \}.
$$
\end{definition}

\begin{proposition}\cite[Lemma 5.4]{Sh1}\label{compldown}
Let $f\colon X\to Y$ be a birational contraction and let $D$ be a
subboundary on $X$. Assume that $K_X+D$ is $n$-complementary for
some $n\in\NN$. Then $K_Y+f(D)$ is also $n$-complementary.
\end{proposition}

\begin{proposition}\cite[Lemma 4.4]{Sh2}\label{complup}
Let $f\colon X\to Z$ be a birational contraction of
varieties and let $D$ be a subboundary on $X$.
Assume that
\begin{enumerate}
\item the divisor $K_X+D$ is $f$-nef;
\item the coefficient of every non-exceptional component of
$D$ meeting $\Exc f$ belongs to $\PPP_n$;
\item the divisor $K_Z+f(D)$ is
$n$-complementary.
\end{enumerate}
Then the divisor $K_X+D$ is also $n$-complementary.
\end{proposition}

\begin{proposition}\cite[Proposition 4.4.1]{PrLect}\label{complind}
Let $f\colon X\to (Z\ni P)$ be a contraction and  $D$ be
a boundary on $X$.
Put $S=\down{D}$ and
$B=\fr{D}$. Assume that
\begin{enumerate}
\item
the divisor $K_X+D$ is purely log terminal;
\item
the divisor $-(K_X+D)$ is $f$-nef and $f$-big;
\item
$S\ne 0$ near $f^{-1}(P)$;
\item
every coefficient of $D$ belongs to $\PPP_n$.
\end{enumerate}
Further, assume that near $f^{-1}(P)\cap S$ there exists an
$n$-complement $K_S+\Diff_S(B)^+$ of $K_S+\Diff_S(B)$. Then near
$f^{-1}(P)$ there exists an $n$-complement $K_X+S+B^+$ of
$K_X+S+B$ such that $\Diff_S(B)^+=\Diff_S(B^+)$.
\end{proposition}

\section{\bf Toric blow-ups}

We refer the reader to \cite{Oda} for the basics
of toric geometry.
\begin{definition}
Let $N$ be the lattice $\ZZ^n$ in the vector linear space $N_{\RR}=N\otimes_{\ZZ}\RR$ and
$M$ be its dual lattice $\Hom_{\ZZ}(N, \ZZ)$ in the vector linear space
$M_{\RR}=M\otimes_{\ZZ}\RR$. We have a canonical pairing
$\langle\ ,\ \rangle \colon N_{\RR}\times M_{\RR}\to \RR$.
\par
For a fan $\Delta$ in $N$ the corresponding
toric variety is denoted by $T_{N}(\D)$.
For a $k$-dimensional cone $\sigma\in\D$ the
closure of corresponding orbit is denoted by
$V(\sigma)$. This is a closed subvariety of
codimension $k$ in  $T_N(\D)$.
\end{definition}

\begin{example}\label{ex1}
1) Let the vectors $e_1,\ldots,e_n$ be a $\ZZ$-basis of $N$, where $n\geq 2$.
Consider the cone
$$
\sigma=\RR_{\geq 0}e_1+\ldots+\RR_{\geq 0}e_{n-1}+\RR_{\geq 0}(a_1e_1+\ldots+a_{n-1}e_{n-1}+re_n).
$$
Let the fan $\D$ consists of the cone $\sigma$ and its faces.
Then the affine toric variety
$T_{N}(\D)$ is the quotient space $(\CC^n\ni 0)/\ZZ_r$ with the action $\frac1r(-a_1,\ldots,-a_{n-1},1)$.
\par
2) Let
$$
\sigma=\langle e_1,e_2,e_3,e_4 \rangle = \langle (1,0,0),(0,1,0),(0,0,1),(1,1,-1) \rangle
$$
for the lattice $N\cong \ZZ^3$. Let the fan $\D$ consists of the cone $\sigma$ and its faces.
The affine toric variety $(X\ni P)=T_{N}(\D)$ is a three-dimensional non-degenerate quadratic cone in $\CC^4$.
Let
$$
\D^1=\{\langle e_1,e_2,e_3\rangle, \langle e_1,e_2,e_4\rangle, \text{their faces}\}
$$
and
$$
\D^2=\{\langle e_1,e_3,e_4\rangle, \langle e_2,e_3,e_4\rangle, \text{their faces}\}.
$$
Then the birational contractions $\psi_i\colon T_N(\D^i)\to T_N(\D)$ are small resolutions for $i=1$, $2$, and $\Exc\psi_1=V(\langle e_1,e_2\rangle)$, 
$\Exc\psi_2=V(\langle e_3,e_4\rangle)$. The birational map $T_N(\D^1)\dashrightarrow T_N(\D^2)$ is a flop.
\par
Let $f\colon (Y,E) \to (X \ni P)$ be a toric blow-up, where $Y$ is $\QQ$-gorenstein, $\Exc f=E$ is an irreducible divisor.
Then $f$ is a plt blow-up.
Let us prove it. The divisor $K_Y+E$ is plt. Let $a=(a_1,a_2,a_3)$ be a primitive vector defining $f$.
Consider any three-dimensional cone $\sigma'$ giving non-$\QQ$-factorial singularity of subdivision of the cone $\sigma$ by $a$. Then the cone $\sigma'$ gives non-$\QQ$-gorenstein singularity by Proposition 4.3 (i) \cite{Reid}, since
there is no any vector $m\in M_{\QQ}$ such that  
$\langle m, e_i \rangle=1$ for every $i$ and $\langle m, a \rangle=1$. Hence $-E$ is a $f$-ample divisor. 
This completes the proof.

Let $f(E)=P$. Then $Y=T_N(\widetilde \D)$ and
$$
\widetilde \D=\{\langle e_1,e_3,a\rangle, \langle e_1,e_4,a\rangle, \langle e_2,e_3,a\rangle, \langle e_2,e_4,a\rangle, \text{their faces}\},
$$
where $a=(a_1,a_2,a_3)$, $\gcd(a_1,a_2,a_3)=1$, $a_1>0$, $a_2>0$, $a_1+a_3>0$ and $a_2+a_3>0$.

Obviously, the converse is also true. Any such vector $a$ defines a plt blow-up.

Let $f(E)=C$ and $\dim C=1$. 
Then, up to a permutation of the faces of the cone $\sigma$ we have $C=\langle e_2,e_3 \rangle$,
$Y=T_N(\widehat \D)$ and
$$
\widehat \D=\{\langle e_2,e_4,a\rangle, \langle e_1,e_3,a\rangle, \langle e_1,e_4,a\rangle, \text{their faces}\},
$$
where $a=(0,a_2,a_3)$, $\gcd(a_2,a_3)=1$, $a_2>0$, $a_3>0$.

Obviously, the converse is also true. Any such vector $a$ defines a plt blow-up. 

The variety $Y$ has the singularities $\frac1{a_3}(0,-a_2,1)$, $\frac1{a_2}(0,1,-a_3)$, $\frac1{a_2+a_3}(-a_3,-a_2,1)$.
The surface $E$ is a toric conic bundle, $\rho(E/C)=2$,
the single singular point of $E$ (with a center of the third singularity of $Y$) 
has type $\AAA_{a_2+a_3-1}$ and $\Diff_E(0)=\frac{a_2-1}{a_2}E_1+\frac{a_3-1}{a_3}E_2$, where $E_1$, $E_2$ are corresponding sections.

We will calculate a structure of 
$f$ by the following way (for convenience). 
Let us consider $(X\ni P)\subset (\CC^4,0)$ as the embedding $\{x_1x_2+x_3x_4=0\}\subset (\CC_{x_1x_2x_3x_4}^4,0)$.
The weighted blow-up of $(\CC^4,0)$ with weights $w=(w_1,w_2,w_3,w_4)$ provided that $w_1+w_2=w_3+w_4$ induces
a toric blow-up $f'\colon (Y',E') \to (X \ni P)$, where 
$$
\Exc f'=E'\cong\{x_1x_2+x_3x_4\subset\PP_{x_1x_2x_3x_4}(w_1,w_2,w_3,w_4)\} -
$$
is an irreducible divisor. If put $w_1=a_1+a_3$, $w_2=a_2$, $w_3=a_2+a_3$ and $w_4=a_1$, then we can easily
compare the natural affine covers of
$Y$ and $Y'$ and prove that $f$ and $f'$ are isomorphic blow-ups.
Note that $C=\{x_1=x_2=x_3=0\}$ in the case $C=f'(E')$.
\end{example}

\begin{proposition}\cite[pages 36-37]{Oda}\label{termtoric}
The following statements are satisfied:
\begin{enumerate}
\item[1)] $(X\ni P)$ is a three-dimensional $\QQ$--factorial toric terminal singularity if and only if
$(X\ni P)\cong (\CC^3\ni 0)/\ZZ_r(q,-1,1)$, where $\gcd(r,q)=1;$
\item[2)] $(X\ni P)$ is a three-dimensional non-$\QQ$--factorial toric terminal singularity if and only if
$(X\ni P)\cong (\{x_1x_2+x_3x_4=0\}\subset (\CC^4_{x_1x_2x_3x_4},0))$.
\end{enumerate}
\end{proposition}

\begin{theorem} \label{cancyclic3dim}\cite{Mor2}
Let $(X\ni P)$ be a three-dimensional cyclic singularity of type $\frac1r(a_1,a_2,a_3)$.
Then $(X\ni P)$ is a canonical singularity if and only if one of the following holds: 
\\
$1) \ a_1+a_2+a_3\equiv 0(\mt{mod} r);$
\\
$2)\ a_i+a_j\equiv 0(\mt{mod} r)$ for some $i\ne j;$
\\
$3)$ $(X\ni P)$ has type $\frac19(1,4,7)$ or type $\frac1{14}(1,9,11)$.
\end{theorem}

\begin{proposition}\label{cantoricdim1}
Let $f\colon (Y,E)\to (X\ni P)$ be a toric canonical blow-up of three-dimensional toric terminal singularity, $f(E)=C$ and $\dim C=1$. Then we have the following statements.
	
$1)$ Let $(X\ni P)$ be a $\QQ$-factorial singularity, that is, it is $(\CC^3_{x_1x_2x_3}\ni 0)/\ZZ_r(-1,-q,1)$, where $\gcd(r,q)=1$, $0<q\leq r-1$ and $r\geq 1$.
Determine the numbers $u$, $v$ by the equality $uq+vr=1$, where $0\leq u\leq r-1$ and $u,v\in\ZZ$.
Consider the cone $\sigma$ defining $(X\ni P)$ $($see example $\ref{ex1}$ $1))$. Let $(w_1,w_2,w_3)$ be a primitive vector defining $f$. 

Then we have one of the two following cases up to permutation of coordinates: either $2A)$ $C=\{x_1=x_2=0\}/\ZZ_r$, $(w_1,w_2,w_3)=(1,w_2,0)$, or $2B)$ $C=\{x_2=x_3=0\}/\ZZ_r$, $(w_1,w_2,w_3)=(0,w_2,1)$.
The variety $Y$ has the singularities
$\frac1r(-1,w_2-q,1)$, $\frac1{rw_2}(-1+uw_2,-uw_2,1)$ in Case $2A)$ and $\frac1r(-1,-w_2-q,1)$, $\frac1{rw_2}(uw_2,-uw_2-1,1)$ in Case $2B)$.

Converse is also true: every such numbers $(w_1,w_2,w_3)$ define a canonical blow-up. 

A general element of the linear system $|-K_Y|$ has Du Val singularities.
	
Let $Q$ be a central point of second singularity in each of the two cases. Then $Q\in \mt{CS}(Y)$ if and only if $r\geq 2$. Therefore $f$ is a terminal blow-up if and only if it is the blow-up of the ideal of the curve $C$ \cite{Kaw1}.
	
$2)$ Let $(X\ni P)$ be a non-$\QQ$-factorial singularity, that is, $(X\ni P)\cong (\{x_1x_2+x_3x_4=0\}\subset (\CC^4_{x_1x_2x_3x_4},0))$. Then $C=\{x_1=x_2=x_3=0\}$ up to permutation of coordinates,  $f$ is induced by the blow-up of $(\CC^4,0)$ with weights $(w_1,w_2,w_1+w_2,0)$, where $w_1=1$, $w_2>0$ or $w_1>0, w_2=1$.
Converse is also true: every such numbers induce a canonical blow-up. A general element of the linear system $|-K_Y|$ has Du Val singularities.
	
The morphism $f$ is a terminal blow-up if and only if
$(w_1,w_2,w_3,w_4)=(1,1,2,0)$.
\begin{proof}
Let us prove 1).
Put $e'_1=e_1$, $e'_2=e_2$ and $e'_3=e_1+qe_2+re_3$ (see Example $\ref{ex1}$ $1))$. Then $w=w_ie'_i+w_je'_j$ for some $i<j$ and $w_i$, $w_j\in\ZZ_{\geq 1}$. We have $Y=T_N(\D)$ and 
$$
\D=\{\langle e'_k,e'_i,w\rangle, \langle e'_k,e'_j,w\rangle, \text{their faces}\},
$$
where $k$ is a third index other than the indices $i$ and $j$.
Consider an induced blow-up of general hyperplane section passing through the general point of $C$. Then $w_1=1$ or $w_2=1$. Now the statement is proved by a simple enumeration of the indices $i$ and $j$. As an example, consider $i=1$, $j=2$.
There are the two possibilities of weights: $(w_1,1,0)$ and $(1,w_2,0)$.
Let $(w_1,1,0)$. The variety $Y$ is covered by two affine charts with singularities of types $\frac1r(-q,qw_1-1,1)$ and $\frac1{rw_1}(-w_1,qw_1-1,1)$.
By Theorem \ref{cancyclic3dim} applied to the second singularity it follows that either $q=1$, or $w_1=1$, or $r=1$.
All these variants are realized, it is Case 2A). The possibility $(1,w_2,0)$ is considered similarly.
		
The proper transform of $\{x_2=0\}/\ZZ_r(-1,-q,1)$ is Du Val element of $|-K_Y|$.
		
The statement $Q\in \mt{CS}(Y)$ is obvious if we consider a blow-up with the weights $(-1+uw_2, (r-u)w_2, 1)$ in Case $2A)$ 
and $(uw_2,(r-u)w_2-1,1)$ in Case $2B)$ provided that $r\geq 2$.
		
Statement 2) obviously follows from Example \ref{ex1} 2). The proper transform of $\{x_1^{w_2}+x_2=0\}|_X$ ($\{x_1+x_2^{w_1}=0\}|_X$) is Du Val element of $|-K_Y|$ for the first (second) possibility.
\end{proof}
\end{proposition}

\begin{proposition}\label{cantoric}
Let $f\colon (Y,E)\to (X\ni P)$ be a toric canonical blow-up of three-dimensional toric terminal point, where $f(E)=P$.
Then we have the following statements.

$1)$ Let $(X\ni P)$ be a smooth point. Then $f$ is a weighted blow-up with weights $(w_1,w_2,1)$, $(l,l-1,2)$, 
$(15,10,6)$, $(12,8,5)$, $(10,7,4)$, $(9,6,4)$, $(8,5,3)$, $(7,5,3)$, $(6,4,3)$, $(5,3,2)$ or $(9,5,2)$ in some coordinate system, where $l \geq 3$.
Converse is also true: every such weights define a canonical blow-up. 
In all cases, except case $(9,5,2)$, a general element of the linear system $|-K_Y|$ has Du Val singularities.
In case $(9,5,2)$ we have 
\begin{gather*}
\min\{m|\exists D\in|-mK_Y|\ \text{such that}\ \\(Y,(1/m)D)\ \text{has canonical singularities}\}=3.
\end{gather*}

The morphism $f$ is a terminal blow-up if and only if it is a weighted blow-up
with weights $(w_1,w_2,1)$ in some coordinate system, where $\gcd(w_1,w_2)=1$.

$2)$ Let $(X\ni P)$ be a $\QQ$-factorial singularity of an index $\ge 2$, that is, it is of type $\frac1r(-1,-q,1)$, where $\gcd(r,q)=1$, $0<q\leq r-1$ and $r\geq 2$.
Let us consider the cone $\sigma$ defining the singularity $(X\ni P)$ $($see Example $\ref{ex1}$ $1))$. 
Determine the numbers $u$, $v$ by the equality $uq+vr=1$, where $0\leq u\leq r-1$ and $u,v\in\ZZ$.
Let $(w_1,w_2,w_3)$ be a primitive vector defining $f$.

Then we have one of the two following cases: either $2A)$ $(w_1,w_2,w_3)=(1,w_2,w_3)$, $w_3\le\min(r-1,\frac{rw_2-1}q)$ 
up to permutation of the numbers $w_1$ and $w_2$ provided that $q=1$,
or $2B)$ 
$(w_1,w_2,w_3)=(w_1,w_2,w_1+w_2-1)$, $w_1\ge 2$, $w_2\ge 2$, $0\le w_1(r-1)-w_2\le r-2$, $q=r-1$.
Converse is also true: every such numbers $(w_1,w_2,w_3)$ define a canonical blow-up. 
A general element of the linear system $|-K_Y|$ has Du Val singularities.

The morphism $f$ is a terminal blow-up if and only if it is a weighted blow-up
with weights $(u,1,r-u)$ \cite{Kaw1}.

$3)$ Let $(X\ni P)$ be a non-$\QQ$-factorial singularity, that is, $(X\ni P)\cong (\{x_1x_2+x_3x_4=0\}\subset (\CC^4_{x_1x_2x_3x_4},0))$.
Then $f$ is induced by the weighted blow-up of $(\CC^4,0)$ with weights $(w_1,w_2,w_3,w_4)$ up to analytical isomorphism of $(\CC^4,0)$, where $1+w_2=w_3+w_4$,
$w_1=1$.
Converse is also true: every such weights induce a canonical blow-up.
A general element of the linear system $|-K_Y|$ has Du Val singularities.

The morphism $f$ is a terminal blow-up if and only if
$(w_1,w_2,w_3,w_4)=(1,1,1,1)$ \cite{Corti}.
\begin{proof} Let us prove 1). Now we classify canonical blow-ups.
To be definite, assume that
$w_1\geq w_2\geq w_3$, where $(w_1,w_2,w_3)$ are primitive weights of $f$.
By $P_1$, $P_2$ and $P_3$ denote the zero-dimensional orbits (points) of $Y$.
These points have types $\frac1{w_1}(w_2,w_3,w_1-1)$, $\frac1{w_2}(w_1,w_3,w_2-1)$ and $\frac1{w_3}(w_1,w_2,w_3-1)$ respectively.

Assume that Cases 1) and 1) of Theorem \ref{cancyclic3dim} are satisfied at the points $P_1$ and $P_2$ respectively.
Then $w_1=w_2+w_3-1$ and $w_2|(2w_3-2)$. Thus we obtain the weights $(l,l,1)$, where $l\geq 1$ and $(3w_3-3,2w_3-2,w_3)$, where $w_3\geq 2$.
For the second possibility, the singularity is of type $\frac1{w_3}(3,2,1)$ at the point $P_3$, therefore $w_3\leq 6$, and it is easy to prove that
every value $w_3=2,\ldots,6$ is realized.

Assume that Cases 1) and 2) of Theorem \ref{cancyclic3dim} are satisfied at the points $P_1$ and $P_2$ respectively.
As above we obtain $w_1=w_2+w_3-1$ and have one of the following possibilities: i1) $w_3=1$, $w_3=2$ or i2) $2w_3-1=w_2$, $w_2=1,\ldots,4$.
These possibilities are realized.

Assume that Cases 1) and 3) of Theorem \ref{cancyclic3dim} are satisfied at the points $P_1$ and $P_2$ respectively.
Then $w_1=w_2+w_3-1$. Let the singularity be of type $\frac19(1,4,7)=\frac19(5,2,8)$ at the point $P_2$, in particular, $w_2=9$. Hence $w_3=2$ or $w_3=5$. 
It follows easily that these possibilities are not realized.
Let the singularity be of type $\frac1{14}(1,9,11)=\frac1{14}(5,3,13)$ at the point $P_2$, in particular, $w_2=14$.
Hence $w_3=3$ or $w_3=5$. It follows easily that these possibilities are not realized.

Assume that Cases 2) and 1) of Theorem \ref{cancyclic3dim} are satisfied at the points $P_1$ and $P_2$ respectively.
Then we obtain the two possibilities: i) $w_1=w_2+w_3$, $w_2=2w_3-1$, $w_3=2,3$ or ii) $w_3=1$.
These possibilities are realized.

Assume that Cases 2) and 2) of Theorem \ref{cancyclic3dim} are satisfied at the points $P_1$ and $P_2$ respectively.
As above it is easy to prove that new weights do not appear.

Assume that Cases 2) and 3) of Theorem \ref{cancyclic3dim} are satisfied at the points $P_1$ and $P_2$ respectively.
As above it is easy to prove that this case is not realized.

Assume that Cases 3) of Theorem \ref{cancyclic3dim} are satisfied at the point $P_1$.
Then $(w_1,w_2,w_3)=(9,5,2)$ or $(14,5,3)$. 
It is obvious that only the first possibility is realized.

For any weights obtained, except case
$(9,5,2)$, 
we can easily find a surface 
$S\subset X$ with Du Val singularity at the point $P$ such that
$a(S,E)=0$. For example, 
the surface $S$ is given (locally at the point $P$) by the equations
$x_1x_2+x_3^{w_1+w_2}=0$ and $x_1^2+x_2^3+x_2x_3^3=0$ for cases $(w_1,w_2,1)$ and $(5,3,2)$
respectively.
Therefore $S_Y\in |-K_Y|$ has Du Val singularities. 

In case $(9,5,2)$ the variety $Y$ has the two non-terminal isolated singularities at the points $P_1$ and $P_2$ ($\mt{CS}(Y)=\{P_1, P_2\}$).
Let $C\subset E=\PP(9,5,2)$ be a curve not passing through the points $P_1$ and $P_2$. Then a 
(quasihomogeneous) degree of $C$ is at least 45 since it must be divided by 9 and 5. Hence $m\geq 3$, and the required element $D$ is the proper transform of 
$x_1^5+x_2^9+x_3^{23}=0$.
The other statements of 1) are obvious.

Let us prove 2). Now we classify canonical blow-ups.
The variety $Y$ is covered by three affine charts with singularities of types
$\frac1{w_3}(-w_1,-w_2,1)$, $\frac1{rw_2-qw_3}(-w_1+uw_2+vw_3,-uw_2-vw_3,1)$ and
$\frac1{rw_1-w_3}(-w_1,qw_1-w_2,1)$ respectively. 
The corresponding zero-dimensional orbits of $Y$ are denoted 
by $P_1$, $P_2$ and $P_3$.
Note that $rw_1-w_3, rw_2-qw_3\in\ZZ_{\geq 1}$. Obviously, $a(S,0)=\frac1r(w_3+rw_2-qw_3+rw_1-w_3)-1$.
The minimal discrepancy of $(X\ni P)$ is equal to $\frac1r$. If $a(S,0)=\frac1r$, that it is easy to calculate that $f$ is a terminal blow-up, that
is, a weighted blow-up with weights $(u,1,r-u)$ \cite{Kaw1}.
Therefore we suppose that $a(S,0)>\frac1r$.

Since $Y$ has canonical singularities, then for some $j\in\{1,2,3\}$ we have the inequality $\frac1r\geq a(S,0)/N_j$ 
and one of the two following
requirements: either $P_j\in\mt{CS}(Z)$, or the singularity at the point $P_j$ is of type $\frac1{N_j}(1,-1,0)$, 
where $N_j\ge 2$, $N_1=w_3$, $N_2=rw_1-w_3$, $N_3=rw_2-qw_3$. 
This is called {\it Property} $R_j$. Note if $j=3$ then $w_1=1$. Therefore we suppose that $j\le 2$.

Let $w_1=\max\{w_1,w_2,w_3\}$. 
Assume that Case 1) of Theorem \ref{cancyclic3dim} is satisfied at the point $P_2$.
Then $q=1$ and $w_2=1$.
Assume that Case 2) of Theorem \ref{cancyclic3dim} is satisfied at the point $P_2$.
Then, either $w_1=w_2=w_3=1$, or $q=2$, $w_1=w_2$, $w_1\ge 2$, $r\ge 3$. 
Since the inequality of Property $R_2$ holds then the second possibility is not realized. 
It is not hard to prove that Case 3) of Theorem \ref{cancyclic3dim} is not realized at the point $P_2$.

Let $w_2=\max\{w_1,w_2,w_3\}$. Property $R_1$ is not realized. Therefore Property $R_2$ holds. Then
$w_2=w_3$, and we have $w_1=1$ by Theorem \ref{cancyclic3dim} for the point $P_1$.

Let us consider the last case $w_3>\max\{w_1,w_2\}$. The possibility $w_1=1$ holds. Therefore we suppose that $w_1\ge 2$.
If $w_2=1$ then Theorem \ref{cancyclic3dim} for the point $P_2$ implies $q=1$. Therefore we suppose that $w_2\ge 2$.

Assume that Case 1) of Theorem \ref{cancyclic3dim} is satisfied at the point $P_1$.
Then $w_1+w_2-1=w_3$. If the inequality of Property $R_1$ holds then $q=r-1$.
Therefore we suppose that Property $R_2$ holds and $N_2>w_3$. 
It is not hard to prove that Case 3) of Theorem \ref{cancyclic3dim} is not realized at the point $P_2$.
If Case 1) of Theorem \ref{cancyclic3dim} is satisfied at the point $P_2$ then the inequality of Property $R_2$ implies that
$(q-1)w_1-w_2+1=0$, but this equality contradicts the same inequality. Therefore the singularity is of type $\frac1{N_2}(1,-1,0)$ at the point
$P_2$. Therefore $w_1=1$. We obtain the contradiction.

Assume that Case 2) of Theorem \ref{cancyclic3dim} is satisfied at the point $P_1$.
Then $w_1+w_2=w_3$ and Property $R_2$ holds.
Let Case 3) of Theorem \ref{cancyclic3dim} be satisfied at the point $P_2$. Then it is not hard to prove that
$(w_1,w_2,w_3,r)=(2,2q+5,2q+7,q+8)$.
We obtain a contradiction with Theorem \ref{cancyclic3dim} for the point $P_3$ since $0<uw_2+vw_3\le N_3$. 
Let Case 1) of Theorem \ref{cancyclic3dim} be satisfied at the point $P_2$.
The inequality of Property $R_2$ implies that
$(q-1)w_1-w_2+1=0$, but this equality contradicts the same inequality. Therefore the singularity is of type $\frac1{N_2}(1,-1,0)$ at the point
$P_2$.
Considering two possibilities: $N_2\le w_1$ and
$N_2>w_1$, it is easy to obtain a contradiction.

Now, applying the blow-up classification obtained, we can prove that the proper transform of the divisor 
$$S_k=\{x_k=0\}/\ZZ_r \subset (\CC^3_{x_1x_2x_3},0)/\ZZ_r(-1,-q,1)$$
is Du Val element of $|-K_Y|$ for some $k$. 
The other statements of 2) are obvious.

Let us prove 3). Consider Example \ref{ex1} 2). Now we classify canonical blow-ups.
Obviously, $a(S,0)=w_1+w_2-1=w_3+w_4-1$.
The variety $Y$ is covered by three affine charts with singularities of types 
$\frac1{w_1}(w_3,w_4,-1)$, $\frac1{w_2}(w_3,w_4,-1)$, $\frac1{w_3}(w_1,w_2,-1)$ and
$\frac1{w_4}(w_1,w_2,-1)$ respectively.
The minimal discrepancy of $(X\ni P)$ is equal to 1.
If $a(S,0)=1$ then it is easy to calculate that $f$ is a terminal blow-up induced
by the weighted blow-up with weights $(1,1,1,1)$ \cite{Corti}. Therefore we suppose that
$a(S,0)>1$.
Since $Y$ has canonical singularities then $1\geq a(S,0)/w_j$ for some $j$.
Hence $w_i=1$ for some $i\ne j$ such that $w_i+w_j-1=a(S,0)$. 
The proper transform of
$\{x_i^{w_j}+x_j=0\}|_X$ is
Du Val element of $|-K_Y|$. 
The other statements of 3) are obvious.
\end{proof}
\end{proposition}

\begin{definition}\label{nondegen}
Let $(X\ni P)$ be an $n$-dimensional $\QQ$-factorial toric singularity. Then
$(X\ni P)\cong (\CC^n\ni 0)/G$, where $G$ is an abelian group acting freely in codimension 1. 
The singularity $(\CC^n\ni 0)/G$ is given by the simplicial cone $\sigma_G$ in the lattice $N=\ZZ^n$.
\par
Let a power series (polynomial) $\varphi=\sum_{m}a_{m}x^m \in \CC[[x_{1},x_{2},\ldots,x_{n}]]$ be $G$-semiinvariant. 

The
{\it Newton polyhedron} $\Gamma_{+}(\varphi)$ in
$\RR^n$ is the convex hull of the set

$$ \bigcup_{x^m\in \varphi}
(m+\sigma_G^{\vee}),\ \text{where}\ \sigma_G^{\vee}\ \text{is a dual cone in}\ M_{\RR}.
$$

For any face
$\gamma$ of $\Gamma_{+}(\varphi)$ we define
\[ \varphi_{\gamma}=\sum_{m\in\gamma}a_{m}x^m. \]

The function $\varphi$ is said to be {\it non-degenerate} if,
for any
compact face
$\gamma$ of the Newton polyhedron, the polynomial equation $\varphi_{\gamma}=0$ defines
a smooth hypersurface in the complement of the set
$x_1x_2\ldots x_{n}=0$. The effective Weil divisor $D$ on $X$ is said to be {\it non-degenerate}
if the $G$-semiinvariant polynomial $\varphi$ defining $D$ in $\CC^n$ is non-degenerate.
\par
For any effective Weil divisor $D$ there exists the fan
$\Delta$ depending on Newton polyhedron
$\Gamma_{+}(\varphi)$ such that
$T_N(\Delta)$ is a smooth variety and a toric birational morphism $\psi \colon T_N(\Delta) \to \CC^n$ is a resolution of
non-degenerate singularities of $D$. So, $\psi$ is said a {\it partial resolution of $(X,D)$}. In particular, if $D$ is a non-degenerate boundary then
$\psi$ is a toric log resolution of the pair $(X,D)$. If $(X\ni P)$ is a smooth variety then this statement was proved in the paper \cite{Var}.
Note that the proof from the paper \cite{Var} is rewritten immediately in our case if we will use  
our Newton polyhedron instead of standard Newton polyhedron.
\end{definition}

The next Theorems \ref{toricplt} and \ref{toriccan} are criteria 
of the characterization
of toric plt and canonical blow-up respectively.
They explicitly show a nature of non-toric contractions.

\begin{theorem}\label{toricplt}
Let $f\colon (Y,E) \to (X \ni P)$ be a plt blow-up of $\QQ$-factorial toric singularity, and let
$f(E)$ be a toric subvariety. Then $f$ is a toric morphism $($under a suitable identification$)$ if and only if there exists an effective non-degenerate Weil divisor
$D$ on $(X\ni P)$ and a number $d>0$ with the following properties:
\begin{enumerate}
\item[1)] $a(E,dD)=-1;$
\item[2)] $E$ is a unique exceptional divisor of $(X,dD)$ with discrepancy $\leq -1$ and $\down{dD}=0$.
\end{enumerate}
\begin{proof}
First let us prove the necessary condition.
Let $D_Y\in |-n(K_Y+E)|$ be a general element for $n\gg 0$. Put $D=f(D_Y)$ and $d=\frac{1}{n}$.
Then $K_Y+E+dD_Y=f^*(K_X+dD)$ is a plt divisor.
Since $D_Y$ is a general divisor by construction, then $D$ is an irreducible reduced non-degenerate divisor. 
\par
Finally let us prove the sufficient condition. Consider the toric log resolution
$\psi\colon Z\to X$ of $(X,dD)$. Write
$$
K_Z+dD_Z+\sum a_iE_i=\psi^*\Big(K_X+dD\Big).
$$ 
By theorem assertion $(Z,dD_Z+\sum a_iE_i)$ is a plt pair. Therefore $E\subset\Exc \psi$. 

Considering corresponding fans (see \cite{Reid})
we have the composition of toric log flips $Z\dashrightarrow Z'$ over $(X\ni P)$ such that the (induced) toric divisorial contraction
$\psi'\colon Z'\to (X\ni P)$ is isomorphic to
$\psi'_1\circ\psi'_2$, where $\psi'_1$, $\psi'_2$ are toric divisorial contractions and $E=\Exc \psi'_1$.  
Therefore $f$ and $\psi'_1$ are isomorphic by Remark \ref{pltpropert} 5).
\end{proof}
\end{theorem}

\begin{theorem}\label{toriccan}
Let $f\colon (Y,E)\to (X\ni P)$ be a canonical blow-up of $\QQ$-factorial toric singularity, and let
$f(E)$ be a toric subvariety.
Then $f$ is a toric morphism $($under a suitable identification$)$ if and only if there exists an effective non-degenerate Weil divisor
$D$ on $(X\ni P)$  and a number $d>0$ with the following properties:
\begin{enumerate}
\item[1)] $a(E,dD)=0;$
\item[2)] $(X,dD)$ has canonical singularities and $\down{2dD}=0$.
\end{enumerate}
\begin{proof}
First let us prove the necessary condition. Let $D_Y\in |-nK_Y|$ be a general element for $n\gg 0$. Put $D=f(D_Y)$ and $d=\frac{1}{n}$.
Then the divisor $K_Y+dD_Y=f^*(K_X+dD)$ has canonical singularities.
Since $D_Y$ is a general divisor by construction, then $D$ is an irreducible reduced non-degenerate divisor.
\par
Finally let us prove the sufficient condition. Consider the toric log resolution
$\psi\colon Z\to X$ of $(X,dD)$. Write
$$
K_Z+dD_Z+\sum a_iE_i=\psi^*\Big(K_X+dD\Big).
$$ 
By theorem assertion $(Z,dD_Z+\sum a_iE_i)$ is a terminal pair. Therefore $E\subset\Exc \psi$. 
Considering corresponding fans (see \cite{Reid})
we have the composition of toric log flips $Z\dashrightarrow Z'$ over $(X\ni P)$ such that the (induced) toric divisorial contraction
$\psi'\colon Z'\to (X\ni P)$ is isomorphic to
$\psi'_1\circ\psi'_2$, where $\psi'_1$, $\psi'_2$ are toric divisorial contractions and $E=\Exc \psi'_1$.  
Therefore $f$ and $\psi'_1$ are isomorphic by Proposition \ref{lemma1}.
\end{proof}
\end{theorem}

\begin{definition}
The subvariety $Y$ is said to be a {\it non-toric subvariety} of the toric pair $(X,D)$, if there is not any toric structure of $X$ such that
$(X,D)$ is a toric pair and $Y$ is a toric subvariety.
\end{definition}

\begin{example}
Consider the toric variety $X=\PP_{x_1x_2x_3}(1,2,3)$. 

1) Let $D=0$. The point $P$ is a non-toric subvariety of $(X,D)$ if and only if $P=(0:1:a)$, where $a\ne 0$. The irreducible curve $C$
is a non-toric subvariety of $(X,D)$ if and only if $C\ne\{x_1=0\}$, $C\ne\{x_2+ax_1^2=0\}$ and $C\ne\{x_3+ax_2x_1+bx_1^3=0\}$.

2) Let $D=\{x_1=0\}+\{x_2=0\}$. The point $P$ is a non-toric subvariety of $(X,D)$ if and only if $P=(0:1:a)$, where $a\ne 0$. The irreducible curve $C$
is a non-toric subvariety of $(X,D)$ if and only if $C\ne\{x_1=0\}$, $C\ne\{x_2=0\}$ and $C\ne\{x_3+ax_2x_1+bx_1^3=0\}$.

3) Let $D=\{x_1=0\}+\{x_2=0\}+\{x_3=0\}$. The point $P$ is a non-toric subvariety of $(X,D)$ if and only if $P\ne(1:0:0)$, $P\ne(0:1:0)$ and
$P\ne(0:0:1)$. 
The irreducible curve $C$ is a non-toric subvariety of $(X,D)$ if and only if $C\ne\{x_1=0\}$, $C\ne\{x_2=0\}$ and $C\ne\{x_3=0\}$.
\end{example}

Next Theorems \ref{dim2plt} and \ref{dim2can} are two-dimensional analogs of main theorems. Their proofs clearly describe 
the main method used in this paper.

\begin{theorem}\label{dim2plt}\cite{PrLect}
Let $f\colon (Y,E) \to (X \ni P)$ be a plt blow-up of two-dimensional toric singularity.
Then $f$ is a toric morphism $($under a suitable identification$)$.
\begin{proof}
A two-dimensional toric singularity is always
$\QQ$-factorial.
Let $f$ be a non-toric morphism (up to identification). 
Let $D_Y\in |-n(K_Y+E)|$ is a general element of $n\gg 0$. Put $D_X=f(D_Y)$ and $d=\frac{1}{n}$.
Then $(X,dD_X)$ is a log canonical pair, $a(E,dD_X)=-1$ and $E$ is a unique exceptional divisor with discrepancy $-1$.

By Criterion \ref{toricplt} there exists a toric divisorial contraction $g\colon Z\to X$ with the following properties.
\begin{enumerate}
\item[A)] The exceptional set $\Exc g=S$ is an irreducible divisor ($S\cong\PP^1$), 
the divisors $S$ and $E$ define the different discrete valuations of the function field
$\mathcal K(X)$.
\item[B)] By $\G$ denote the center of $E$ on $S$. Then the point $\G$ is a non-toric subvariety of $Z$
for any toric structure of $(X \ni P)$.
In the other words, $\G$ is a non-toric subvariety of the toric pair 
$(S,\Diff_S(0))$.
\end{enumerate}

Condition B) implies that the surface $Z$ has the two singular points
$P_1$ and $P_2$, which lie on the curve $S$. Also $\G$ is a non-toric point of 
$(S,\Diff_S(0))\cong (\PP^1, \frac{n_1-1}{n_1}P_1+\frac{n_2-1}{n_2}P_2)$, where $n_1\geq 2, n_2 \geq 2$. Write

$$
K_Z+dD_Z+aS=g^*\big(K_X+dD_X\big),
$$
where $a<1$. Hence
$$
a\big(E,S+dD_Z\big)< a\big(E,aS+dD_Z\big)=-1.
$$
Therefore $K_Z+S+dD_Z$ is not a log canonical divisor at the point $\G$ and is an anti-ample over $X$ divisor. Hence, by the inversion of adjunction,
$K_S+\Diff_S(dD_Z)$ is not a log canonical divisor at the point $\G$ and is an anti-ample divisor. We obtain the contradiction
\begin{gather*}
0>\deg\big(K_S+\Diff_S(dD_Z)\big)>-2+\frac{n_1-1}{n_1}+\frac{n_2-1}{n_2}+1\geq 0.
\end{gather*}
 
\end{proof}
\end{theorem}

\begin{theorem}\label{dim2can}\cite{Mor}
Let $f\colon (Y,E)\to (X\ni P)$ be a canonical blow-up of two-dimensional toric singularity.
Then $(X\ni P)$ is a smooth point, and $f$ is a weighted blow-up with weights $(1,\alpha)$ $($under a suitable identification$)$.
\begin{proof}
Theorem assertion implies that $(X\ni P)$ is a terminal point, therefore it is smooth.

Assume that $f$ is a toric morphism then $f$ is a weighted blow-up of the smooth point with weights
$(\beta,\alpha)$. Since $Y$ is Du Val surface then
$\alpha=1$ or $\beta=1$.

Let $f$ be a non-toric morphism (up to identification). 
Let $D_Y\in |-nK_Y|$ be a general element for $n\gg 0$. Put $D_X=f(D_Y)$ and $d=\frac{1}{n}$.
The pair $(X,dD_X)$ has canonical singularities and $a(E,dD_X)=0$.

By Criterion \ref{toriccan} there exists a toric divisorial contraction $g\colon Z\to X$ with the following properties.
\begin{enumerate}
\item[A)] The exceptional set $\Exc g=S$ is an irreducible divisor ($S\cong\PP^1$), 
the divisors $S$ and $E$ define the different discrete valuations of the function field
$\mathcal K(X)$.
\item[B)] By $\G$ denote the center of $E$ on $S$. Then the point $\G$ is a non-toric subvariety of $Z$
for any toric structure of $(X \ni P)$.
In the other words, $\G$ is a non-toric subvariety of the toric pair 
$(S,\Diff_S(0))$.
\end{enumerate}

Condition B) implies that the surface $Z$ has the two singular points
$P_1$ and $P_2$, which lie on the curve $S$. Also $\G$ is a non-toric point of 
$(S,\Diff_S(0))\cong (\PP^1, \frac{n_1-1}{n_1}P_1+\frac{n_2-1}{n_2}P_2)$, where $n_1\geq 2, n_2 \geq 2$. Write

$$
K_Z+dD_Z+S=g^*\big(K_X+dD_X\big)+(a(S,dD_X)+1)S,
$$
where $a(S,dD_X)\geq 0$.
Since $S$ is (locally) Cartier divisor at the point $\G$, then
$$
a\big(E,S+dD_Z\big)\leq a\big(E,dD_X\big)-1=-1.
$$

Therefore $K_Z+S+dD_Z$ is not a plt divisor at the point $\G$ and is an anti-ample divisor over $X$. 
Hence, by the inversion of adjunction
$K_S+\Diff_S(dD_Z)$ is not a klt divisor at the point $\G$ and is an anti-ample divisor. We obtain the contradiction

\begin{gather*}
0>\deg\big(K_S+\Diff_S(dD_Z)\big)\geq -2+\frac{n_1-1}{n_1}+\frac{n_2-1}{n_2}+1\geq 0.
\end{gather*}

\end{proof}
\end{theorem}

\begin{example}
Theorems \ref{dim2plt} and \ref{dim2can} cannot be generalized in dimension at least three for divisorial contraction to a point.
Consider the blow-up $g\colon Z\to (X\ni P)$ with the weights $(1,\ldots,1)$, where $(X\ni P)\cong(\CC^n_{x_1\ldots x_n}\ni 0)$ and
consider the divisors $D=\{x_1^2+\ldots+x_n^2=0\}$, $T^i=\{x_i=0\}$, where $i=1,\ldots ,n$ and
$n\ge 3$. The exceptional set $\Exc g=S$ is isomorphic to $\PP^{n-1}$, $Q=S\cap D_Z$ is a smooth quadric. 
Let
$\widetilde g\colon \widetilde Z\to Z$ be the standard blow-up of the ideal $I_Q$. By the base point free theorem \cite{KMM} the linear system $|mD_{\widetilde Z}|$ 
gives a divisorial contraction
$\varphi\colon\widetilde Z\to Y$, which contracts the divisor $S_{\widetilde Z}\cong \PP^{n-1}$ for $m\gg 0$. Since the divisor
$K_{\widetilde Z}+S_{\widetilde Z}+\sum_{i=1}^nT^i_{\widetilde Z}\sim 0/Y$ has log canonical singularities, then by
Shokurov's criterion on the characterization
of toric varieties
for divisorial contractions to a $\QQ$-factorial singularity \cite[Theorem 18.22]{Koetal},
the morphism $\varphi$ is toric. Hence $Y$ has only one singularity and its type is $\frac1r(1,\ldots,1)$.
Let $l$ be a straight line in a general position in $S_{\widetilde Z}$. Considering $\varphi$ we have $S_{\widetilde Z}\cdot l=-r$, and considering $g\circ \widetilde g$ 
we have
$S_{\widetilde Z}\cdot l=-3$, hence $r=3$. 
\par
We obtain a non-toric divisorial contraction $f\colon Y\to (X\ni P)$. The variety $Y$ has only one singularity and its type is
$\frac13(1,\ldots,1)$.
Thus, if $n\ge 4$, then $Y$ is a terminal variety, and if
$n=3$, then $Y$ is a canonical non-terminal variety (cf. \cite{Kawakita1}). The blow-up $f$ is plt
since the exceptional set $\Exc f$ is a cone over a smooth $(n-2)$-dimensional quadric.
\end{example}

We will apply the following special case of Shokurov's criterion on the characterization
of toric varieties.
\begin{proposition}\label{toriccrit}
Let $f\colon (X,D)\to (Z\ni P)$ be a small contraction of the
$\QQ$-factorial threefold $X$. 
Assume that $D=\sum_{i=1}^rD_i$, where $D_i$ is a prime divisor for each $i$.
Assume that $K_X+D$ is a log canonical divisor,
$-(K_X+D)$ is a $f$-nef divisor and $\Exc f=C$ is an irreducible curve $(\rho(X/Z)=1)$.
Then $r\leq 4$. 
Moreover, the equality holds if and only
if the pair $(X/Z\ni P,D)$ is analytically isomorphic to a toric
pair, in particular, $K_X+D\sim 0/Z$.
\begin{proof}
If the pair $(X/Z\ni P,D)$ is analytically isomorphic to a toric
pair then all statements immediately follow from the description of toric log flips
\cite{Reid}.
Let $r\geq 4$. 
Let the divisor $K_X+D'$ be a $\QQ$-complement
of $K_X+D$.
It exists, since we can add to the divisor $D$ the necessary
number of general hyperplane sections of $X$.
So, by abundance theorem \cite[Theorem 8.4]{Koetal} the
$\QQ$-complement $D'$ required is constructed for our contraction $(X/Z\ni P,D)$.

Put $D'=\sum d_iD_i'$. We will prove that $D'=D$. For any $\QQ$-Weil divisor $B=\sum b_iB_i$ we define
$||B||=\sum b_i$. Put

$$
D^{\mt{hor}}=\sum_{i\colon D'_i\cdot C>0} d_iD'_i \ \ \ \text{and}\ \ \
D^{\mt{vert}}=\sum_{i\colon D'_i\cdot C\le 0} d_iD'_i.
$$

Let $f^+\colon X^+\to Z$ be a log flip of $f$ and
$\Exc f^+=C^+$.

\begin{lemma}\cite[Lemma 2.10]{PrCbf}
We have $||D^{\mt{hor}}||=||D^{\mt{vert}}||=2$. Hence,
$D=D'$. Moreover, $C\not \subset \Supp D^{\mt{hor}}$,
$C^+\not\subset \Supp (D^{\mt{vert}})^+$ and $D'_i\cdot C\ne 0$ for all $i$.
\begin{proof}
Since $K_X+D$ is a log canonical divisor then
$||D^{\mt{vert}}||\le 2$.
Since $K_{X^+}+D^+$ is a log canonical divisor then
$||D^{\mt{hor}}||\le 2$. The statements remained are obvious.
\end{proof}
\end{lemma}

Let $S$ be an irreducible component of the divisor
$D^{\mt{vert}}$ and let
$F=D-S$. The divisorial log contraction
$(S,\Diff_S(F))\to (f(S)\ni P)$ is toric by the two-dimensional Shokurov's criterion on the characterization
of toric varieties \cite[Theorem 6.4]{Sh2}.
In particular, it is a toric blow-up of cyclic singularity.
Thus, the singularities of
$X$ are toric by three-dimensional Shokurov's criterion on the characterization
of toric varieties for $\QQ$-factorial singularities \cite[Theorem 18.22]{Koetal}.
Replacing $X$ by $X^+$ it can be assumed that $-(K_X+S)$ is a $f$-ample divisor and
$S\cdot C<0$. 

In order to prove the proposition we will apply some modification, which is a toric one by its nature. After it
we will get some small contraction, which is analytically isomorphic to a small toric contraction of Example \ref{ex1} 2). 
Therefore the initial contraction is a toric up to analytical isomorphism.

Now, taking toric blow-ups of $X$ (every time we take an one blow-up with a unique exceptional divisor that has a minimal
discrepancy of a singularity considered and consider 
two extremal rays on a variety obtained), it can be assumed that
$S$ is a smooth surface, and $X$ is a smooth variety outside the curves $C$.
The condition that $-(K_X+S)$ is $f$-ample holds is preserved, since the discrepancies of exceptional divisors of 
$(X,S)$ are less than and equal to 0.
In some analytical neighborhood of every point of $C$ the variety $X$ is analytically isomorphic to $\frac1{k}(q,1)\times \CC^1$, where $(k,q)=1$.

Assume that $k\ge 2$.
Consider a natural cyclic cover $\psi\colon \overline X\to X$ of degree $k$.
Put $\overline C=\psi^{-1}(C)$ and let $\overline Z$ be the normalization of $Z$ in the function
field of $\overline X$.
Let $\overline f\colon\overline X\to (\overline Z\ni \overline P)$ be the
induced small contraction of the curve
$\overline C$. Thus we can assume that $k=1$, that is, $X$ is a smooth variety.

Since 
$-K_S$ is a $f$-ample divisor then
$f\colon S\to f(S)$ is the contraction of the
(--1) curve $C$ and $(K_X+S)\cdot C=-1$.
We have $S\cdot C=-m$, $K_X\cdot C=m-1$ for some $m\in \ZZ_{\geq 1}$. 

Let $m\geq 2$.
Using the natural section of $\OO_X(S)$ we can construct
a degree $m$-cyclic cover
$\varphi\colon \widetilde X\to X$ ramified along $S$ (cf. \cite[Theorem 5.4]{Koetal}).
Let $\widetilde C=\varphi^{-1}(C)$ and let $\widetilde Z$ be the normalization of $Z$ in the function
field of $\widetilde X$.
Let $\widetilde f\colon\widetilde X\to (\widetilde Z\ni \widetilde P)$ be the
induced small contraction of the curve
$\widetilde C$. By the ramification formula

$$
K_{\widetilde X}\cdot\widetilde C=\varphi^*\Big(K_X+\frac{m-1}m S\Big)
\cdot \widetilde C=
K_X\cdot C+\frac{m-1}mS\cdot C=0.
$$

Thus we can assume that $f$ is a small flopping contraction with respect to $K_X$ ($K_X\cdot C=0$), that is, 
we can assume that $m=1$.

Since the minimal discrepancy of three-dimensional terminal non-cDV singularity is strict less than 1 then
$(Z\ni P)\cong (g=0\subset (\CC^4,0))$ is an isolated cDV (terminal) singularity.
Note that $(D_1+D_2)\cdot C=(D_3+D_4)\cdot C=0$
up to permutation of components of $D$. Hence $L_1$ and $L_2$ are Cartier divisors, where
$L_1=f(D_1)+f(D_2)$ and $L_2=f(D_3)+f(D_4)$.
By Bertini theorem \cite[Theorem 4.8]{Kollar}
the pair $(Z\ni P,H+L_i)$ is log canonical for any $i=1,2$, where
$H$ is a general hyperplane section passing through the point $P$. By the inversion of adjunction $(H\ni P,L_i|_H)$ is a log canonical pair.
Thus, the classification of two-dimensional log canonical pairs \cite{Koetal} implies that
$(H\ni P)$ is a cyclic singularity at the point
$P$, that is, it has type $\AAA_k$.
By the paper \cite{KaM} or the paper \cite{Kaw0} the singularity
$(H\ni P)$ is of type $\AAA_1$.
Thus 
$$
(Z\ni P)\cong (xy+z^2+t^{2l}=0\subset (\CC^4,0))
$$
and
$f(D)=\{x=0\}|_Z+\{y=0\}|_Z$. Since
$(Z\ni P, f(D))$ is a log canonical pair then we can take the weighted blow-up of $(\CC^4,0)$ with the weights $(l,l,l,1)$ and
obtain $l=1$. This completes the proof.
\end{proof}
\end{proposition}

\begin{remark}
Let $\rho(P)$ be a rank of local analytic group of Weil divisors at the point $P$.
Then the Proposition \ref{toriccrit} implies easily 
Shokurov's criterion on the characterization
of toric varieties for three-dimensional singularities $(Z\ni P)$ if $\rho(P)=1$, and hence the same criterion
for three-dimensional divisorial contractions $f\colon X\to (Z\ni P)$ if $\rho(P)=1$.
\end{remark}

\section{\bf Three-dimensional blow-ups. Case of curve} \label{chpt}

\begin{example} \label{pltdim1}
Now we construct the examples of three-dimensional non-toric plt blow-ups $f\colon (Y,E)\to (X\supset C\ni P)$ provided that $(X\ni P)$ is a
$\QQ$-gorenstein toric singularity,
$\dim f(E)=1$ and the curve $C=f(E)$ is a toric (smooth) subvariety. Depending on a type of $(X\ni P)$ we consider two Cases {\bf A1)} and {\bf A2)}.
\par {\bf A1)} Let $(X\ni P)$ be a $\QQ$-factorial toric singularity, that is, $(X\ni P)\cong (\CC^3\ni 0)/G$, 
where $G$ is an abelian group acting freely in codimension 1. 

All plt blow-ups are constructed by the procedure illustrated on the next diagram (Fig. 1) and defined below.

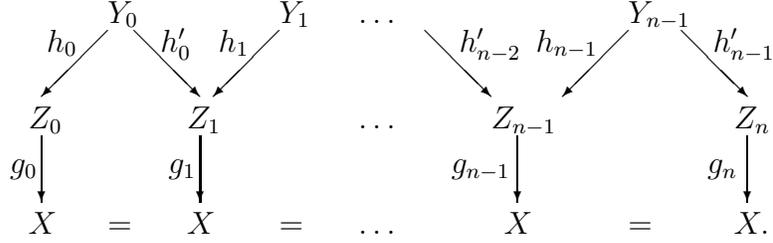
\[
\xymatrix{& Y_0 \ar[dl]_{h_0} \ar[d]^{h'_0} & Y_1 \ar[l]_{h_1} \ar@{-->}[d]^{h'_1} & \ar[l]_{h_2} Y_2 \ar@{-->}[d]^{h'_2}\\
Z_0 \ar[d]_{g_0}& Z_1 \ar[d]_{g_1} & Z_2 \ar[d]_{g_2} & Z_3	\ar[d]_{g_3}\\
X& X& X & X
}
\]\begin{center}
	Fig. $1$
\end{center}

{\it First step}.
Let $g_0\colon (Z_0,S_0)\to (X\supset C\ni P)$ be a toric blow-up, where $\Exc g_0=S_0$ is an irreducible divisor 
and $g_0(S_0)=C$.  
Recall that $g_0$ is a plt blow-up,
the surface $S_0$ is a toric conic bundle, $\rho(S_0/C)=1$ 
and $\Diff_{S_0}(0)=\frac{w^1_0-1}{w^1_0}E^1_0+\frac{w^2_0-1}{w^2_0}E^2_0+\frac{d_0-1}{d_0}F_0$, where
$E^1_0$, $E^2_0$ are some sections of conic bundle, $F_0$ is a fiber over $P$ and
$w^1_0, w^2_0, d_0\in\ZZ_{\geq 1}$. Let us remark that the numbers $w^1_0, w^2_0$ determine $g_0$.
Moreover, $d_0=1$ if $(X\ni P)$ is a smooth point.

Assume that there exists a curve $\G_0\subset S_0$ with the following two properties:
1) $K_{S_0}+\Diff_{S_0}(0)+\G_0$ is a plt and $g_0$-anti-ample divisor; 2) $\G_0$ is a non-toric subvariety in any  
analytical neighborhood of the fiber $F_0$
on the toric variety $Z_0$ for any toric structure of $(X\ni P)$, that is, the curve $\G_0$ is a non-toric subvariety of $(S_0,\Diff_{S_0}(0))$ in any analytical
neighborhood of $F_0$ on $S_0$.

By considering the general fiber over a general point of $C$ we obtain 
$w^i_0=1$ for some $i=1,2$.
To be definite, put $w^1_0=1$ and let
$Q_0=E^2_0\cap F_0$. Applying the adjunction formula it is easy to prove that $\G_0\cap F_0=Q_0$, 
$w^2_0\ge 2$, $d_0=1$, $(S_0\ni Q_0)$ is of type $\frac1{r_0}(1,1)$ ($r_0\ge 1$) and $\G_0\cdot F_0=\G\cdot E_0^2=\frac1{r_0}$.

\begin{remark}
Let $(X\ni P)$ be a terminal singularity, that is, $(X\ni P)\cong (\CC^3_{x_1, x_2, x_3}\ni 0)/\ZZ_r(-1, -q, 1)$. Then $r=r_0$ and one of the following cases holds by simple calculation.

1) $C=\{x_1=x_2=0\}$, $g_0$ is a blow-up with weights $(w_0^2, 1, 0)$, $r_0|w_0^2$ or $(1, w_0^2, 0)$, $r_0|(w_0^2-q+1)$.

2) $C=\{x_1=x_3=0\}$, $g_0$ is a blow-up with weights $(w_0^2, 0, 1)$, $r_0|(w_0^2+1+q)$ or $(1, 0, w_0^2)$, $r_0|(w_0^2-q+1)$.

3) $C=\{x_2=x_3=0\}$, $g_0$ is a blow-up with weights $(0, w_0^2, 1)$, $r_0|(w_0^2+1+q)$ or $(0, 1, w_0^2)$, $r_0|w_0^2$.
\end{remark}

Consider an arbitrary toric structure of $Z_0$ in any neighborhood of the point $Q_0$ such that 
$\G_0$ is also a toric subvariety of $Z_0$. 
Let $h_0\colon (Y_0,(S_1)_{Y_0})\to (Z_0\supset \G_0\ni Q_0)$ be an arbitrary toric blow-up of the curve 
$\G_0$ with an unique exceptional divisor ($\Exc h_0=(S_1)_{Y_0}$).
The structures of $h_0$ and $g_0$ are similar, in particular, 
$h_0$ is determined by some numbers $w^1_1$ and $w^2_1$, $(S_0)_{Y_0}\cong S_0$.

The set of all possible blow-ups $h_0$ for any toric structure of $(Z_0\ni Q_0, \G_0)$ is denoted by $\mathcal H_0$.

Let $(D_0)_{Z_0}$ be a toric Weil divisor of $(Z_0\ni Q_0)$ such that $(D_0)_{Z_0}|_{S_0}=\G_0$ and 
$a((S_1)_{Y_0},(D_0)_{Z_0}+S_0)=-1$.
Let $T$ be a toric Weil divisor of $(X\ni P)$ such that
$T_{Z_0}\cap S_0=E^2_0$.
Then
$K_{Y_0}+(S_1)_{Y_0}+(S_0)_{Y_0}+(D_0)_{Y_0}+T_{Y_0}\sim 0$ is lc by Inversion of Adjunction.
The ray $\RR_+[(F_0)_{Y_0}]$ gives the divisorial contraction of
$(S_0)_{Y_0}$ onto a curve, denoted by $h'_0$ in our diagram.
We obtain a non-toric blow-up $g_1\colon (Z_1,S_1)\to (X\supset C\ni P)$, where $S_1=\Exc g_1$, $g_1(S_1)=C$ and $(S_1)_{Y_0}\cong S_1$.
Since $g_1$ be a toric blow-up (under identification) in some neighborhood of any point other than $P$, then 
$\Diff_{S_1}(0)=\frac{w^3_1-1}{w^3_1}E^2_1+\frac{w^j_1-1}{w^j_1}E^1_1+\frac{d_1-1}{d_1}(F_1)_{Z_0}$, $j\in\{1,2\}$, $E^2_1=h'_0((S_0)_{Y_0})$ and $E^1_1$ are some sections,
$F_1$ is a fiber over $P$, $w^3_1\in\ZZ_{\ge 3}$ and $d_1\in\ZZ_{\geq 1}$. Hence $g_1$ is a plt blow-up.

{\it Second step}.
Assume that there exists a curve $\G_1\subset (S_1)_{Y_0}$ with the following two properties:
1) $K_{(S_1)_{Y_0}}+\Diff_{(S_1)_{Y_0}}(0)+\G_1$ is a plt and $h_0$-anti-ample divisor, $h_0\colon \G_1\to\G_0$ is a surjective morphism and
2) $\G_1$ is not a center of any blow-up of
${\mathcal H_0}$, that is, $\G_1$ is a non-toric subvariety of $((S_1)_{Y_0},\Diff_{(S_1)_{Y_0}}(0))$ in any analytical
neighborhood of the fiber $(F_1)_{Y_0}$ over $P$.

The triples $((S_1)_{Y_0},\Diff_{(S_1)_{Y_0}}(0), \G_1)$ and
$(S_0,\Diff_{S_0}(0), \G_0)$
have the same structures and (with similar notation)
$w^1_1=1$, 
$Q_1=(E^2_1)_{Y_0}\cap (F_1)_{Y_0}$, $\G_1\cap (F_1)_{Y_0}=Q_1$, 
$w^2_1\ge 1$, $d_1=1$, $((S_1)_{Y_0}\ni Q_1)$ is of type $\frac1{r_1}(1,1)$ ($r_1\ge 1$) and $\G_1\cdot (F_1)_{Y_0}=\G_1\cdot (E_1^2)_{Y_0}=\frac1{r_1}$.

Consider an arbitrary toric structure of $Y_0$ in any neighborhood of the point $Q_1$ such that 
$\G_1$ is also a toric subvariety of $Y_0$. 
Let $h_1\colon (Y_1,(S_2)_{Y_1})\to (Y_0\supset \G_1\ni Q_1)$ be an arbitrary toric blow-up of the curve 
$\G_1$ with an unique exceptional divisor ($\Exc h_1=(S_2)_{Y_1}$),
$(S_1)_{Y_1}\cong (S_1)_{Y_0}$.

The set of all possible blow-ups $h_1$ for any toric structure of $(Y_0\ni Q_1, \G_1)$ is denoted by $\mathcal H_1$.

Let $(D_1)_{Y_0}$ be a toric Weil divisor of $(Y_0\ni Q_1)$ such that $(D_1)_{Y_0}|_{S_1}=\G_1$ and $a((S_2)_{Y_1},(D_1)_{Y_0}+(S_0)_{Y_0}+(S_1)_{Y_0})=-1$. 
We have 1-complement $K_{Y_1}+(S_2)_{Y_1}+(S_1)_{Y_1}+(S_0)_{Y_1}+(D_1)_{Y_1}\sim 0/X$ by Inversion of Adjunction applied to the surfaces
$(S_i)_{Y_1}$. By the cone theorem we have:

1) there exists a divisorial contraction $h'_{1,1}\colon Y_1\to Y_{1,1}$ of $(S_1)_{Y_1}$
onto a curve, $(S_2)_{Y_1}\cong (S_2)_{Y_{1,1}}$;

2) there exists a small contraction $\varphi_{1,1}$ of an extremal ray generated  by $(F_0)_{Y_{1,1}}$. Let
$\varphi^+_{1,1}$ be a log flip of $\varphi_{1,1}$, 
$\Exc \varphi^+_{1,1}=(F^+_0)_{Y_{1,2}}$,
$h'_{1,2}\colon Y_{1,1}\dashrightarrow Y_{1,2}$ be a corresponding birational map;

3) there exists a divisorial contraction $h'_{1,3}\colon Y_{1,2}\to Z_2$ of $(S_0)_{Y_{1,2}}$
onto a curve.

Thus we obtain a birational map 
$h'_1=h'_{1,3}\circ h'_{1,2}\circ h'_{1,1}\colon Y_1\dashrightarrow Z_2$. Put $S_2=(S_2)_{Z_2}$.
Since $(E^2_0)_{Y_{1,1}}\cap (F_0)_{Y_{1,1}}=(Q_0)_{Y_{1,1}}$ then
$(D_1)_{Y_{1,1}}\cdot (F_0)_{Y_{1,1}}>0$ and
the divisor $(D_1)_{Z_2}$ contains the fiber $(F^+_0)_{Z_2}$ and two sections of the local conic bundle $S_2\to C$, $\rho(S_2/C)=1$, $K_{Z_2}+S_2+(D_1)_{Z_2}\sim 0/X$ is lc.
By Shokurov's criterion on the
characterization of toric varieties $(S_2, \Diff_{S_2}(0))\to C$ is a toric conic bundle \cite{Sh2}.
We obtain a non-toric plt blow-up $g_2\colon (Z_2,S_2)\to (X\supset C\ni P)$.

We prove the following proposition.
\begin{proposition}
The pair $(S_i,\Diff_{S_i}(0))$ is klt and local toric conic bundle $(1$-complementary$)$, $\rho(S_i/C)=1$, $g_i$ is a non-toric plt blow-up for $i=1, 2$.
\end{proposition}

{\it Third step}.
Assume that there exists a curve $\G_2\subset (S_2)_{Y_1}$ with the following two properties:
1) $K_{(S_2)_{Y_1}}+\Diff_{(S_2)_{Y_1}}(0)+\G_2$ is a plt and $h_1$-anti-ample divisor, $h_1\colon \G_2\to\G_1$ is a surjective morphism and
2) $\G_2$ is not a center of any blow-up of
${\mathcal H_1}$, that is, $\G_2$ is a non-toric subvariety of $((S_2)_{Y_1},\Diff_{(S_2)_{Y_1}}(0))$ in any analytical
neighborhood of the central fiber $F_2$ of $(S_2)_{Y_1}$ over $P$.

The triple $((S_2)_{Y_1},\Diff_{(S_2)_{Y_1}}(0), \G_2)$ has the same structures as the previous ones. In particular (with similar notation), $w^1_2=1$ and $w^2_2\ge 1$.

\begin{proposition}
There is no any blow-up 
$h_2\colon (Y_2, (S_3)_{Y_2})\to (Y_1\supset \G_2)$ of the curve $\G_2$ with unique exceptional divisor such that $(S_3)_{Y_2}$ is realized by some plt blow-up $g_3\colon (Z_3, (S_3)_{Z_3})\to (X\supset C\ni P)$.
\begin{proof}
Assume the converse. Consider a general point of $C$.
Let $F_3$ be a fiber of $(S_3)_{Y_2}$ over $P$.
Put $\Theta=\Diff_{(S_3)_{Z_3}}(0)$ for simplicity.
Since $w^2_0+w^2_1+w^2_2+1\ge 5$ then $\Theta$ has some component (a section of conic bundle) with a coefficient $\ge 4/5$. 

We claim that $K_{(S_3)_{Z_3}}+\Theta$ is 1 or 2-complementary. Assume that $K_{(S_3)_{Z_3}}+\Theta$ is not 1-complementary. Then the divisor $K_{(S_3)_{Z_3}}+\alpha F_3+\Theta$ is lc, but not plt for some $\alpha\le 1$, and consider its inductive blow-up $\sigma\colon \widetilde X\to (S_3)_{Z_3}$ with exceptional divisor $\widetilde E$. The curve $(F_3)_{\widetilde X}$ can be contracted in the appropriate MMP over $C$. Denote this contraction by $\widetilde X\to\overline X$. The divisor $K_{\overline X}+\overline E+\Theta_{\overline X}$ is plt. 

Let $K_{\widetilde X}+\widetilde E+\Theta_{\widetilde X}$ be nonnegative on $(F_3)_{\widetilde X}$. We can extend complement of $K_{\overline E}+\Diff_{\overline E}(\Theta_{\overline X})$ on $\overline X$, pull back on $\widetilde X$ and push-down them on $(S_3)_{Z_3}$.
There are only two cases: 1) $\Diff_{\overline E}(\Theta_{\overline X})=1/2P_1+1/2P_2+(1-1/m)P_3$ and 2) $\Diff_{\overline E}(\Theta_{\overline X})=1/2P_1+2/3P_2+4/5P_3$, where $\{P_i\}$ are some points, $m\ge 5$. We obtain 2- or 6-complement.  

Let $K_{\widetilde X}+\widetilde E+\Theta_{\widetilde X}$ be negative on $(F_3)_{\widetilde X}$. The divisor   
$-(K_{\widetilde X}+\widetilde E+\Theta_{\widetilde X})$ is ample over $C$. Similarly 2- or 6-complement of $K_{\widetilde E}+\Diff_{\widetilde E}(\Theta_{\widetilde X})$ can be extended on $\widetilde X$ and we have 2- or 6-complement $D^+$ of $K_X$
with $a((S_3)_{Y_2}, D^+)=-1$. 

Consider the case of 6-complement. Since $a((S_3)_{Y_2}, D^+)=-1$ then there is one possibility $a((S_0)_{Y_2}, D^+)=-1/2$, $D^+|_{S_0}=(7/6)\G_0+...$ and $a((S_1)_{Y_2}, D^+)\le -2/3$.
Since $F_3\subset (S_i)_{Y_2}$ for $i=0, 1$ then $K_{Y_2}+a((S_0)_{Y_2}, D^+)(S_0)_{Y_2}+a((S_1)_{Y_2}, D^+)(S_1)_{Y_2}+(S_3)_{Y_2}$ is not lc, the contradiction. 

Thus we have 1- or 2-complement. Therefore the coefficients of $D^+$ are equal 1 or 1/2 and $a((S_0)_{Y_2}, D^+)\le-1/2$. We have the same contradiction as above.
\end{proof}
\end{proposition}

\par {\bf A2)} Let $(X\ni P)$ be a non-$\QQ$-factorial terminal toric three-dimensional singularity, that is, 
$(X\ni P)\cong (\{x_1x_2+x_3x_4=0\}\subset (\CC^4_{x_1x_2x_3x_4},0))$ by Proposition \ref{termtoric}.

Let $f\colon (Y, E)\to (X\supset C\ni P)$ be some plt blow-up.
Let $\varphi_i\colon X_i\to (X\ni P)$ be any of two $\QQ$-factorializations, $\Exc\varphi_i=C_i$ ($i=1, 2$).
Let $\psi_i\colon(Y_i, E_i)\to (X_i\supset C_{X_i}\ni P_{X_i})$ be a plt blow-up of $C_{X_i}$ such that $E_i$ and $E$ define the same discrete valuation of the function field $\mathcal K(X)$, $\rho(E_i/C)=1$. The blow-up $\psi_i$ was constructed in the previous case of $\QQ$-factorial singularities. Let $Y_i\dashrightarrow Y$ be a log flip for the curve $(C_i)_{Y_i}$. Thus $f$ has constructed and $\rho(E/C)=2$.  

We give another construction and prove that $(E,\Diff_E(0))\to C$ is a toric conic bundle by the procedure illustrated on the next diagram (Fig. 2) and defined below.

\[
\xymatrix{& Y_0 \ar[dl]_{h_0} \ar@{-->}[d]^{h'_0} & Y_1 \ar[l]_{h_1} \ar@{-->}[d]^{h'_1}\\
	Z_0 \ar[d]_{g_0}& Z_1 \ar[d]_{g_1} & Z_2 \ar[d]_{g_2}\\
	X& X& X
}
\]\begin{center}
	Fig. $2$
\end{center}

{\it First step}.
Let $g_0\colon (Z_0,S_0)\to (X\supset C\ni P)$ be any toric plt blow-up, where $g_0(S_0)=C$. Its description is given in example \ref{ex1} 2), whose notation is used.
Let $F_0=F^1_0+F^2_0$ be a fiber over the point $P$. Put $Q_0=F^1_0\cap F^2_0$. 

{\it Second step}.
Assume that there exists a curve $\G_0\subset S_0$ with the following two properties:
1) $K_{S_0}+\Diff_{S_0}(0)+\G_0$ is a plt and $g_0$-anti-ample divisor; 2) $\G_0$ is a non-toric subvariety in any  
analytical neighborhood of the fiber $F_0$
on the toric variety $Z_0$ for any toric structure of $(X\ni P)$, that is, the curve $\G_0$ is a non-toric subvariety of $(S_0,\Diff_{S_0}(0))$ in any analytical
neighborhood of $F_0$ on $S_0$.

Considering a fiber over a general point of $C$ we have $a_2=1$ or $a_3=1$. To be definite, put $a_2=1$ and $F_0^2\cap E_2\ne \emptyset$. By simple calculations $\G_0\cap (F^1_0\cup F^2_0)=Q_0$, $F_0^1\cdot \G_0=\frac{a_3}{a_3+1}$ and $F^2_0\cdot \G_0=\frac1{a_3+1}$.

Consider an arbitrary toric structure of $Z_0$ in any neighborhood of the point $Q_0$ such that 
$\G_0$ is a toric subvariety of $Z_0$ also. Let $h_0\colon (Y_0,(S_1)_{Y_0})\to (Z_0\supset \G_0\ni Q_0)$ be an arbitrary toric blow-up of the curve 
$\G_0$ with an unique exceptional divisor ($\Exc h_0=(S_1)_{Y_0}$), $(S_0)_{Y_0}\cong S_0$.

The set of all possible blow-ups $h_0$ for any toric structure of $(Z_0\ni Q_0, \G_0)$ is denoted by $\mathcal H_0$.

Let $(D_0)_{Z_0}$ be a toric Weil divisor of $(Z_0\ni Q_0)$ such that $(D_0)_{Z_0}|_{S_0}=\G_0$ and $a((S_1)_{Y_0}, (D_0)_{Z_0}+S_0)=-1$.
Let $T_1$ and $T_2$ be toric Weil divisors of $(X\ni P)$ such that
$$
K_{S_0}+\Diff_{S_0}((T_1+T_2)_{Z_0}+(D_0)_{Z_0})=K_{S_0}+F^2_0+E_2+\G_0\sim 0.
$$

The pairs $(X\ni P, T_1+T_2+(D_0)_X)$ and $((S_1)_{Y_0},\Diff_{(S_1)_{Y_0}}((T_1+T_2)_{Y_0}+(D_0)_{Y_0}+(S_0)_{Y_0}))$ are lc. 
Since $T_1+T_2$ is Cartier divisor then $(D_0)_X$ is Cartier divisor. The curves $(F^1_0)_{Y_0}$ and $(F^2_0)_{Y_0}$ generate extremal rays of $\NE(Y_0/X)$ that give small contractions.
Let $h'_{0,1}\colon Y_0\dashrightarrow Y_{0,1}$ be any of two log flips. Since our pairs are lc then $\rho((S_0)_{Y_{0,1}}/C)=1$.
Let $h'_{0,2}\colon Y_{0,1}\to Z_1$ be a divisorial contraction of $(S_0)_{Y_{0,1}}$ onto a curve.

Thus we obtain a birational map $h'_0=h'_{0,2}\circ h'_{0,1}\colon Y_0\dashrightarrow Z_1$ and a non-toric blow-up $g_1\colon (Z_1, S_1)\to (X\supset C\ni P)$, $\rho(S_1/C)=2$. It can be proved by direct computation that $-S_1$ is $g_1$-ample divisor, but if we consider the construction of $g_1$ through two $\QQ$-factorializations of $(X\ni P)$ as done above, then it is obvious that the divisor $-S_1$ is $g_1$-ample.
The divisor $\Diff_{S_1}((T_1+T_2)_{Z_1}+(D_0)_{Z_1})$ consists of four curves and gives 1-complement of $K_{S_1}+\Diff_{S_1}(0)$.
By Shokurov's criterion on the characterization of toric varieties $(S_1, \Diff_{S_1}((T_1+T_2)_{Z_1}+(D_0)_{Z_1})\to C$ is a toric conic bundle \cite{Sh2}. Thus $g_1$ is a plt blow-up.

{\it Third step}.
Assume that there exists a curve $\G_1\subset (S_1)_{Y_0}$ with the following two properties:
1) $K_{(S_1)_{Y_0}}+\Diff_{(S_1)_{Y_0}}(0)+\G_1$ is plt and $h_0$-anti-ample divisor, $h_0\colon \G_1\to(\G_0)_{Z_0}$ is a surjective morphism and
2) $\G_1$ is not a center of any blow-up of
${\mathcal H_0}$, that is, $\G_1$ is a non-toric subvariety of $((S_1)_{Y_0},\Diff_{(S_1)_{Y_0}}(0))$ in any analytical
neighborhood of the central fiber $F_1$ of $(S_1)_{Y_0}$ over $P$.

The triple $((S_1)_{Y_0},\Diff_{(S_1)_{Y_0}}(0), \G_1)$ has the same structures as in the previous case of $\QQ$-factorial singularities, and we use its notation.

Consider an arbitrary toric structure of $Y_0$ in any neighborhood of the point $Q_1$ such that 
$\G_1$ is also a toric subvariety of $Y_0$. 
Let $h_1\colon (Y_1,(S_2)_{Y_1})\to (Y_0\supset \G_1\ni Q_1)$ be an arbitrary toric blow-up of the curve 
$\G_1$ with an unique exceptional divisor ($\Exc h_1=(S_2)_{Y_1}$), $(S_1)_{Y_1}\cong (S_1)_{Y_0}$.

Let $(D_1)_{Y_0}$ be a toric Weil divisor of $(Y_0\ni Q_1)$ such that $(D_1)_{Y_0}|_{S_1}=\G_1$ and $a((S_2)_{Y_1},(D_1)_{Y_0}+(S_0)_{Y_0}+(S_1)_{Y_0})=-1$. Considering the case of $\QQ$-factorial singularities and construction of $g_0\circ h_0$ through $\QQ$-factorializations of $(X\ni P)$ we have $(E_2)_{Y_0}\subset (D_1)_{Y_0}$ and hence $F^2_0\subset (D_1)_{Y_0}$. Thus we have 1-complement $K_{Y_1}+(S_2)_{Y_1}+(S_1)_{Y_1}+(S_0)_{Y_1}+(D_1)_{Y_1}\sim 0/X$ by Inversion of Adjunction applied to the surfaces
$(S_i)_{Y_1}$. By the cone theorem we have:

1) there exists a divisorial contraction $h'_{1,1}\colon Y_1\to Y_{1,1}$ of $(S_1)_{Y_1}$
onto a curve, $(S_2)_{Y_1}\cong (S_2)_{Y_{1,1}}$;

2) there exists a small contraction of $(F^1_0)_{Y_{1,1}}$, $h'_{1,2}\colon Y_{1,1}\dashrightarrow Y_{1,2}$ is a corresponding log flip;

3) there exists a small contraction of $(F^2_0)_{Y_{1,2}}$, $h'_{1,3}\colon Y_{1,2}\dashrightarrow Y_{1,3}$ is a corresponding log flip;

4) there exists an divisorial contraction $h'_{1,4}\colon Y_{1,3}\to Z_2$ of $(S_0)_{Y_{1,3}}$
onto a curve.

Thus we obtain a birational map 
$h'_1=h'_{1,4}\circ h'_{1,3}\circ h'_{1,2}\circ h'_{1,1}\colon Y_1\dashrightarrow Z_2$, the local conic bundle $(S_2)_{Z_2}\to C$, $\rho((S_2)_{Z_2}/C)=2$ and $K_{Z_2}+(S_2)_{Z_2}+(D_1)_{Z_2}\sim 0/X$ is lc. Let $F_2=F^1_2+F^2_2$ be a fiber over $P$ and the curves $F^1_2$, $F^2_2$ appear due to log flips $h'_{1,2}$, $h'_{1,3}$ respectively. By the construction the divisor $(D_1)_{Z_2}$ contains two sections of $(S_2)_{Z_2}$ and 
$F^1_2$.

If we consider this construction through two $\QQ$-factorializations of $(X\ni P)$ then $(S_2)_{Z_2}$ is anti-ample over $C$ and $(F^1_0)_{Y_{1,1}}\cap (F^2_0)_{Y_{1,1}}=(Q_0)_{Y_{1,1}}$.
Since $(F_0^2)_{Y_{1,2}}\cdot (F_0^2)_{Y_{1,2}}=0$, $K_{Y_{1,2}}+(S_0)_{Y_{1,2}}+(S_2)_{Y_{1,2}}+(D_1)_{Y_{1,2}}\sim 0$ then for some $e>0$ we have $(D_1)_{Y_{1,2}}\cdot (F^2_0)_{Y_{1,2}}=e(E_2)_{Y_{1,2}}\cdot (F_0^2)_{Y_{1,2}}>0$ and
$(D_1)_{Z_2}$ contains $F^2_2$.

By Shokurov's criterion on the
characterization of toric varieties $((S_2)_{Z_2}, \Diff_{(S_2)_{Z_2}}(0))\to C$ is a toric conic bundle \cite{Sh2}.
We obtain a non-toric plt blow-up $g_2\colon (Z_2,S_2)\to (X\supset C\ni P)$, where $S_2=(S_2)_{Z_2}$.

We prove the following proposition.
\begin{proposition}
The pair $(S_i,\Diff_{S_i}(0))$ is klt and local toric conic bundle $(1$-complementary$)$, $\rho(S_i/C)=1$, $g_i$ is a non-toric plt blow-up for $i=1, 2$.
\end{proposition}

\end{example}

\begin{example} \label{candim1}
Let us describe the non-toric canonical blow-ups (they will be non-terminal blow-ups always)
$g\colon (Y,E)\to (X\supset C\ni P)$ provided that $(X\ni P)$ is
a toric terminal singularity, $C=g(E)$ is a toric (smooth) subvariety and $\dim C=1$.
Depending on a type of $(X\ni P)$ we consider two Cases {\bf B1)} and {\bf B2)}.

{\bf B1)} 
Let $(X\ni P)$ be a $\QQ$-factorial terminal singularity.
Let $g\colon (Z,S)\to (X\supset C\ni P)$ be any toric canonical blow-up (see Proposition \ref{cantoricdim1}).

Assume that there exists a curve $\G\subset S$ with the following two properties:
1) $K_S+\Diff_S(0)+\G$ is $g$-anti-ample divisor, and $\G$ does not contain any center of canonical singularities of $Z$; 2) $\G$ is a non-toric subvariety in any  
analytical neighborhood of the fiber $F$ (over $P$) on the toric variety $Z$ for any toric structure of $(X\ni P)$, that is, the curve $\G$ is a non-toric subvariety of $(S,\Diff_S(0))$ in any analytical neighborhood of $F$ on $S$.

Thus $(X\ni P)$ is a smooth point, $S$ is a smooth surface,
$\Diff_S(0)=\frac{k-1}{k}E$, where $k\geq 2$ and $E$ is some section by Proposition \ref{cantoricdim1}.
By adjunction formula $\G$ is smooth,
$Q=\G\cap F\cap E$, $\G\cdot F=1$. 

Let $(X\ni P, D)$ be any pair with canonical singularities such that $D$ is a boundary, $\G\in \mt{CS}(Z,D_Z-a(S,D)S)$.
Obviously, $D_Z|_S=\G+aF$ and $a(S, D)=0$, where $a\ge 0$.

Considering the blow-up $(\CC^3_{x_1x_2x_3}\ni 0)\cong (X\supset C\ni P)$ with weights $(k,1,0)$, $C=\{x_1=x_2=0\}$ and the divisor given by the equation
$x_1^2+x_1x_2+x_1x_3^m+bx_2^k=0$, then clearly, there is a divisor $D$ for any such curve $\G$.

By Theorem \ref{caninductive} there exists a divisorial contraction 
$h\colon (\widetilde Y,\widetilde E)\to (Z\supset \G)$ such that $a(\widetilde E, D)=0$,
$\Exc h=\widetilde E$ is an irreducible divisor and $h(\widetilde E)=\G$.
Applly $K_{\widetilde Y}+D_{\widetilde Y}+\varepsilon \widetilde S$--MMP. Since $\rho(\widetilde Y/X)=2$ and  
$K_{\widetilde Y}+D_{\widetilde Y}+\varepsilon \widetilde S\equiv \varepsilon \widetilde S$ over $X$, then after log flips 
$\widetilde Y\dashrightarrow \overline Y$ (perhaps their lack) 
we obtain a divisorial contraction $h'\colon \overline Y\to Y$, which contracts $\overline S$ onto a curve $C_Y$. 

Thus we obtain a non-toric canonical blow-up $f$.
Since $C_Y\in\mt{CS}(Y)$ by the construction then $f$ is not a terminal blow-up.

{\bf B2)} Let $(X\ni P)$ be a non-$\QQ$-factorial terminal toric three-dimensional singularity, that is,
$(X\ni P)\cong (\{x_1x_2+x_3x_4=0\}\subset (\CC^4_{x_1x_2x_3x_4},0))$.
Consider a $\QQ$-factorialization $g\colon\widetilde X\to X$, $\widetilde T=\Exc g$ and $\widetilde P=\widetilde T\cap \widetilde C$. 
We apply the construction from {\bf B1)} for the curve
$\widetilde C\subset (\widetilde X\ni \widetilde P)$ such that the divisor $D$ from the construction has the form $g^*D_X$, where $D_X$ is a $\QQ$-Cartier divisor. We obtain a non-toric canonical blow-up $f\colon Y^+\to \widetilde X$.
Let $Y^+\dashrightarrow Y$ be a log flip for the curve $T_{Y^+}$. 
Thus we obtain a required non-toric canonical
blow-up $f$ (anti-amplness of $E$ is proved as in case {\bf A2)}).

Let us describe the curves $\G$. Let $ g\colon (Z, S) \to (\widetilde X \ni \widetilde P) $ be a toric canonical blow-up
obtained in the first step of the construction. Let $\psi \colon Z \dashrightarrow Z^+$ be a toric log flip for the curve $T_Z$. So
$g^+ \colon (Z^+, S^+) \to (X \ni P)$ is a toric canonical blow-up.
The structure of the curve $ \G_{S^+}$ is completely identical to the structure of the curve $\G$ considered in case {\bf A2)}.
To prove that any such curve $\G_{S^+}$ is realizable, it suffices to consider a divisor of the form $x_ {i_1}+bx_{i_2}^k=0$ on $(X\ni P)$ for some $b$, $k$, $\{i_1, i_2\}=\{1,2\}$ or $\{3,4\}$.
\end{example}

\begin{theorem}\label{dim31_plt}
Let $f\colon (Y,E) \to (X\supset C\ni P)$ be a plt blow-up of three-dimensional toric terminal singularity, where
$\dim f(E)=1$. Assume that the curve $C=f(E)$ is a toric subvariety of $(X\ni P)$.
Then, either $f$ is a toric morphism $($see Example $\ref{ex1})$, or $f$ is a non-toric morphism described in Example $\ref{pltdim1}$.
\begin{proof}
By Example \ref{pltdim1} we must only consider the case when $(X\ni P)$ is a $\QQ$--factorial singularity.
Let $f$ be a non-toric morphism (up to analytic isomorphism). 
Let $D_Y\in |-n(K_Y+E)|$ be a general element for $n\gg 0$. Put $D_X=f(D_Y)$ and $d=\frac{1}{n}$.
The pair $(X,dD_X)$ is log canonical, $a(E,dD_X)=-1$, and $E$ is a unique exceptional divisor with discrepancy $-1$.

By the construction of partial resolution of $(X,dD_X)$ (see Definition \ref{nondegen} and the paper \cite{Var}) and by Criterion \ref{toricplt},
there exists a toric divisorial contraction
$g\colon Z\to X$ dominated by partial resolution 
of $(X,dD_X)$ (up to toric log flips) and the following properties are fulfilled.
\begin{enumerate}
\item[A)] The exceptional set $\Exc g=S$ is an irreducible divisor, 
the divisors $S$ and $E$ define the different discrete valuations of the function field
$\mathcal K(X)$, and $g(S)=C$.
\item[B)] By $\G$ denote the center of $E$ on the surface $S$. Then the curve $\G$ is a non-toric subvariety of $Z$.
In the other words, $\G$ is a non-toric subvariety of
$(S,\Diff_S(0))$.
\end{enumerate}

Obviously, $a(S_0, dD_X)<0$.
By Example \ref{pltdim1} (in its notation) we must prove only that the anti-ample over $X$ divisor $K_{S_0}+\Diff_{S_0}(0)+\G_0$ is plt
in some analytical neighborhood of the fiber $F_0\subset S_0$.
We can choose the divisor $dD_X$ such that $\Supp (dD_X|_{S_0})\subset \G_0\cup F\cup \G'_0 \cup E^2_0$, where $\G'_0$ is a general divisor on $S_0$.

Assume that $K_{S_0}+\Diff_{S_0}(0)+\G_0$ is not a plt divisor.
By the adjunction formula the curve $\G_0$ is smooth.
By connectedness lemma
$K_{S_0}+\Diff_{S_0}(0)+\G_0$ is not a plt divisor at unique point, and denote this point by $G_0$.
The point $G_0$ is a non-toric subvariety of $(S_0,\Diff_{S_0}(0))$. Moreover, the curve $\G_0$ is locally a non-toric subvariety at the point $G_0$ only.
By the construction of partial resolution
\cite{Var} (in a small analytical neighborhood of the point $G_0$) there exists a divisorial toric contraction  
$\widehat g_0\colon\widehat Z_0\to Z_0$ such that
$\Exc \widehat g_0=S''_0$ is an irreducible divisor, $\widehat g(S''_0)=G_0$ and the two following conditions are satisfied.
\par 1)  Put $S'_0=(S_0)_{\widehat Z_0}$ and $C_0=S'_0\cap S''_0$. Let $c(\G_0)$ be the log canonical threshold of $\G_0$ for the pair $(S_0,\Diff_{S_0}(0))$. 
Then ${\widehat g_0}|_{S'_0}\colon S'_0\to S_0$ is the toric inductive blow-up of $K_{S_0}+\Diff_{S_0}(0)+c(\G_0)\G_0$ 
(see Theorems \ref{inductive} and \ref{dim2plt}), and
the point $\widehat G_0=C_0\cap(\G_0)_{S'_0}$ is a non-toric subvariety of $(S''_0,\Diff_{S''_0}(0))$.
\par 2) The divisor $\Diff_{S''_0}(dD_{\widehat Z_0}+a(S_0,dD_X)S'_0)$ 
is a boundary in some small analytical neighborhood of the point $\widehat G_0$.

Let $H$ be a general hyperplane section of sufficiently large degree passing through the point $P$ such that it does not contain the curve $C$.
Then there exists a number $h>0$ such that $a(S''_0,dD_X+hH)>-1$, and the point $\widehat G_0$ is a center of 
$(S''_0,\Diff_{S''_0}(dD_{\widehat Z_0}+a(S_0,dD_X)S'_0+hH_{\widehat Z_0}))$.
Therefore we obtain a contradiction for the pair
$(S''_0,\Diff_{S''_0}(dD_{\widehat Z_0}+a(S_0,dD_X)S'_0+hH_{\widehat Z_0}))$ and the point $\widehat G_0$
by Theorem \ref{point}.
\end{proof}
\end{theorem}

We have proved the next theorem too.
\begin{theorem}\label{dim32_plt}
Let $f\colon (Y,E) \to (X\supset C \ni P)$ be a plt blow-up of three-dimensional toric $\QQ$-factorial singularity, where
$\dim f(E)=1$. Assume that the curve $C=f(E)$ is a toric subvariety of $(X\ni P)$.
Then, either $f$ is a toric morphism $($see Example $\ref{ex1})$, or $f$ is a non-toric morphism described in Example $\ref{pltdim1}$.
\end{theorem}

\begin{theorem}\label{dim31_can}
Let $f\colon (Y,E) \to (X\supset C \ni P)$ be a canonical blow-up of three-dimensional toric terminal singularity, where
$\dim f(E)=1$. Assume that the curve $C=f(E)$ is a toric subvariety of $(X\ni P)$.
Then, either $f$ is a toric morphism $($see Proposition $\ref{cantoricdim1})$, or $f$ is a non-toric morphism and described in Example $\ref{candim1}$.
\begin{proof}
Let $f$ be a non-toric morphism (up to analytic isomorphism).  
Let $D_Y\in |-nK_Y|$ be a general element for $n\gg 0$. Put $D_X=f(D_Y)$ and $d=\frac{1}{n}$.
The pair $(X,dD_X)$ has canonical singularities and $a(E,dD_X)=0$. Now the arguments of the proof of Theorem
\ref{dim31_plt} can be obviously applied, and we have $a(S,dD_X)=0$, this completes the proof.
\end{proof}
\end{theorem}

\begin{corollary}\label{dim32_can}
Under the same assumption as in Theorem
$\ref{dim31_can}$ the two following statements are satisfied:
\par $1)$ \cite{Kaw1} if $f$ is a terminal blow-up then the $($toric$)$ morphism $f$ is isomorphic to the blow-up of the ideal of the curve
$C$ and an index of $(X\ni P)$ is equal to $1$, that is, either $(X\ni P)$ is a smooth point or 
$(X\ni P)\cong (\{x_1x_2+x_3x_4=0\}\subset (\CC^4_{x_1x_2x_3x_4},0));$
\par $2)$ if $f$ is a non-toric morphism then an index of $(X\ni P)$ is equal to $1$.
\end{corollary}

\section{\bf Toric log surfaces}
\begin{definition}\label{odpdef}
Let $\PP({\bf w})=\PP_{x_1x_2x_3x_4}(w_1,w_2,w_3,w_4)$, where $w_1+w_2=w_3+w_4$ and $\gcd(w_1,w_2,w_3,w_4)=1$. 
Put $(w_1,w_2,w_3,w_4)=(a_1d_{23}d_{24}, a_2d_{13}d_{14},
a_3d_{14}d_{24}, a_4d_{13}d_{23})$, where $d_{ij}=\gcd(w_k,w_l)$ and $i,j,k,l$ are mutually distinct indices from $1$ to $4$.
The toric pair
$$\big(S,D\big)=\big(x_1x_2+x_3x_4\subset\PP({\bf w}),\Diff_{S/\PP({\bf w})}(0)\big)$$
is called {\it an ODP pair}, $\rho(S)=2$. By Proposition 1.6 of \cite{Kud2} we have
$D=\sum_{i<j, 1\leq i\leq 2}\frac{d_{ij}-1}{d_{ij}}C_{ij}$, where $C_{ij}=\{x_i=x_j=0\}\cap S$.

Let $f\colon (Y,E)\to (X\ni P)$ be a toric plt blow-up of three-dimensional ordinary double point. 
Then $(E,\Diff_E(0))$ is an ODP pair by Example \ref{ex1}. 
Converse is also true: every ODP pair is realized by some toric plt blow-up of three-dimensional ordinary double point. 

To be definite, assume that $w_1\leq w_2$, $w_3\leq w_4$, $w_2\leq w_4$, $P_1=(1:0:0:0)$, $\ldots$, $P_4=(0:0:0:1)$. 
The surface $S$ has a cyclic singularity at the point $P_i$ for every
$i=1,2,3,4$ (see Fig. 3).
\\
\begin{center}
\begin{picture}(300,140)(20,0)
\put(0,30){\line(1,0){280}}
\put(30,0){\line(0,1){140}}
\put(110,10){$C_{14}=\{x_1=x_4=0\}$}
\put(30,30){\circle*{8}}
\put(35,15){$\frac1{a_2}\big{(}\frac{w_3}{d_{14}},\frac{w_4}{d_{13}}\big{)}$}
\put(250,0){\line(0,1){140}}
\put(250,30){\circle*{8}}
\put(255,15){$\frac1{a_3}\big{(}\frac{w_1}{d_{24}},\frac{w_2}{d_{14}}\big{)}$}
\put(40,60){$C_{13}=\{x_1=x_3=0\}$}
\put(260,60){$C_{24}=\{x_2=x_4=0\}$}
\put(255,95){$\frac1{a_1}\big{(}\frac{w_3}{d_{24}},\frac{w_4}{d_{23}}\big{)}$}
\put(250,110){\circle*{8}}
\put(30,110){\circle*{8}}
\put(0,110){\line(1,0){280}}
\put(95,123){$C_{23}=\{x_2=x_3=0\}$}
\put(35,95){$\frac1{a_4}\big{(}\frac{w_1}{d_{23}},\frac{w_2}{d_{13}}\big{)}$}
\put(10,95){$P_4$}
\put(10,15){$P_2$}
\put(231,95){$P_1$}
\put(231,15){$P_3$}
\end{picture}

Fig. $3$
\end{center}

Since $\OO_{\PP({\bf w})}(w_i)|_S=\{x_i=0\}|_S=\frac1{d_{ik}}C_{ik}+\frac1{d_{il}}C_{il}$ for the corresponding different 
indices $k$ and $l$, then it is easy to calculate that  
$C_{13}^2=d_{13}^2(w_3-w_2)/(w_2w_4)\leq 0$, $C_{23}^2=d_{23}^2(w_2-w_4)/(w_1w_4)\leq 0$, $C_{14}^2=d_{14}^2(w_4-w_2)/(w_2w_3)\geq 0$ and
$C_{24}^2=d_{24}^2(w_2-w_3)/(w_1w_3)\geq 0$. In particular, Mori cone $\NE(S)$ is generated by the two rays $\RR_+[C_{13}]$, $\RR_+[C_{23}]$.
\end{definition}

Now we prove a two-dimensional {\it non-toric point theorem}.  
An one-dimensional analog ($\dim S=1$) of Theorem \ref{point} 1) is obvious (see the proofs of Theorems \ref{dim2plt} and \ref{dim2can} also).

\begin{theorem} \label{point}
Let $(S,D)$ be a toric pair, where $S$ is a normal projective surface.
Assume that $D=\sum_{i=1}^rd_iD_i$, where $D_i$ is a prime divisor and $\frac12 \leq d_i\leq 1$ for each $i$.
Assume that there exists the boundary $T$ such that $T\geq D$ and $-(K_S+T)$ is an ample divisor.
Assume that some point $\Gamma$ is 
a center of $\LCS(S,T)$, and there exists the analytical neighborhood $U$ of $\Gamma$ such that
$K_S+T$ is a log canonical divisor in the punctured neighborhood $U\backslash\Gamma$. Then
the point $\Gamma$ is a toric subvariety of $(S,D)$ if one of the two following conditions is satisfied:

$1)$ $\rho(S)=1$; 

$2)$ $\rho(S)=2$, two different extremal rays of $\NE(S)$ give two toric conic bundles;

$3)$ $(S,D)$ is ODP pair.
\begin{proof} Let the point $\Gamma$ be a non-toric subvariety of $(S,D)$. We will obtain a contradiction.

Consider Condition 1).
It is clear that this theorem is sufficient to prove in the case
$d_i=\frac12$ for all $i$.

Since $-(K_S+T)$ is an ample divisor, then replacing $T$ by some divisor we can assume that 
$\LCS(S,T)\cap U=\Gamma$. Hence, connectedness lemma implies that $\LCS(S,T)=\Gamma$.
\par
The toric projective surface $S$ (with Picard number $\rho(S)=1$) is determined by the fan $\Delta$ in the lattice $N\cong \ZZ^2$, where
$$
\Delta=\big\{\langle n_1, n_2\rangle,\langle n_2, n_3\rangle, \langle n_1, n_3\rangle, \text{their faces} \big\}.
$$
Thus surface $S$ has at most three singular points. 
If the number of singularities is less than or equal to two, then there exists an isomorphism of the lattice
$N$ such that
$n_1=(1,0)$, $n_2=(0,1)$, and therefore $S\cong \PP_{x_1x_2x_3}(a_1,a_2,1)$.

Suppose that the point $\Gamma$ is a non-toric subvariety of $(S,D')$, where $D'=D-\frac12D_j=\sum_{i\ne j}\frac12D_i$. Then the divisor $D$ can be replaced
by the other divisor $D'<D$. Therefore we have the four possibilities for the pair $(S, D)$ and the point $\G$.
\par {\bf A)} $S$ has three singular points and $D=0$. In this possibility $\G\not \in \Supp(\Sing S)$.
\par {\bf B)} $\Gamma\not \in D_{i_1}\cup D_{i_2}$, where $i_1\ne i_2$. To be definite, let $D_{i_1}-D_{i_2}$ be a nef divisor. 
\par {\bf C)} $S$ has two singular points, that is, $S\cong\PP(a_1,a_2,1)$, where $a_1\geq 3$, $a_2\geq 2$ and $\G=(b:1:0)$, where $b\ne 0$.
\par {\bf D)} $S\cong \PP(a_1,a_2,1)$, $D=\frac12\{x_1=0\}+\frac12\{x_2=0\}$, $a_1\geq 2$, $a_2\geq 1$ and $\G=(1:0:b)$, where $b\ne 0$.

\par
Possibility {\bf B)} is impossible since $\LCS(S,T-\frac12D_{i_1}+\frac12D_{i_2})=\G\cup D_{i_2}$, that is, we have the contradiction with connectedness lemma. 
Possibility {\bf D)} is impossible since $\LCS(S,T-\frac12\{x_1=0\}+\{x_3=0\})=\G\cup \{x_3=0\}$, that is, we have the contradiction with connectedness lemma.
Consider possibility {\bf C)}. Write $T=a\{x_3=0\}+T'$, where $\{x_3=0\}\not\subset\Supp(T')$ and $0\leq a<1$. The divisor $K_S+\{x_3=0\}+T'$ is
not log canonical at the point $\G$, therefore
by the inversion of adjunction we have $\big(\{x_3=0\}\cdot T'\big)_{\G}>1$. We obtain the contradiction
$$
1<\big(\{x_3=0\}\cdot T'\big)_{\G}<\{x_3=0\}\cdot (-K_S)=\frac{a_1+a_2+1}{a_1a_2}\leq 1.
$$
Consider possibility {\bf A)}. Let $f\colon (Y,E) \to (S \ni \G)$ be an inductive blow-up of $(S,T)$ (see Theorem \ref{inductive}). 
By Theorem \ref{dim2plt} the morphism $f$ is a weighted blow-up of smooth point with weights $(\alpha_1,\alpha_2)$. 
Write $K_Y+E+T_Y=f^*(K_S+T)$.

\begin{lemma} \label{compl1}
The divisor $K_S$ has a $1$-complement $B^{+}$ such that $\G$ is a center of $\LCS(S,B^{+})$. 
\begin{proof}
The divisor $K_Y+E+(1-\delta)T_Y$ is plt and anti-ample for $0<\delta\ll 1$.
Since $\rho(Y)=2$ then the cone $\NE(Y)$ is degenerated by two extremal rays. By $R_1$ and $R_2$ denote these two rays. To be definite, let $R_1$ gives
the contraction $f$.
If $-(K_Y+E)$ is a nef divisor  then a $1$--complement of 
$K_E+\Diff_E(0)=K_E+\frac{\alpha_1-1}{\alpha_1}P_1+\frac{\alpha_2-1}{\alpha_2}P_2$ is extended to a $1$--complement of $K_Y+E$
by Proposition \ref{complind}, therefore we obtain the required $1$--complement of $K_S$
by Proposition \ref{compldown}. 

Consider the last possibility: $(K_Y+E)\cdot R_2>0$, $T_Y\cdot R_2<0$.
Let $L(\delta)\in |-n(K_Y+E+(1-\delta)T_Y)|$ be a general element for $n\gg 0$ and let $M=(1-\delta)T_Y+\frac1nL(\delta)$, where $\delta>0$ is a sufficiently
small fixed rational number. By construction, $K_Y+E+(1+\varepsilon)M\equiv \varepsilon M$, $K_Y+E+(1+\varepsilon)M$ is a plt divisor. Therefore,
applying $(K_Y+E+(1+\varepsilon)M)$--MMP is a contraction of the ray 
$R_2$ for $0<\varepsilon\ll 1$. The corresponding divisorial contraction is denoted by $h\colon Y\to \overline S$, and the image of
$E$ on the surface $\overline S$ is denoted by $\overline E$, put $\Exc h=C_Y$ and $C_S=f(C_Y)$. The divisor $K_{\overline S}+\overline E$ is plt and anti-ample.
Therefore, if $1$--complement of $K_{\overline E}+\Diff_{\overline E}(0)$ exists then we consistently apply Theorems \ref{complind}, \ref{complup} and
\ref{compldown} and obtain the required $1$--complement of $K_S$. 

Suppose that there does not exist any $1$--complement of $K_{\overline E}+\Diff_{\overline E}(0)$. It is possible if and only if
there are three singular points of $\overline S$ lying on the curve $\overline E$. It implies that $\alpha_1\geq 2$, $\alpha_2\geq 2$, the curve
$C_Y$ is contracted to a cyclic singularity, and the curve $C_S$ passes through at most one singularity of $S$ (see \cite[Chapter 3]{Koetal}).
Let us apply Corollary 9.2 of the paper \cite{KeM} for $K_{\overline S}+\overline E$. We obtain that $S$ has the two singularities of type $\AAA_1$, which 
do not lie on the curve $C_S$. 
Let $V(\langle n_1\rangle)$ be the closure of one-dimensional orbit passing through the two singular points of type $\AAA_1$. 
Then there exists an isomorphism of the lattice 
$N$ such that $n_1=(1,0)$, $n_2=(1,2)$, and therefore $n_3=(-2n+1,-2)$, where $n\ge 2$.
By considering the cone $\langle n_2, n_3 \rangle$ we obtain that the third singularity of
$S$ is of type $\frac1{4n-4}(2n-1,1)$, its minimal resolution graph consists of three exceptional curve chain with the self-intersection indices
$-2$, $-n$ and $-2$ respectively. The following two cases are possible: i) $\G \in V(\langle n_2\rangle)\cup V(\langle n_3\rangle)$ and ii)
$\G \not\in V(\langle n_2\rangle)\cup V(\langle n_3\rangle)$. 

Consider former Case i). To be definite, let  $\G \in V(\langle n_2\rangle)$, then
$V(\langle n_2\rangle)\cdot (-K_S)=\frac{n}{2n-2}\leq 1$, and therefore we obtain a contradiction for the same reason as in Case {\bf C)}. 

Consider latter Case ii).
Let $g\colon S^{\min}\to S$ be a minimal resolution. Let us contract all curves of $\Exc g$, except the exceptional curve of the singularity
$\frac1{4n-4}(2n-1,1)$ with the self-intersection index
$-n$. We obtain the divisorial contractions $S^{\min}\to \widetilde S$
and $\widetilde S\to S$. Note that $\rho(\widetilde S)=2$ and $\widetilde S=T_N(\widetilde \Delta)$, where the fan $\widetilde \Delta$ is given by $\Delta$ 
with the help of subdivision of the cone 
$\langle n_2, n_3 \rangle$
into the two cones $\langle n_2, n_4 \rangle$, $\langle n_4, n_3 \rangle$, where $n_4=(-1,0)$. The surface $\widetilde S$ is a conic bundle with
irreducible fibers, and its two fibers are non-reduced. These two fibers are the curves 
$V(\langle n_2\rangle)$, $V(\langle n_3\rangle)$, and every such curve contains the
two singularities of type $\AAA_1$. By $\widetilde\G$ denote the transform of $\G$ on the surface $\widetilde S$. We have 
$K_{\widetilde S}+\widetilde B^+_1+\widetilde B^+_2+V(\langle n_4\rangle)\sim 0$, where $\widetilde B^+_1\sim V(\langle n_2\rangle)+V(\langle n_3\rangle)$ is the
fiber passing through the point
$\widetilde\G$, and 
$\widetilde B^+_2\sim V(\langle n_1\rangle)$ is the section passing through the point $\widetilde\G$. By Proposition \ref{compldown} 
we obtain the required $1$--complement of $K_S$.
\end{proof}
\end{lemma}

Assume that $B^+=B_1^+ + B^{+'}$, where the irreducible curve $B_1^+$ has an ordinary double point singularity at the point $\G$.
By the inversion of adjunction we have $B^{+'}=0$, $B_1^+\cap \Supp(\Sing S)=\emptyset$ and $K_S+B_1^+\sim 0$, therefore $K_S$ is Cartier divisor.
Classification of Del Pezzo surfaces with Du Val singularities (in our case Du Val singularities are cyclic), with Picard number 1 and with 
three singular points implies $K^2_S\leq 4$ \cite{Fur}.
Write $T=aB_1^++T'$, where $B_1^+\not\subset\Supp(T')$ and $0\leq a<1$. Since $0\sim K_Y+E+\widetilde{B_1^+}=f^*(K_S+B_1^+)$ then we obtain the contradiction
\begin{gather*}
0>(K_Y+E+T_Y)\cdot \widetilde{B_1^+}\geq (-1+a)\Big(\widetilde{B_1^+}\Big)^2=\\ =(-1+a)\Big(K_S^2-\frac{(\alpha_1+\alpha_2)^2}{\alpha_1\alpha_2}\Big)\geq 0.
\end{gather*}

Consider the last case $B^+=B_1^+ +B_2^+ +B^{+'}$, where the irreducible curves $B_1^+$ and $B_2^+$ have a simple normal crossing at the point $\G$.
We have $(B_1^+\cup B_2^+)\supset \Supp(\Sing S)$ according to Corollary 9.2 of the paper \cite{KeM} applied for $K_{S}+B_1^+ +B_2^+$. To be definite,
let the curve $B_1^+$ contains two singular points of
$S$. By the inversion of adjunction, $\deg\Diff_{B_1^+}(0)\leq 1$, and therefore the curve $B_1^+$ passes through two singular points only, and they are of type
$\AAA_1$. Such surfaces were classified in the proof of Lemma \ref{compl1}, and therefore it can be assumed that the third singularity of $S$ is of 
type $\frac1{4n-4}(2n-1,1)$, 
$B^{+'}=0$, $B^+_1\cap B^+_2=\G$, $(B_1^+)^2=n-1$ and $(B_2^+)^2=\frac1{n-1}$, where $n\geq 2$.
To be definite, assume that $f^*(B_1^+)=\widetilde{B_1^+}+\alpha_1E$ and
$f^*(B_2^+)=\widetilde{B_2^+}+\alpha_2E$. Thus $(\widetilde{B_1^+})^2=n-1-\alpha_1/\alpha_2$, $(\widetilde{B_2^+})^2=\frac1{n-1}-\alpha_2/\alpha_1$, and therefore
$(\widetilde{B_k^+})^2\leq 0$ for either $k=1$ or $k=2$. 
Write $T=a_1B_1^++a_2B_2^++T'$, where $B_1^+, B_2^+\not\subset\Supp(T')$, $0\leq a_1< 1$, $0\leq a_2<1$.
Since $0\sim K_Y+E+\widetilde{B_1^+}+\widetilde{B_2^+}=f^*(K_S+B_1^++B_2^+)$, then we obtain the contradiction
\begin{gather*}
0>(K_Y+E+T_Y)\cdot \widetilde{B_k^+}= (-1+a_k)\Big(\widetilde{B_k^+}\Big)^2+T'_Y\cdot\widetilde{B_k^+}\geq\\ \geq (-1+a_k)\Big(\widetilde{B_k^+}\Big)^2\geq 0.
\end{gather*}

Consider Condition 2). Such toric surface is determined by the fan $\Delta$ in the lattice $N\cong \ZZ^2$, where
$$
\Delta=\big\{\langle m_1, m_2\rangle,\langle m_2, m_3\rangle, \langle m_3, m_4,\rangle, \langle m_4, m_1\rangle, \text{their faces} \big\},
$$
$m_1=(1,0)$, $m_2=(q,r)$, $m_3=(-1,0)$, $m_4=(-q,-r)$, $q\ge 1$, $r\ge 1$ and $\gcd(q,r)=1$.
Therefore $S$ has four singularities of types $\frac1r(1,-q)$, $\frac1r(1,q)$, $\frac1r(1,-q)$ and $\frac1r(1,q)$ respectively. 

Two different fibers passing through the point $\G$ are denoted by $F_1$ and $F_2$. Since $T\cdot F_i\ge 1$ by Lemma \ref{adjnlc} for $i=1,2$, then
$T-F_1-F_2$ is nef. 

\begin{lemma}\label{adjnlc}
Let $O$ be a smooth point of the surface $M$. Assume $(M,N)$ is not a log canonical pair at the point $O$, where
$N=dI+\Sigma\geq 0$, $I\not\subset\Supp\Sigma$, $d\leq 1$, $I$ is an irreducible curve which is a smooth at the point $O$. Then $(\Sigma\cdot I)_O>1$.
\begin{proof}
The proof follows by the inversion of adjunction, see, for example, \cite[Theorem 7.5]{Kollar}.
\end{proof}
\end{lemma}

Consider the index $j$ such that $F_j$ is non-toric subvariety of $(S,D)$. 
Let $F'$ and $F''$ be the closures of one-dimensional toric orbits provided that
$F'\sim_{\QQ} F''\sim_{\QQ}\frac1r F_j$.
We obtain the contradiction
$(K_S+T)\cdot F_i\ge(-F'-F''+D+F_j)\cdot F_i\geq 0$, where the index $i\in\{1,2\}$ satisfies the condition $i\ne j$.

Consider Condition 3). Assume that either $\G\in C_{13}$, or $\G\in C_{23}$. Let us consider the first possibility. The second possibility is considered similarly.
If $\G$ is a non-toric point of $(C_{13},\Diff_{C_{13}}(D))$ then we have a contradiction with one-dimensional analog of this theorem since $C_{13}^2<0$.
Therefore, $a_4=1$ and $d_{23}=1$. The case  
$C^2_{23}=0$ is impossible also (in this case the surface $S$ is a toric conic bundle, and we use its structure).
Thus $C_{23}^2<0$ and consider the contraction $\psi\colon S\to S'$ of $C_{23}$. If $\psi(\G)$ is a non-toric point of $(S',\psi(D))$ then
we have a contradiction with this theorem under Condition 1). Therefore the curve $C_{23}$ is contracted to a smooth point and $d_{24}=1$. We obtain the contradiction
$a_2d_{13}d_{14}=w_2>w_4=d_{13}$. 

Assume that $\G\in C_{13}\cup C_{23}$. Let $C_{23}^2=0$.
Then $(w_3,w_4)=(w_1,w_2)$, $(S,D)\cong\big(\FFF_{w_1-w_2},\frac{w_2-1}{w_2}C_{13}+\frac{w_1-1}{w_1}C_{24}\big)$ and $2\le w_1<w_2$.
By $F_j$ denote a fiber of
$S$ passing through $\G$. Then $T'\cdot F_j\ge 1$ by Lemma \ref{adjnlc}, 
where $T=T'+\alpha F_j$, $F_j\not\subset \Supp(T')$, and we have the contradiction
$0>(K_S+T)\cdot F_j\ge (K_S+F_j+T')\cdot F_j\ge 0$.

Therefore $C_{23}^2<0$, $C_{13}^2<0$. 
Considering case by case the contractions of the curves
$C_{13}$ and $C_{23}$, we obtain that these curves are contracted to smooth points and 
$d_{13}=d_{23}=d_{24}=d_{14}=1$.
Since $C_{13}^2=-\frac1{a_2a_4}$, $C_{23}^2=-\frac1{a_1a_4}$ then $(w_1,w_2,w_3,w_4)=(a_2,a_2,a_2-1,a_2+1)$, $a_2\ge 3$.
It is easy to find a birational map
$$
S\dashrightarrow S'\bigg(\frac1{a_2-1}(1,-1)+\frac1{a_2-1}(1,1)+\frac1{a_2-1}(1,-1)+\frac1{a_2-1}(1,1)\bigg),
$$
where $\rho(S')=2$, and in result of this map we obtain a contradiction with this theorem under Condition 2).
To find this map it is enough to consider two (required) toric blow-ups at the points $P_2$, $P_4$ and a contraction of proper transforms of $C_{13}$ and $C_{23}$.
\end{proof}
\end{theorem}

\begin{remark} 
Theorem \ref{point} 1) can not be generalized to the case $\rho(S)\ge 2$. 
Consider the toric pair$(S,D)=(\FFF_1,\frac12E_0)$ and the divisor $T=\frac12E_0+E'_0+F+\delta E_{\infty}$ provided that $F\cap E'_0\not\in E_0\cup E_{\infty}$, where  
$E_0$, $E'_0$ are two different zero sections, $E_{\infty}$ is the infinity section, $F$ is a fiber and $0<\delta<\frac12$.
Put 
$\G=F\cap E'_0$. Then $\G$ is a non-toric point of $(S,D)$, $T\ge D$, $K_S+T$ is anti-ample log canonical divisor and $\G\in\LCS(S,T)$.

Nevertheless, it is expected that Theorem 
\ref{point} can be generalized to every dimension and every Picard number $\rho(S)$, if we require the following condition, instead of Conditions 1), 2) and 3):
$(S,D)=(E,\Diff_E(0))$, where $f\colon (Y,E)\to (X\ni P)$ is a toric plt blow-up of some toric singularity.
\end{remark}

\begin{definition}
Let $(\G,D_{\G})\cong(\PP^1,\sum_{i=1}^r\frac{m_i-1}{m_i}P_i)$. Assume that $-(K_{\G}+D_{\G})$ is an ample divisor. Then, for set
$(m_1,\ldots,m_r)$ we have one of the following cases up to permutations:
$(m_1,m_2)$, it is of type $A$; $(2,2,m)$, $m\geq 2$, it is of type $D_{m+2}$; $(2,3,3)$, it is of type $E_6$; $(2,3,4)$, it is of type $E_7$;
$(2,3,5)$, it is of type $E_8$.
In Propositions \ref{defS} and \ref{defS2} the classification according to types corresponds to the types of $(\G,D_{\G})=(\G,\Diff_{\G}(D))$.
\end{definition}

\begin{proposition}\label{defS}
Let $(S,D)$ be a toric pair, where $S$ is a normal projective surface with $\rho(S)=1$, and let
$D$ be a divisor with standard coefficients. Assume that there exists a curve $\G$ such that
$-(K_S+D+\G)$ is an ample divisor and $(S,D+\G)$ is a plt non-toric pair.
Let us denote a hypersurface of degree $d$ in a weighted projective space by $X_d$.
Then one of the following cases is satisfied.
\par $1)$ $(S,D,\G)\cong (\PP^2_{x_1x_2x_3},\frac{d_1-1}{d_1}\{x_1=0\},X_2)$ and $d_1\geq 1$. It is of type $A$.
\par $2)$ $(S,D,\G)\cong (\PP^2_{x_1x_2x_3},\sum_{i=1}^3\frac{d_i-1}{d_i}\{x_i=0\},X_1)$, 
the integer number triple $(d_1,d_2,d_3)$ is either $(2,2,k)$, $(2,3,3)$, $(2,3,4)$ or $(2,3,5)$, where $k\geq 2$. 
They are of types $D_{k+2}$, $E_6$, $E_7$ and $E_8$ respectively.
\par $3)$ $(S,D,\G)\cong (\PP_{x_1x_2x_3}(a_1,1,1), \sum_{i=1}^2\frac{d_i-1}{d_i}\{x_i=0\},X_{a_1})$, 
the integer number triple $(a_1,d_1,d_2)$ is either $(2,2,k_1)$, $(2,3,k_2)$, $(2,k_3,1)$ or $(3,2,1)$, 
where $k_1\geq 1$, $1\leq k_2\leq 2$, $k_3\geq 4$. In the first possibility, if $k_1\geq 2$ then it is of type $D_{k_1+2}$. In the second possibility, if $k_2=2$ 
then it is of type $E_6$. The other possibilities are of type $A$ always. 
\par $4)$ $(S,D,\G)\cong (\PP_{x_1x_2x_3}(a_1,1,1), \frac{d_1-1}{d_1}\{x_2=0\},X_{a_1+1})$, $a_1\geq 2$ and $d_1\geq 1$. It is of type $A$.
\par $5)$ $(S,D,\G)\cong (\PP_{x_1x_2x_3}(a_2+1,a_2,1), \sum_{i=1}^2\frac{d_i-1}{d_i}\{x_i=0\},X_{a_2+1})$,
the integer number triple $(a_2,d_1,d_2)$ is either $(2,2,k_1)$, $(k_2,2,k_3)$ or $(k_4,k_5,1)$, where $k_1\leq 3$, $k_2\geq 3$, $k_3\leq 2$, $k_4\geq 2$ and
$k_5\geq 3$. In the first possibility, if $k_1=2$ then it is of type $D_6$, and, if $k_1=3$ then it is of type $E_7$. In the second possibility, 
if $k_3=2$ then it is of type $D_{2k_2+2}$.
The other possibilities are of type $A$ always.
\par $6)$ $(S,D,\G)\cong (\PP_{x_1x_2x_3}(2a_2+1,a_2,1), \frac12\{x_1=0\},X_{2a_2+1})$, $a_2\geq 2$. It is of type $D_{2a_2+2}$. 
\par $7)$ $(S,D,\G)\cong (\PP_{x_1x_2x_3}(la_2-1,a_2,1), \sum_{i=1}^2\frac{d_i-1}{d_i}\{x_i=0\},X_{la_2})$, $a_2\geq 2$,
the integer number triple $(l,d_1,d_2)$ is either $(2,2,1)$ or $(k_1,1,k_2)$, where $k_1\geq 2$ and $k_2\geq 1$. They are of types $D_{2a_2+1}$ and $A$ respectively.
\par $8)$ $(S,D,\G)\cong (\PP_{x_1x_2x_3}(a_1,a_2,1), \frac{d_1-1}{d_1}\{x_3=0\},X_{a_1+a_2})$, $a_1>a_2\geq 2$ and $d_1\geq 1$. It is of type $A$.
\par $9)$ $(S,D)\cong(S(\frac1{r_1}(1,1)+\frac1{r_2}(1,1)+\AAA_{r_1+r_2-1}),\frac{d_1-1}{d_1}D_3)$, $\G\sim_{\QQ}D_3$ is an irreducible curve being different from $D_3$, 
where $D_3$ is the closure of one-dimensional orbit passing through the first and second singular points,
$d_1\geq 2$ and $r_1$, $r_2\geq 2$. It is of type $A$.
\par $10)$ $(S,D)\cong(S(\frac1{r_1}(l,1)+\frac1{r_2}(l,1)+\AAA_{(r_1+r_2)/l-1}),\frac{d_1-1}{d_1}D_3)$, the surface $S$ has three singular points,
$\G\sim D_1+D_2$,  
where $D_i$ is the closure of one-dimensional orbit not passing through the $i$-th singular point 
of $S$, $d_1\geq 1$, $l\geq 2$ and $l|(r_1+r_2)$. It is of type $A$.
\begin{proof}
By the adjunction formula the curve $\G$ is smooth and irreducible.
It follows easily that, if $P\in \Supp D\cap\G$ then $(S,D+\G)$ is a toric pair in a sufficiently small analytical neighborhood of $P$. 
If $S$ is a smooth surface then $S\cong \PP^2$ and we have two Cases 1) and 2). 
\par
Assume that $S$ is a non-smooth surface having at most two singular points. Then we have
$S\cong \PP_{x_1x_2x_3}(a_1,a_2,1)$ as before in the
proof of Theorem \ref{point}. At first let us consider the case of one singular point, that is, $a_1\geq 2$ and $a_2=1$. 
Then either $\G\sim \OO_S(1)$, $\OO_S(a_1)$ or $\OO_S(a_1+1)$.
The variant $\G\sim \OO_S(1)$ is impossible since $K_S+D+\G$ is not a plt divisor at the point $(1:0:0)$.
The other variants lead us to Cases 3) and 4) respectively. 
At second let us consider the case of two singular points, that is, $a_1>a_2\geq 2$. Put $\G=\{\psi(x_1,x_2,x_3)=0\}$. Suppose that 
$\G\not\sim\OO_S(a_1+a_2)$, $\OO_S(a_1)$, $\OO_S(a_2)$, $\OO_S(1)$ then
$\psi(x_1,x_2,x_3)=bx_1x_3^l+\varphi(x_2,x_3)$, and by considering the point $(1:0:0)$ we obtain $b\ne 0$, $l=1$, $\G\sim\OO_S(a_1+1)$ and 
$x_2^m\in \varphi(x_2,x_3)$. It leads us to Case 7). If $\G\sim\OO_S(a_1)$ then by considering the point $(0:1:0)$ we obtain $x_1,x_2^lx_3\in\psi(x_1,x_2,x_3)$. 
It leads us to Cases 5) and 6).
It is easy to prove that cases $\G\sim\OO_S(a_2)$ and $\G\sim\OO_S(1)$ are not realized.
If $\G\sim\OO_S(a_1+a_2)$ then $x_1x_2,x_3^{a_1+a_2}\in\psi(x_1,x_2,x_3)$, and we have Case 8). 
\par
Assume that $S$ is a surface having three singular points (it is the last possibility for $S$).
According to Corollary 9.2 of the paper \cite{KeM} for the divisor $K_S+\G$, we obtain that the curve $\G$ contains a singular point of $S$.

Suppose that the curve $\G$ contains only one singular point of $S$, then arguing as above in the proof of Theorem
\ref{point}, we obtain $S=S(2\AAA_1+\frac1{4n-4}(2n-1,1))$, where $n\geq 2$, and $\G$ is locally a toric subvariety of
$(S\ni P)$, where $(S\ni P)$ is of type $\frac1{4n-4}(2n-1,1)$. By $T_1$ and $T_2$ denote the closures of one-dimensional orbits passing through the singular point $P$. 
Since 
$T_1\sim T_2$ and $(\G\cdot T_1)_P\ne (\G\cdot T_2)_P$ then $\G\cdot T_i>1$. Therefore $\G-(4n-4)T_1$ is an ample divisor, and we obtain the contradiction
with ampleness of  
$-(K_S+\G)\sim 2nT_1-\G$. Thus this possibility is not realized. 

Suppose that the curve $\G$ passes through the two singular points $P_1$ and $P_2$
of $S$ only.
There exists a 1--complement of $K_{\G}+\Diff_{\G}(0)$, and we obtain the 1-complement $K_S+\G+T\sim 0$ of $K_S+\G$
by Proposition \ref{complind}.
There are two Cases {\bf A)} and {\bf B)}.
\par {\bf A)} Let $T$ is a reducible divisor. By the two-dimensional criterion on the characterization
of toric varieties
\cite[Theorem 6.4]{Sh2} we have $T=T_1+T_2$, $\G\sim T_3$, $D=\frac{d_1-1}{d_1}T_3$, 
the singularities at the points $P_j$ are of type $\frac1{r_j}(1,1)$,  
where $d_1\geq 2$, $r_j\geq 2$ and $T_i$ are the closures of one-dimensional orbits, and $P_1\in T_1$.
Let $f\colon\widetilde S\to S$ be a minimal resolution at the points $P_1$ and $P_2$ only. By $E_1$ denote the curve such that $f(E_1)=P_1$.
By the inversion of adjunction $\G\cdot T_3=\frac1{r_1}+\frac1{r_2}$, hence
$(\G_{\widetilde S})^2=\G_{\widetilde S}\cdot(T_3)_{\widetilde S}=0$, and the linear system $|E_1+m\G_{\widetilde S}|$ gives the birational morphism
$g\colon\widetilde S\to \FFF_{r_1}$ for $m\gg 0$ \cite[Proposition 1.10]{KudLd} such that the curve $(T_2)_{\widetilde S}$ is contracted to a smooth point.
The morphism $g$ is toric and the third singularity of
$S$ is of type $\AAA_{r_1+r_2-1}$. We obtain Case 9).
\par {\bf B)} Let $T$ is an irreducible divisor. To be definite, let $D_i$ be the closures of one-dimensional orbits not passing through the $i$-th singular point
of $S=S(\frac1{r_1}(a_1,1)+\frac1{r_2}(a_2,1)+\frac1{r_3}(a_3,1))$. We have $\frac1{r_1}D_1\equiv\frac1{r_2}D_2\equiv\frac1{r_3}D_3$. 
To be definite, the curve $\G$ passes through the first and second singular point of $S$. 
By the definition of
1--complement we obtain $\G\cdot T=\frac1{r_1}+\frac1{r_2}$, $\G+T\sim \sum_{i=1}^3D_i$. Hence, either $\G\sim D_1+D_2$, $T\sim D_3$ or 
$\G\sim D_3$, $T\sim D_1+D_2$. Since 1--complement not passing through the third singular point of $S$ then it is of type $\AAA_{r_3-1}$.
The case $\G\sim D_3$ was considered in Case A). Since the curve $\G$ does not pass through the third singular point then we have to consider 
the possibility remained:
$\G\sim D_1+D_2\sim lD_3$, where $l\geq 2$, $l\in\ZZ$. We obtain Case 10). 
\par
Suppose that the curve $\G$ passes through three singular points of $S$ with the indices $r_1$, $r_2$ and $r_3$ respectively.
By the inversion of adjunction the triple $(r_1,r_2,r_3)$ is either $(2,2,k)$, $(2,3,3)$, $(2,3,4)$ or $(2,3,5)$, where $k\geq 2$. 
For the second and third variants there does not exist any surface 
$S$.
For the first and fourth variants we have
$S=S(2\AAA_1+\frac1{4n-4}(2n-1,1))$ and $S\cong\PP(2,3,5)$ respectively, where $n\geq 2$. These variants are considered
as above mentioned case, when the curve $\G$ contains only one singular point of $S$. 
\end{proof}
\end{proposition} 

\begin{proposition}\label{defS2}
Let $(S,D)$ be ODP pair. 
Assume that there exist a curve $\G$ and an effective $\QQ$-divisor $\G'$ such that
$K_S+D+\G+\G'$ is an anti-ample and plt divisor, and $(S,D+\G)$ is a non-toric pair.
Then $d_{23}=d_{24}=1$, $a_1|a_2$ and $\G\sim\OO_{\PP({\bf w})}(w_2)|_S$ up to permutation of the coordinates.
In particular, $-(K_S+D+\G)$ is an ample divisor and $w_1|w_2$. It is of type $A$.
\begin{proof}
The sets $\G\cap C_{13}$, $\G\cap C_{23}$ consist of at most one point by the adjunction formula. Moreover, we may assume that $\G'=\gamma_1C_{13}+\gamma_2C_{23}$, where
$\gamma_1<1$ and $\gamma_2<1$. If $C_{i3}^2=0$ then $\gamma_i=0$, where $i=1,2$. 

Let us prove that $\G\cdot C_{13}>0$ and $\G\cdot C_{23}>0$. 
Assuming the converse: $\G\cdot C_{13}=0$, that is, $\G\sim dC_{24}$. The possibility $\G\cdot C_{23}=0$ is considered similarly. 
Since $C_{23}\cdot C_{24}=\frac1{a_1}$, $a_1(C_{23}\cdot\G)\in\ZZ_{>0}$ then $d\in\ZZ_{>0}$. 
The divisor $C_{24}-\gamma C_{13}$ is nef for $0\leq\gamma\leq \frac1{d_{13}}$, hence it is semiample by the base point free theorem \cite{KMM}.
Therefore, if $d\geq 2$ then we have a contradiction with connectedness lemma, since there exists a $\QQ$--divisor $\G''$ such that $\down{\G''}=0$ and
$D+\G+\G'\sim_{\QQ} C_{24}+C_{13}+\G''$. Thus, $d=1$. Since the curve $\G$ is a non-toric subvariety of $(S,D)$ then
$d_{24}\geq 2$, and we have $d_{13}=1$ by connectedness lemma again. We obtain the contradiction
\begin{gather*}
0>(K_S+D+\G+\G')\cdot C_{23}\geq \\ \geq\Big(\frac{d_{24}-1}{d_{24}}C_{24}-C_{13}-C_{23}-C_{14}+\G'\Big)\cdot C_{23}\geq
\\
\geq\frac{d_{24}-1}{d_{24}}C_{24}\cdot C_{23}-C_{13}\cdot C_{23}=d_{23}\Big(\frac{d_{24}-1}{w_1}-\frac1{w_4}\Big)\geq 0.
\end{gather*}

Thus, we proved that the sets $\G\cap C_{13}$ and $\G\cap C_{23}$ consist of one point only.

Suppose that $P_4\not\in\G$. Then $\G\sim_{\QQ} \alpha_1C_{14}+\alpha_2C_{24}$, $\alpha_1=a_2(\G\cdot C_{13})\in \ZZ_{>0}$ and
$\alpha_2=a_1(\G\cdot C_{23})\in \ZZ_{>0}$. By applying connectedness lemma we have $\alpha_1=\alpha_2=1$. Let us prove that $d_{14}=d_{24}=1$.
Assuming the converse: $d_{14}\geq 2$. The possibility $d_{24}\geq 2$ is considered similarly. In order to apply connectedness lemma and obtain
a contradiction (for the disjoint curves $C_{14}$, $C_{23}$) we must only prove that
$D_1=\frac{d_{14}-1}{d_{14}}C_{14}+C_{24}+\frac{d_{24}-1}{d_{24}}C_{24}-\frac1{d_{23}}C_{23}$ 
is a semiample divisor. Since $D_1\cdot C_{23}>0$ and $D_1\cdot C_{13}=d_{13}(\frac{d_{14}-1}{w_2}-\frac1{w_4})\geq 0$ then
$D_1$ is a nef divisor and it is semiample by the base point free theorem \cite{KMM}. 
Finally, since $K_S+\G+C_{13}+C_{23}\sim 0$ then $K_S$ is Cartier divisor at the point $P_3$, and the singularity at the point $P_3$
is Du Val of type $\frac1{w_3}(w_1, w_2)$.
Therefore $w_3+w_4=w_1+w_2\equiv 0 (\mt{mod} w_3)$, $w_3|w_4$ and $a_3|a_4$. 

Suppose that $P_4\in\G$. Since the curve $\G$ is a (locally) toric orbit in some analytical neighborhood of $P_4$ then either 
$\G\cdot C_{13}=\frac1{a_4}$ or $\G\cdot C_{23}=\frac1{a_4}$. Let us consider the former case. The latter case is considered similarly.
Write $\G\sim_{\QQ} \alpha_1C_{23}+\alpha_2C_{24}$, $\alpha_1=a_4(\G\cdot C_{13})=1$ and
$\alpha_2=a_3(\G\cdot C_{14})\in \ZZ_{>0}$. Arguing as above, we see that $\alpha_2=1$, $d_{24}=1$. If $d_{23}=1$ then this proposition is proved. 
Let $d_{23}\geq 2$. By the plt assumption of this proposition $\G\cdot C_{23}=\frac1{a_4}$ and $d_{13}=1$.
Considering $\G\sim_{\QQ} C_{13}+\alpha_2'C_{14}$ we obtain $\alpha_2'=1$, $d_{14}=1$. This completes the proof.
\end{proof}
\end{proposition}

\begin{definition}\label{triple}
The triple $(S,D,\G)$ determined by the assertions of Propositions \ref{defS} or \ref{defS2} is said to be a {\it purely log terminal triple}.
\end{definition}

The following problem is important for the classification of plt blow-ups of three-dimensional toric non-$\QQ$-factorial
singularity
(if we follow the method described in this paper).
\begin{problem}
Let $(S,D)=(E,\Diff_E(0))$, where $f\colon (Y,E)\to (X\ni P)$ is a toric plt blow-up of some toric three-dimensional 
(non-$\QQ$-factorial) singularity. Assume that
there exist a curve $\G$ and an effective $\QQ$-divisor $\G'$ such that
$K_S+D+\G+\G'$ is an anti-ample plt divisor, and $\G$ is a non-toric subvariety of $(S,D)$.
Classify the triples $(S,D,\G)$.
\end{problem}

\section{\bf Non-toric three-dimensional blow-ups. Case of point}\label{sectionex}
\begin{example}
Now we construct the examples of three-dimensional non-toric plt blow-ups $f\colon (Y,E)\to (X\ni P)$ provided that $(X\ni P)$ is a
$\QQ$-gorenstein toric singularity and
$P=f(E)$. Depending on a type of $(X\ni P)$ we consider two Cases {\bf A1)} and {\bf A2)}.

\par {\bf A1).} Let $(X\ni P)$ be a $\QQ$-factorial toric singularity, that is, $(X\ni P)\cong (\CC^3\ni 0)/G$, 
where $G$ is an abelian group acting freely in codimension 1.
All plt blow-ups are constructed by the procedure illustrated on the next diagram (Fig. 4) and defined below.

\[
\xymatrix{& Y_0 \ar[dl]_{h_0} \ar[d]^{h'_0} & Y_1 \ar[l]_{h_1} \ar@{-->}[d]^{h'_1}\\
	Z_0 \ar[d]_{g_0}& Z_1 \ar[d]_{g_1} & Z_2 \ar[d]_{g_2}\\
	X& X& X
}
\]\begin{center}
	Fig. $4$
\end{center}

{\it First step.}
Let $g_0\colon (Z_0,S_0)\to (X\ni P)$ be a toric blow-up, where $\Exc g_0=S_0$ is an irreducible divisor
and $g_0(S_0)=P$.
Assume that there exists a curve $\G_0\subset S_0$ such that $(S_0,\Diff_{S_0}(0),\G_0)$ is a plt triple (see Definition \ref{triple}).
Such triples are classified in Proposition \ref{defS} and are divided into the five types: $A$, $D_l$, $E_6$, $E_7$ and $E_8$.

\begin{remark}\label{rem1}
There exists an irreducible reduced Weil divisor $\Omega$ on $X$ such that $\Omega_{Z_0}|_{S_0}=\G_0$.
The surface $\Omega$ has a log terminal singularity at the point $P$. A singularity type coincides with a type of the triple $(S_0,\Diff_{S_0}(0),\G_0)$. 
In particular, if $\psi$ is a 
$G$--semi-invariant polynomial in $\CC^3$ determining $\Omega$ then Du Val singularity $\{\psi=0\}\subset(\CC^3 \ni 0)$ is of the same type.
\end{remark}

The following lemma gives a restriction on the triple $(S_0,\Diff_{S_0}(0),\G_0)$ in the case of terminal singularities.
\begin{lemma}\label{termrestr}
Let $(X\ni P)$ be a terminal singularity, that is, it is of type $\frac1r(-1,-q,1)$, where $\gcd(r,q)=1$ and $1\leq q\leq r$. 
Write $\Diff_{S_0}(0)=\sum_{i=1}^3\frac{d_i-1}{d_i}D_i$, where $D_i$ are the closures of corresponding one-dimensional orbits of the toric surface $S_0$.
Then $\gcd(d_i,d_j)=1$ for $i\ne j$.
\begin{proof}
It is sufficient to prove that the singularities of $Z_0$ are cyclic.
Consider the cone $\sigma$ determining the singularity $(X\ni P)$ (see Example \ref{ex1} 1)). By $(w_1,w_2,w_3)$ denote the primitive vector defining the blow-up $g_0$.
Then $Z_0$ is covered by three affine charts with the singularities of types $\frac1{w_3}(-w_1,-w_2,1)$, $\frac1{rw_2-qw_3}(-w_1+uw_2+vw_3,-uw_2-vw_3,1)$ and
$\frac1{rw_1-w_3}(-w_1,qw_1-w_2,1)$, where $uq+vr=1$ and $u$, $v\in\ZZ$. 
\end{proof}
\end{lemma}

According to Proposition \ref{defS} 
the curve $\G_0$ is locally a toric subvariety of $Z_0$ in every sufficiently small analytic neighborhood of each point of $\G_0$.
Note also that $Z_0$ is a smooth variety at a general point of
$\G_0$. 

Let $h_0\colon (Y_0,\widetilde S_1)\to (Z_0\supset \G_0)$ be an arbitrary blow-up of the
curve $\G_0$ with an unique exceptional divisor ($\Exc h_0=\widetilde S_1$) for which the following three conditions are satisfied.

1) The morphism $h_0$ is locally toric at every point of $\G_0$. In particular, $\widetilde S_0\cong S_0$, $\rho(\widetilde S_0)=1$. 

2) Let $H_0$ be a general hyperplane section of $Z_0$ passing through the general point $Q_0\in\G_0$. Then the morphism $h_0$
induces a weighted blow-up of the smooth point $(H_0\ni Q_0)$ with weights $(\beta^1_0,\beta^2_0)$.

3) $h_0^*S_0=\widetilde S_0+\beta^2_0\widetilde S_1$. 

The set of all possible blow-ups $h_0$ is denoted by $\mathcal H_0$.
The morphism $h'_0$ gives the divisorial contraction $h'_0\colon Y_0\to Z_1$ which contracts the divisor $\widetilde S_0$ to a point.
We obtain a non-toric blow-up
$g_1\colon (Z_1,S_1)\to (X\ni P)$, where $\Exc g_1=S_1$ is an irreducible divisor and $g_1(S_1)=P$.

\begin{lemma}\label{lemm1}
Let $\widetilde\Gamma_0=\widetilde S_0\cap\widetilde S_1$. Then
$$
(\widetilde\Gamma_0^2)_{\widetilde S_1}=\beta^1_0\frac{\big(K_{S_0}+\Diff_{S_0}(0)\big)\cdot\Gamma_0}{a(S_0,0)+1}
- \beta^2_0(\Gamma_0^2)_{S_0}.
$$
\begin{proof}
This formula follows from the following equalities
\begin{gather*}
(\widetilde\Gamma_0^2)_{\widetilde S_1}=\beta^1_0\widetilde S_0\cdot\widetilde\Gamma_0=\beta^1_0(S_0\cdot\Gamma_0-
\beta^2_0\widetilde S_1\cdot \widetilde \Gamma_0)= \beta^1_0 S_0\cdot\Gamma_0 -\\-
\beta^2_0(\widetilde\Gamma_0^2)_{\widetilde S_0}=\beta^1_0 S_0\cdot\Gamma_0  -\beta^2_0(\Gamma_0^2)_{S_0}
=\\= \beta^1_0((K_{Z_0}+S_0)\cdot \G_0)/(a(S_0,0)+1) - \beta^2_0(\Gamma_0^2)_{S_0}.
\end{gather*}
\end{proof}
\end{lemma}

In next Proposition \ref{pltEDiff} we will describe the pair $(S_1,\Diff_{S_1}(0))$.
The surface $\widetilde S_1$ is a conic bundle with $\rho(\widetilde S_1)=2$,
in particular, every geometric fiber is irreducible. 
If we contract the section $\widetilde\G_0=\widetilde S_0\cap\widetilde S_1$ of $\widetilde S_1$ then
we obtain
the surface $S_1$. 
The curve $\G_0$ passes through a finite number of the singular points $Q_1,\ldots,Q_r$ of $Z_0$ ($r\leq 3$), and by $\widetilde F_1,\ldots,\widetilde F_r$
denote the fibers of
$\widetilde S_1$ over these points. 
In small analytic neighborhoods of a general point of $\widetilde \G_0$ and 
a general point of some section $\widetilde E_0$
the variety $Y_0$ has the singularities of types $\CC^1\times\frac1{\beta^1_0}(-\beta^2_0,1)$ and $\CC^1\times\frac1{\beta^2_0}(-\beta^1_0,1)$
respectively. By $F_1,\ldots,F_r$, $E_0$ denote the transforms of $\widetilde F_1,\ldots,\widetilde F_r$, $\widetilde E_0$ on the surface $S_1$ respectively. 
The empty circles are $\widetilde F_1,\ldots,\widetilde F_r$ in the figures of Proposition \ref{pltEDiff}.
The singularities of $\widetilde S_1$ are into ovals.
Note that the self-intersection index $(\widetilde\Gamma_0^2)_{\widetilde S_1}$ was calculated in Lemma \ref{lemm1}.

\begin{proposition}\label{pltEDiff}
Depending on a type of the triple $(S_0,\Diff_{S_0}(0),\G_0)$ we have the following structure
of $(S_1,\Diff_{S_1}(0))$.
\\$1)$ Type $A$, 

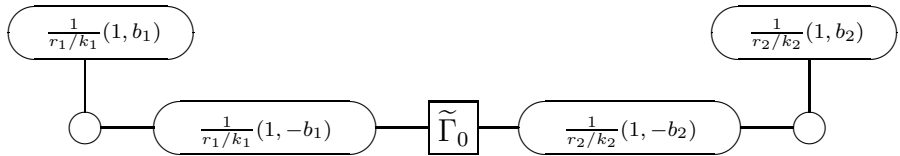
\begin{center}
\begin{picture}(350,70)(0,0)
\put(170,10){\fbox{$\widetilde\Gamma_0$}}
\put(170,15){\line(-1,0){20}}
\put(83,12){\tiny{$\frac1{r_1/k_1}(1,-b_1)$}}
\put(108,15){\oval(84,20)}
\put(66,15){\line(-1,0){20}}
\put(40,15){\circle{12}}
\put(40,21){\line(0,1){20}}
\put(25,49){\tiny{$\frac1{r_1/k_1}(1,b_1)$}}
\put(46,51){\oval(70,20)}
\put(189,15){\line(1,0){15}}
\put(221,12){\tiny{$\frac1{r_2/k_2}(1,-b_2)$}}
\put(246,15){\oval(84,20)}
\put(288,15){\line(1,0){20}}
\put(314,15){\circle{12}}
\put(314,21){\line(0,1){20}}
\put(291,49){\tiny{$\frac1{r_2/k_2}(1,b_2)$}}
\put(312,51){\oval(70,20)}
\end{picture}
Fig. 5
\end{center}

and
$$
\Diff_{S_1}(0)=\frac{k_1-1}{k_1}F_1+\frac{k_2-1}{k_2}F_2+\frac{\beta^2_0-1}{\beta^2_0}E_0.
$$

The pair $(S_1,\Diff_{S_1}(0))$ is toric.\\
$2)$ Type $D_l$ $(l\geq 4)$, 

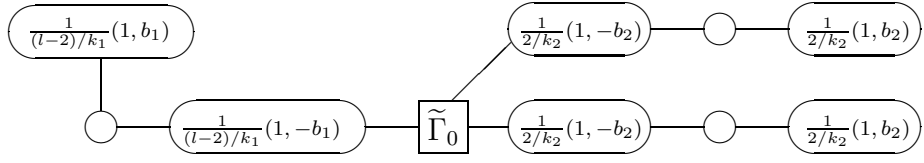
\begin{center}
\begin{picture}(350,65)(0,0)
\put(160,10){\fbox{$\widetilde\Gamma_0$}}
\put(160,15){\line(-1,0){20}}
\put(70,12){\tiny{$\frac1{(l-2)/k_1}(1,-b_1)$}}
\put(103,15){\oval(74,20)}
\put(66,15){\line(-1,0){20}}
\put(40,15){\circle{12}}
\put(40,21){\line(0,1){20}}
\put(12,49){\tiny{$\frac1{(l-2)/k_1}(1,b_1)$}}
\put(40,51){\oval(70,20)}
\put(179,15){\line(1,0){15}}
\put(198,12){\tiny{$\frac1{2/k_2}(1,-b_2)$}}
\put(221,15){\oval(54,20)}
\put(248,15){\line(1,0){20}}
\put(274,15){\circle{12}}
\put(280,15){\line(1,0){20}}
\put(306,12){\tiny{$\frac1{2/k_2}(1,b_2)$}}
\put(325,15){\oval(50,20)}
\put(173,25){\line(1,1){22}}
\put(198,49){\tiny{$\frac1{2/k_2}(1,-b_2)$}}
\put(221,52){\oval(54,20)}
\put(248,52){\line(1,0){20}}
\put(274,52){\circle{12}}
\put(280,52){\line(1,0){20}}
\put(306,49){\tiny{$\frac1{2/k_2}(1,b_2)$}}
\put(325,52){\oval(50,20)}
\end{picture}
Fig. 6
\end{center}

and
$$
\Diff_{S_1}(0)=\frac{k_1-1}{k_1}F_1+\frac{k_2-1}{k_2}F_2+\frac{k_2-1}{k_2}F_3+
\frac{\beta^2_0-1}{\beta^2_0}E_0.
$$

$3)$ Type $E_6$,

\begin{center}
\begin{picture}(350,65)(0,0)
\put(160,10){\fbox{$\widetilde\Gamma_0$}}
\put(160,15){\line(-1,0){20}}
\put(78,12){\tiny{$\frac1{2/k_1}(1,-b_1)$}}
\put(103,15){\oval(74,20)}
\put(66,15){\line(-1,0){20}}
\put(40,15){\circle{12}}
\put(40,21){\line(0,1){20}}
\put(22,49){\tiny{$\frac1{2/k_1}(1,b_1)$}}
\put(40,51){\oval(70,20)}
\put(179,15){\line(1,0){15}}
\put(198,12){\tiny{$\frac1{3/k_2}(1,-b_2)$}}
\put(221,15){\oval(54,20)}
\put(248,15){\line(1,0){20}}
\put(274,15){\circle{12}}
\put(280,15){\line(1,0){20}}
\put(306,12){\tiny{$\frac1{3/k_2}(1,b_2)$}}
\put(325,15){\oval(50,20)}
\put(173,25){\line(1,1){22}}
\put(198,49){\tiny{$\frac1{3/k_3}(1,-b_3)$}}
\put(221,52){\oval(54,20)}
\put(248,52){\line(1,0){20}}
\put(274,52){\circle{12}}
\put(280,52){\line(1,0){20}}
\put(306,49){\tiny{$\frac1{3/k_3}(1,b_3)$}}
\put(325,52){\oval(50,20)}
\end{picture}
Fig. 7
\end{center}

and
$$
\Diff_{S_1}(0)=\frac{k_1-1}{k_1}F_1+\frac{k_2-1}{k_2}F_2+\frac{k_3-1}{k_3}F_3+
\frac{\beta^2_0-1}{\beta^2_0}E_0.
$$

$4)$ Type $E_7$,

\begin{center}
\begin{picture}(350,65)(0,0)
\put(160,10){\fbox{$\widetilde\Gamma_0$}}
\put(160,15){\line(-1,0){20}}
\put(78,12){\tiny{$\frac1{4/k_1}(1,-b_1)$}}
\put(103,15){\oval(74,20)}
\put(66,15){\line(-1,0){20}}
\put(40,15){\circle{12}}
\put(40,21){\line(0,1){20}}
\put(22,49){\tiny{$\frac1{4/k_1}(1,b_1)$}}
\put(40,51){\oval(70,20)}
\put(179,15){\line(1,0){15}}
\put(198,12){\tiny{$\frac1{2/k_2}(1,-b_2)$}}
\put(221,15){\oval(54,20)}
\put(248,15){\line(1,0){20}}
\put(274,15){\circle{12}}
\put(280,15){\line(1,0){20}}
\put(306,12){\tiny{$\frac1{2/k_2}(1,b_2)$}}
\put(325,15){\oval(50,20)}
\put(173,25){\line(1,1){22}}
\put(198,49){\tiny{$\frac1{3/k_3}(1,-b_3)$}}
\put(221,52){\oval(54,20)}
\put(248,52){\line(1,0){20}}
\put(274,52){\circle{12}}
\put(280,52){\line(1,0){20}}
\put(306,49){\tiny{$\frac1{3/k_3}(1,b_3)$}}
\put(325,52){\oval(50,20)}
\end{picture}
Fig. 8
\end{center}

and
$$
\Diff_{S_1}(0)=\frac{k_1-1}{k_1}F_1+\frac{k_2-1}{k_2}F_2+\frac{k_3-1}{k_3}F_3+
\frac{\beta^2_0-1}{\beta^2_0}E_0.
$$

$5)$ Type $E_8$,

\begin{center}
\begin{picture}(350,65)(0,0)
\put(160,10){\fbox{$\widetilde\Gamma_0$}}
\put(160,15){\line(-1,0){20}}
\put(78,12){\tiny{$\frac1{5/k_1}(1,-b_1)$}}
\put(103,15){\oval(74,20)}
\put(66,15){\line(-1,0){20}}
\put(40,15){\circle{12}}
\put(40,21){\line(0,1){20}}
\put(22,49){\tiny{$\frac1{5/k_1}(1,b_1)$}}
\put(40,51){\oval(70,20)}
\put(179,15){\line(1,0){15}}
\put(198,12){\tiny{$\frac1{2/k_2}(1,-b_2)$}}
\put(221,15){\oval(54,20)}
\put(248,15){\line(1,0){20}}
\put(274,15){\circle{12}}
\put(280,15){\line(1,0){20}}
\put(306,12){\tiny{$\frac1{2/k_2}(1,b_2)$}}
\put(325,15){\oval(50,20)}
\put(173,25){\line(1,1){22}}
\put(198,49){\tiny{$\frac1{3/k_3}(1,-b_3)$}}
\put(221,52){\oval(54,20)}
\put(248,52){\line(1,0){20}}
\put(274,52){\circle{12}}
\put(280,52){\line(1,0){20}}
\put(306,49){\tiny{$\frac1{3/k_3}(1,b_3)$}}
\put(325,52){\oval(50,20)}
\end{picture}
Fig. 9
\end{center}

and
$$
\Diff_{S_1}(0)=\frac{k_1-1}{k_1}F_1+\frac{k_2-1}{k_2}F_2+\frac{k_3-1}{k_3}F_3+
\frac{\beta^2_0-1}{\beta^2_0}E_0.
$$

The pair $(S_1, \Diff_{S_1}(0))$ is klt, therefore $g_1\colon (Z_1,S_1)\to (X\ni P)$ is a non-toric plt blow-up.

In cases $A$, $D_l$, $E_6$, $E_7$ and $E_8$ we have a non-plt $1$-, $2$-, $3$-, $4$- and $6$-complement of
$(S_1,\Diff_{S_1}(0))$ respectively.
\begin{proof}
By the construction, the morphism $h_0|_{\widetilde S_1}\colon\widetilde S_1\to\G_0$ is locally toric. 
Therefore, the surface $\widetilde S_1$ has either no singularities in a fiber or only two singularities of types
$\frac1{r_1}(1,b_1)$ and $\frac1{r_1}(1,-b_1)$. Let us show the local calculations.
Consider the singularity at the point $Q_1$ of $Z_0$ such that the curve $\G_0$ contains it.
Let the cone $\langle e_1,e_2,e_3\rangle$ determines locally the variety $Z_0$ in some analytical neighborhood of $Q_1$, $\G_0=V(\langle e_2,e_3\rangle)$ and
$S_0=V(\langle e_3\rangle)$.
According to Proposition \ref{defS} we may assume $e_1=(1,0,0)$.
We locally have $Y_0=T_N(\Delta')$, where
$$
\Delta'=\{\langle \beta,e_1,e_2 \rangle, \langle \beta,e_1,e_3 \rangle,\ \text{their faces}\},
$$
$\beta=\beta^1_0e_2+\beta^2_0e_3$ and $N\cong \ZZ^3$.
Note that $V(\langle \beta\rangle)=\widetilde S_1$ and $\widetilde F_1=V(\langle \beta, e_1 \rangle)$ is the fiber of $\widetilde S_1$ over the point $Q_1$.
Write $(Z_0\ni Q_1)\cong (\CC^3\ni 0)/G$, $(Y_0\ni Q_1')\cong (\CC^3\ni 0)/G_1$, $(Y_0\ni Q_1'')\cong (\CC^3\ni 0)/G_2$,
where $Q_1'=\widetilde F_1\cap\widetilde E_0$, $Q_1''=\widetilde F_1\cap\widetilde S_0$, and 
$G$, $G_1$, $G_2$ are the abelian groups acting freely in codimension 1.
Hence, $\beta^2_0|G|=|G_1|$ and $\beta^1_0|G|=|G_2|$.

Finally, a corresponding complement of the pair 
$(E_0, \Diff_{E_0}(\Diff_{S_1}(0)))$ is extended to a required complement of $(S_1, \Diff_{S_1}(0))$ 
by Proposition \ref{complind}.
\end{proof}
\end{proposition}

{\it Second step.}
Assume that there exists a curve $\G_1\subset S_1$ with the following two properties:
1) $K_{S_1}+\Diff_{S_1}(0)+\G_1$ is an anti-ample divisor, $h_0\colon (\G_1)_{\widetilde S_1}\to\G_0$ is a surjective morphism and
2) $\G_1$ is not a center of any blow-up of
${\mathcal H_0}$, in particular, if $(S_1, \Diff_{S_1}(0))$ is a toric pair then $\G_1$ is its non-toric subvariety.
For convenience, we put $\widetilde \G_1=(\G_1)_{\widetilde S_1}$. 

\begin{lemma}\label{typeA}
The triples $(S_0,\Diff_{S_0}(0),\G_0)$ and $(S_1,\Diff_{S_1}(0),\G_1)$
are of type $A$. Moreover, $\G_1\sim E_0+F_j$ for some index $j$ and 
$\beta^2_0=1$ $($that is, $E_0\not\subset \Supp(\Diff_{S_1}(0)))$.
\begin{proof}
Let us remember that the pairs $(S_1,\Diff_{S_1}(0))$ were classified in Proposition \ref{pltEDiff}, and we will use the same notation.

Put $M=(K_{\widetilde S_1}+\Diff_{\widetilde S_1}(0)+\widetilde\G_1)\cdot \widetilde E_0$. Note that $M<0$.
There are two possibilities:

1) $\widetilde \G_1\sim \widetilde E_0$, $\widetilde E_0\subset \Supp(\Diff_{\widetilde S_1}(0))$ and $\widetilde \G_1\not = \widetilde E_0$;

2) $\widetilde \G_1\not\sim \widetilde E_0$, $\widetilde \G_1\sim a_0\widetilde E_0+\sum_{i=1}^r a_i\widetilde F_i$, where $a_i\in\ZZ_{\geq 0}$ and $a_0\geq 1$.

Suppose that the triple $(S_0,\Diff_{S_0}(0),\G_0)$ does not have type $A$. We will prove that it is impossible.
Proposition \ref{defS} and Lemma \ref{lemm1} imply that $(\widetilde\Gamma_0^2)_{\widetilde S_1}<- \beta^2_0(\Gamma_0^2)_{S_0}\leq - \beta^2_0\leq -1$.
Hence the proper transform of $\widetilde\Gamma_0$ has the self-intersection index $\leq -2$ on the minimal resolution of $\widetilde S_1$.
Consider possibility 1). Then  
$M=-2+\deg(\Diff_{\widetilde E_0}(0))+\frac12\widetilde E_0^2=1-\sum_{i=1}^3\frac1{n_i}+\frac12\widetilde E_0^2$, where $n_i\geq 2$ for all $i$. 
Since the linear system $|\widetilde E_0|$ is movable then 
$\widetilde E_0^2=\widetilde E_0\cdot\widetilde \G_1\geq \frac1{n_{i_1}}+\frac1{n_{i_2}}$ (it is possible that $i_1=i_2$), and hence $M\geq 0$.
Consider possibility 2). If $a_i\geq 1$ for some $i\geq 1$ then it is obvious that $M\geq 0$. Therefore we have to consider the last case 
$\widetilde \G_1\sim a_0\widetilde E_0$, where $a_0\geq 2$. Arguing as in possibility 1) and in its notation we have 
$\widetilde E_0^2=\frac1{a_0}\widetilde E_0\cdot\widetilde \G_1\geq \frac2{a_0}\sum_{k=1}^{a_0} \frac1{n_{i_k}}$, where $i_k\in\{1,2,3\}$,
and hence $M\geq 0$.

Suppose that the triple $(S_0,\Diff_{S_0}(0),\G_0)$ is of type $A$.
We will prove that possibility 1) is not realized, and $a_0=1$, $r=1$, $a_1=1$ in possibility 2).

Let $m_i=r_i/k_i$ be an index of the singularity at the point $\widetilde F_i\cap\widetilde E_0\in \widetilde S_1$, where $i=1,2$.
Lemma \ref{lemm1} implies that
\begin{equation}\label{lab1}
(\widetilde\Gamma_0^2)_{\widetilde S_1}<- \beta^2_0(\Gamma_0^2)_{S_0}\leq - \beta^2_0\Big(\frac1{m_1k_1}+\frac1{m_2k_2}\Big).
\end{equation}
The morphism $h'_0|_{\widetilde S_1}\colon \widetilde S_1\to S_1$ contracts $\widetilde\G_0$ to a point of type $\frac1{m_3}(m_1,m_2)$ and 
$h'_0|_{\widetilde S_1}$ is a toric blow-up corresponding to the weights $(m_1,m_2)$. Hence
\begin{equation}\label{lab2}
(\widetilde \G_0^2)_{\widetilde S_1}=-\frac{m_3}{m_1m_2}.
\end{equation} 
Therefore $m_3>\beta^2_0(m_1/k_2+m_2/k_1)$. 
The toric surface $S_1$ is completely determined by the triple $(m_1,m_2,m_3)$.
For possibility 1) (recall that $\beta_0^2\geq 2$) we obtain the contradiction

\begin{gather*}
M\geq -2+\deg\Big(\Diff_{\widetilde E_0}\Big(\frac{k_1-1}{k_1}\widetilde F_1+\frac{k_2-1}{k_2}\widetilde F_2\Big)\Big)+\frac12\widetilde E_0^2=\\
=-\frac1{m_1k_1}-\frac1{m_2k_2}+\frac{m_3}{2m_1m_2}>0.
\end{gather*}
The same calculations for possibility 2) imply $a_0=1$, and since $\widetilde\G_1$ is an irreducible curve that the same calculations imply
$r=1$ and $a_1=1$.

In order to prove the lemma we must prove only that the plt triple $(S_1,\Diff_{S_1}(0),\G_1)$ is of type $A$.
Assuming the converse: its type differs from type $A$. For instance, let us consider Case
6) of Proposition \ref{defS}, the other cases are considered similarly.
Thus $(S_1,\Diff_{S_1}(0),\G_1)=(\PP_{x_1x_2x_3}(2b_2+1,b_2,1),\frac12\{x_1=0\},\OO_{S_1}(2b_2+1))$, where $b_2\geq 2$. 
Since $\widetilde S_1\to \G_0$ is a toric conic bundle then there are one possibility only: 
$\widetilde S_1\to S_1$ is the weighted blow-up of singularity of type
$\frac1{b_2}(1,1)$ at the point $(0:1:0)$ with the weights $(2b_2+1,1)$. Now $(\widetilde \G_0^2)_{\widetilde S_1}=-\frac{b_2}{2b_2+1}$ 
by equality (\ref{lab2})
and $(\widetilde \G_0^2)_{\widetilde S_1}\leq -(\frac12+\frac1{2b_2+1})$ by inequality (\ref{lab1}). This contradiction concludes the proof.
\end{proof}
\end{lemma}

\begin{remark}
A klt singularity is called {\it weakly exceptional} if there exists its unique plt blow-up (see \cite{Pr2}, \cite{Kud1}).
A two-dimensional klt singularity is weakly exceptional if and only if it is of type $\DDD_n$, $\EEE_6$, $\EEE_7$ or $\EEE_8$.
Lemma \ref{typeA} shows the interesting correspondence of the types. 
\end{remark}

Let $h_1\colon (Y_1,(S_2)_{Y_1})\to (Y_0\supset \widetilde\G_1)$ be a blow-up of the curve 
$\widetilde \G_1$ with an unique exceptional divisor ($\Exc h_1=(S_2)_{Y_1}$),
$(S_1)_{Y_1}\cong (S_1)_{Y_0}$ and the same structure as $h_0$.
The set of all possible blow-ups $h_1$ is denoted by $\mathcal H_1$.

By Proposition \ref{defS} there is 1-complement of $K_{S_1}+\Diff_{S_1}(0)+\widetilde \G_1$ that extends to 1-complement of $K_{Z_1}+S_1$. Therefore we have 1-complement $K_{Y_0}+\widetilde S_1+\widetilde S_0+(D_1)_{Y_0}\sim 0$.
Since $(D_1)_X=(\psi=0\subset (\CC^3\ni 0))/G$ we can slightly change the function $\psi$ keeping all properties. Therefore there is at least a pencil of $(D_1)_{Y_1}$ by proof of Proposition 4.4.1 \cite{PrLect}, and we can assume that $a((S_2)_{Y_1}, (D_1)_X)=-1$.

If $a(S_0, (D_1)_X)\ge 0$ then $S_0\cdot (D_1)_{Z_0}\ge 2\G_0$, hence $K_{S_0}+\Diff_{S_0}((D_1)_{Z_0})$ is nef by Proposition \ref{defS} and $a(S_0, (D_1)_X)\le -1$.

So we have 1-complement $K_{Y_1}+(S_2)_{Y_1}+(S_1)_{Y_1}+(S_0)_{Y_1}+(D_1)_{Y_1}\sim 0$. 
By the cone theorem we have:

1) there exists an divisorial contraction $h'_{1,1}\colon Y_1\to Y_{1,1}$ of $(S_1)_{Y_1}$
onto a curve, $(S_2)_{Y_1}\cong (S_2)_{Y_{1,1}}$;

2) apply $K_{Y_{1,1}}+(S_0)_{Y_{1,1}}+(S_2)_{Y_{1,1}}$-MMP to contract small extremal ray by a small contraction $\varphi_{1,1}$. Put $\Exc \varphi_{1,1}=(F_0)_{Y_{1,1}}$.
Let $\varphi^+_{1,1}$ be a log flip of $\varphi_{1,1}$, $\Exc \varphi^+_{1,1}=(F^+_0)_{Y_{1,2}}$,
$h'_{1,2}\colon Y_{1,1}\dashrightarrow Y_{1,2}$ be a corresponding birational map;

3) there exists a divisorial contraction $h'_{1,3}\colon Y_{1,2}\to Z_2$ of $(S_0)_{Y_{1,2}}$
to a point.

Thus we obtain a birational map 
$h'_1=h'_{1,3}\circ h'_{1,2}\circ h'_{1,1}\colon Y_1\dashrightarrow Z_2$.
Since $(D_1)_{Y_{1,1}}\cdot (F_0)_{Y_{1,1}}=-(K_{Y_{1,1}}+(S_0)_{Y_{1,1}}+(S_2)_{Y_{1,1}})\cdot(F_0)_{Y_{1,1}}>0$, $(D_1)_{Y_{1,1}}$ contains a some fiber of $(S_2)_{Y_{1,1}}$ and $(D_1)_{Y_{1,1}}\not\supset (F_0)_{Y_{1,1}}$ by Proposition \ref{defS}, then the divisor $(D_1)_{Z_2}$ contains the fiber $(F^+_0)_{Z_2}$ and $((S_2)_{Z_2}, \Diff_{(S_2)_{Z_2}}(0))$ is a toric pair
by Shokurov's criterion on the
characterization of toric varieties \cite{Sh2}.
We obtain a non-toric blow-up $g_2\colon (Z_2,S_2)\to (X\ni P)$.

We prove the following proposition.
\begin{proposition}
The pair $(S_2,\Diff_{S_2}(0))$ is toric  $(1$-complementary$)$ with the structure described in Proposition $\ref{pltEDiff}$ $($Type A$)$, $g_2$ is a non-toric plt blow-up.
\end{proposition}

{\it Third step.}
Assume that there exists a curve $\G_2\subset S_2$ with the following two properties:
1) $K_{S_2}+\Diff_{S_2}(0)+\G_2$ is an anti-ample divisor, $h_0\circ h_1\colon (\G_2)_{Y_1}\to\G_0$ is a surjective morphism and
2) $\G_2$ is not a center of any blow-up of
${\mathcal H_1}$, in particular, $\G_2$ is a non-toric subvariety of
$(S_2, \Diff_{S_2}(0))$.

\begin{proposition}
There is no any blow-up 
$h_2\colon (Y_2, (S_3)_{Y_2})\to (Y_1\supset (\G_2)_{Y_1})$ of the curve $(\G_2)_{Y_1}$ with unique exceptional divisor such that $(S_3)_{Y_2}$ is realized by some plt blow-up $g_3\colon Z_3\to (X\ni P)$.
\begin{proof}
Assume the converse. 
Repeat the procedure described in Diagram 4, but with one change, replace the blow-up $g_0\colon Z_0\to X$ by the blow-up
$g_1\colon Z_1\to X$.
Therefore, returning to the main procedure, we can assume that there is 
1-complement $K_{Y_2}+(S_3)_{Y_2}+(S_2)_{Y_2}+(S_1)_{Y_2}+(S_0)_{Y_2}+(D_2)_{Y_2}\sim 0$.
Apply MMP to contract $S_1$ and $S_2$. Let $Y_2\dashrightarrow Y_{2,2}$ be a corresponding birational map. 
If $(S_0)_{Y_{2,2}}$ contains one fiber of $(S_3)_{Y_{2,2}}$ then 
$(S_1)_{Y_2}$ and $(S_0)_{Y_2}$ contain a fiber of $(S_3)_{Y_2}$, a contradiction with log canonicity.
Therefore $(S_0)_{Y_{2,2}}$ contains two fibers of $(S_3)_{Y_{2,2}}$. Then we obtain the contradiction $(K_{(S_3)_{Y_{2,2}}}+\Diff_{(S_3)_{Y_{2,2}}}((S_0)_{Y_{2,2}}+(D_2)_{Y_{2,2}}))\cdot C>0$, where $C$ is any section of the conic bundle $(S_3)_{Y_{2,2}}$.
\end{proof}
\end{proposition}

\par {\bf A2).} Let $(X\ni P)$ be a non-$\QQ$-factorial terminal toric three-dimensional singularity, that is, 
$(X\ni P)\cong (\{x_1x_2+x_3x_4=0\}\subset (\CC^4_{x_1x_2x_3x_4},0))$.

Let $f\colon (Y, E)\to (X\ni P)$ be some non-toric plt blow-up.
Let $\varphi_i\colon X_i\to (X\ni P)$ be two $\QQ$-factorializations, $\Exc\varphi_i=C_i$ ($i=1, 2$).
Let $\psi_i\colon(Y_i, E_i)\to (X_i\ni Q_i)$ be a plt blow-up for some $i$ such that $E_i$ and $E$ define the same discrete valuation of the function field $\mathcal K(X)$, $Q_i$ is a point. The blow-up $\psi_i$ was constructed in the previous case of $\QQ$-factorial singularities, $\rho(E_i)=1$.

Let $Y_i\dashrightarrow \overline Y_i$ be a log flip for the curve $(C_i)_{Y_i}$. 
Considering another value of $i$ we see that $-(E_i)_{\overline Y_i}$ is ample. Therefore $\overline Y_i=Y$ and $\rho(E)=2$.  

We give another construction and prove that $(E,\Diff_E(0))$ is a toric pair by the procedure illustrated on the next diagram (Fig. 10) and defined below.

\[
\xymatrix{& Y_0 \ar[dl]_{h_0} \ar@{-->}[d]^{h'_0} & Y_1 \ar[l]_{h_1} \ar@{-->}[d]^{h'_1}\\
	Z_0 \ar[d]_{g_0}& Z_1 \ar[d]_{g_1} & Z_2 \ar[d]_{g_2}\\
	X& X& X
}
\]\begin{center}
	Fig. $10$
\end{center}

{\it First step.}
Let $g_0\colon (Z_0,S_0)\to (X\ni P)$ be a toric plt blow-up, where $\Exc g_0=S_0$
and $g_0(S_0)=P$ (see Definition \ref{odpdef} and its notation).
Assume that there exists a curve $\G_0\subset S_0$ such that $(S_0,\Diff_{S_0}(0),\G_0)$ is a plt triple (see Definition \ref{triple}).
Such triples are classified in Proposition \ref{defS2}.

\begin{remark}
Note that there exists the divisor $\Omega=\{x_2+\gamma x_1^{w_2/w_1}+\ldots=0\}|_X$ such that $\Omega_{Z}|_S=\G_0$, and it has Du Val singularity of type $\AAA_{w_2/w_1}$, where $\gamma\ne 0$.
\end{remark}

Let $h_0\colon (Y_0,\widetilde S_1)\to (Z_0\supset \G_0)$ be an arbitrary blow-up of the
curve $\G_0$ with an unique exceptional divisor ($\Exc h_0=\widetilde S_1$) as in case {\bf A1)}.
The set of all possible blow-ups $h_0$ is denoted by $\mathcal H_0$.

There are two possibilities. The first possibility is as follows. There is a divisorial contraction of $\widetilde S_0$ to a curve: $h'_0\colon Y_0\to Z_1$, and we obtain a non-toric plt blow-up
$g_1\colon (Z_1,S_1)\to (X\ni P)$, where $\Exc g_1=S_1$ and $g_1(S_1)=P$. The pair $(S_1, \Diff_{S_1}(0))$ is toric as in Proposition \ref{pltEDiff} 1).

The second possibility is when the first possibility is not realized.
The curves $(C_{13})_{Y_0}$ and $(C_{23})_{Y_0}$ (see Definition \ref{odpdef}) generate extremal rays of $\NE(Y_0/X)$ that give small contractions. Let us contract the second one and 
$h'_{0,1}\colon Y_0\dashrightarrow Y_{0,1}$ be a log flip. Let $h'_{0,2}\colon Y_{0,1}\to Z_1$ be a divisorial contraction of $(S_0)_{Y_{0,1}}$ to a point.
Thus we obtain a birational map $h'_0=h'_{0,2}\circ h'_{0,1}\colon Y_0\dashrightarrow Z_1$.
As in case {\bf A1)} 1-complement $K_{S_0}+C_{13}+C_{23}+\G_0$ of $K_{S_0}+\Diff_{S_0}(0)$ extends to 1-complement $K_{Z_0}+S_0+(D_0)_{Z_0}$ 
such that $a((S_1)_{Y_0}, (D_0)_{Z_0}+S_0)=-1$.
Therefore
the divisor $\Diff_{S_1}((D_0)_{Z_1})$ consists of four curves and is 1-complement of $K_{S_1}+\Diff_{S_1}(0)$.
By Shokurov's criterion on the characterization of toric varieties $(S_1, \Diff_{S_1}(0))$ is a toric pair. Thus $g_1\colon Z_1\to (X\ni P)$ is a non-toric plt blow-up.

{\it Second step}.
Assume that there exists a curve $\G_1\subset S_1$ with the following two properties:
1) $K_{S_1}+\Diff_{S_1}(0)+\G_1$ is an anti-ample divisor, $h_0\colon (\G_1)_{\widetilde S_1}\to\G_0$ is a surjective morphism and
2) $\G_1$ is not a center of any blow-up of
${\mathcal H_0}$, $\G_1$ is a non-toric subvariety of $(S_1, \Diff_{S_1}(0))$.

The self-intersection index $\G_0^2$ is calculated by Proposition \ref{defS2}. Lemmas \ref{lemm1} and \ref{typeA} are also true in this case.
So we have 1-complement $K_{Y_1}+(S_2)_{Y_1}+(S_1)_{Y_1}+(S_0)_{Y_1}+(D_1)_{Y_1}\sim 0$. 
By the cone theorem we have:

1) there exists an divisorial contraction $h'_{1,1}\colon Y_1\to Y_{1,1}$ of $(S_1)_{Y_1}$
onto a curve, $(S_2)_{Y_1}\cong (S_2)_{Y_{1,1}}$;

2) apply $K_{Y_{1,1}}+(S_0)_{Y_{1,1}}+(S_2)_{Y_{1,1}}$-MMP to contract small extremal ray, let
$h'_{1,2}\colon Y_{1,1}\dashrightarrow Y_{1,2}$ be a corresponding log flip; 

3) apply $K_{Y_{1,2}}+(S_0)_{Y_{1,2}}+(S_2)_{Y_{1,2}}$-MMP to contract either small extremal ray or
the divisor $(S_0)_{Y_{1,2}}$ onto a curve; we obtain a birational map $h'_{1,3}\colon Y_{1,2}\dashrightarrow Y_{1,3}$ or a morphism $h'_{1,4}\colon Y_{1,3}\to Z_2$ respectively; 
 
4) in the first case of 3) there exists a divisorial contraction $h'_{1,3}\colon Y_{1,3}\to Z_2$ of $(S_0)_{Y_{1,2}}$
to a point.

Thus we obtain a birational map $h'_1\colon Y_1\dashrightarrow Z_2$
and a non-toric blow-up $g_2\colon (Z_2,S_2)\to (X\ni P)$.
The pair
$(S_2, \Diff_{S_2}(0))$ is toric by the same arguments as in case {\bf A1)}.

We prove the following proposition.
\begin{proposition}
The pair $(S_i,\Diff_{S_i}(0))$ is klt and toric $(1$-complementary$)$, $\rho(S_i)=2$, $g_i$ is a non-toric plt blow-up for $i=1, 2$.
\end{proposition}
\end{example}

\begin{example}
In this case we will construct examples of non-toric canonical blow-ups and prove that they are not terminal blow-ups.
Depending on a type of $(X\ni P)$ there are two Cases {\bf B1)} and {\bf B2)}.
\par {\bf B1).} Let $(X\ni P)\cong(\CC_{x_1x_2x_3}^3\ni 0)$.
Let us consider a weighted blow-up $g\colon (Z,S)\to (X\ni P)$ with weights $(w_1,w_2,w_3)$ such that 
$g(S)=P$ (that is, $w_i>0$ for all $i=1,2,3$), where $\gcd(w_1,w_2,w_3)=1$. 
Write $(w_1,w_2,w_3)=(a_1q_2q_3, a_2q_1q_3,a_3q_1q_2)$, where $q_i=\gcd(w_k,w_l)$ and $i,k,l$ are mutually distinct indices from 1 to 3. Then
$$
\Big(S,\Diff_S(0)\Big)\cong
\Big(\PP_{x_1x_2x_3}\big(a_1,a_2,a_3\big),\sum_{i=1}^3
\frac{q_i-1}{q_i}\{x_i=0\}\Big).
$$ 

Assume that $g$ is a canonical blow-up.
\begin{proposition} \label{lemm2}
Let the curve $\G$ be a non-toric subvariety of $(S,\Diff_S(0))$. Assume that 
$\G$ does not contain any center of canonical singularities of $Z$ and  
$-(K_S+\Diff_S(0)+\G)$ is an ample divisor.
Then we have one of the following possibilities for weights $(w_1,w_2,w_3)$ up to permutation of coordinates.

Type $\AAA)$. $(w_1,w_2,w_3)=(a_1q_3,a_2q_3,1)$, $\G\sim\OO_S(a_1+a_2)$.

Type $\DDD)$. $(w_1,w_2,w_3)=(l,l-1,2)$, $(l+1,l,1)$, $(l,l,1)$ and $\G\sim\OO_S(l)$, $\OO_S(2l)$, $\OO_S(2)$ respectively, where 
$l\geq 2$.

Type $\EEE_6)$. $(w_1,w_2,w_3)=(3,2,2)$, $(6,4,3)$, $(5,3,2)$, $(4,2,1)$ and $\G\sim\OO_S(3)$, $\OO_S(2)$, $\OO_S(9)$, 
$\OO_S(3)$ respectively.

Type $\EEE_7)$. $(w_1,w_2,w_3)=(3,2,2)$, $(6,4,3)$, $(9,6,4)$, $(3,3,1)$, $(5,4,2)$, $(7,5,3)$, $(5,3,2)$ and 
$\G\sim\OO_S(3)$, $\OO_S(2)$, $\OO_S(3)$, $\OO_S(2)$, $\OO_S(5)$, $\OO_S(14)$, $\OO_S(6)$ respectively.

Type $\EEE_8)$. $(w_1,w_2,w_3)=(3,2,2)$, $(6,4,3)$, $(9,6,4)$, $(12,8,5)$, $(15,10,6)$, $(5,4,2)$, $(10,7,4)$, $(8,5,3)$ and
$\G\sim\OO_S(3)$, $\OO_S(2)$, $\OO_S(3)$, $\OO_S(6)$, $\OO_S(1)$, $\OO_S(5)$, $\OO_S(10)$, $\OO_S(15)$ respectively.

In all possibilities there is Du Val element $\Omega_Z\in |-K_Z|$ such that $\Omega_Z|_S=\G+\sum_{i=1}^r\gamma_i\G_i$.
Moreover, $\Omega_Z|_S=\G$, except the two possibilities: $(l+1,l,1)$, $\G\sim\OO_S(2l) ($type $\DDD)$ and 
$(5,3,2)$, $\G\sim\OO_S(6) ($type $\EEE_7)$. In these two possibilities we have $\Omega_Z|_S=\G+\G_1$, where $\G_1\sim \OO_S(1)$ and $\OO_S(3)$ respectively.
\begin{proof}
The proof follows from Proposition \ref{cantoric} by enumeration of cases.
\end{proof}
\end{proposition}

\begin{remark}
Proposition \ref{lemm2} is similar to Proposition \ref{defS}. 
Note that there is one-to-one correspondence between the sets $(w_1,w_2,w_3,\G)$ and
the exceptional curves of minimal resolution of Du Val singularity $(\Omega\ni P)$, where $\Omega=g(\Omega_Z)$.
Types in Proposition \ref{lemm2} correspond to Du Val types of the singularity $(\Omega\ni P)$. 
\end{remark}

By Theorem \ref{caninductive} there exists a divisorial contraction 
$h\colon (\widetilde Y,\widetilde E)\to (Z\supset \G)$ for any weights $(\beta_1,1)$ such that
\par 1) $\Exc h=\widetilde E$ is an irreducible divisor and $h(\widetilde E)=\G$;
\par 2) the morphism $h$ is locally toric for a general point of $\G$;
\par 3) if $H$ is a general hyperplane section passing through the general point $Q\in\G$, then $h$ induces the weighted blow-up
of the smooth point $(H\ni Q)$ with weights $(\beta_1,1)$;
\par 4) $h^*S=\widetilde S+\widetilde E$ and $h^*\Omega_{Z}=\Omega_{\widetilde Y}+\beta_1\widetilde E$. 

Apply $K_{\widetilde Y}+\Omega_{\widetilde Y}+\varepsilon \widetilde S$--MMP. Since $\rho(\widetilde Y/X)=2$ and  
$K_{\widetilde Y}+\Omega_{\widetilde Y}+\varepsilon \widetilde S\equiv \varepsilon \widetilde S$ over $X$, then we obtain a sequence of log flips 
$\widetilde Y\dashrightarrow \overline Y$, and after it we obtain the divisorial contraction $h'\colon \overline Y\to Y$ which contracts the proper transform $\overline S$
of $\widetilde S$. 

Thus we obtain a required non-toric blow-up $f\colon (Y,E)\to (X\ni P)$, where $\Exc f=E$ is an irreducible divisor and $f(E)=P$.
Since $K_Y+\Omega_Y=f^*(K_X+\Omega)$ then $f$ is a canonical blow-up.

Finally let us prove that $f$ is a non-terminal blow-up, that is, the singularities of $Y$ are non-terminal. 
We must prove only that
the center of $\overline S$ on $Y$ does not lie in $\Omega_Y$, since $0=a(S,\Omega)$.
Let $\widetilde Y=\overline Y_1\dashrightarrow \overline Y_2\dashrightarrow\ldots\dashrightarrow\overline Y_n=\overline Y$ be a decomposition
of log flip sequence into elementary steps. 
If $\Omega_{\overline Y_i}$ is a nef divisor then by the base point free theorem \cite{KMM} the linear system  
$|m\Omega_{\overline Y_i}|$ gives the birational contraction $h'$ for $m\gg 0$. It contracts the proper transform of
$\widetilde S$ to a point, $i=n$, and this completes the proof. 
Suppose that $\Omega_{\overline Y_i}$ is not a nef divisor.
The cone $\NE(\overline Y_i/X)$ is generated by two extremal rays. 
By $Q_i$, $R_i$ denote them, and to be definite, assume that the ray $R_i$ determines the next step of MMP. By construction, we have $\Omega_{\overline Y_i}\cdot Q_i>0$, 
and hence
$-K_{\overline Y_i}\cdot R_i=\Omega_{\overline Y_i}\cdot R_i<0$. Since $K_{\overline Y_i}\cdot R_i>0$ and the singularities of MMP are canonical,
then the ray $R_i$ gives a log flip (that is, $i<n$), and after it we have $\Omega_{\overline Y_{i+1}}\cdot Q_{i+1}>0$. At the end we obtain that
$\Omega_{\overline Y_j}$ is a nef divisor for some $j$. This completes the proof.

\par {\bf B2).} Let $(X\ni P)\cong (\{x_1x_2+x_3x_4=0\}\subset (\CC^4_{x_1x_2x_3x_4},0))$.
Let us consider a toric canonical blow-up $g\colon (Z,S)\to (X\ni P)$ (see Proposition \ref{cantoric} 3)).

\begin{proposition}\label{lemmB2}
Let a curve $\G$ be a non-toric subvariety of $(S,\Diff_S(0))$. Assume that 
$\G$ does not contain any center of canonical singularities of $Z$ and  
$-(K_S+\Diff_S(0)+\G+\G')$ is an ample divisor, where
$\G'$ is some effective $\QQ$-divisor.
Then $w_1=1$ and $\G\sim\OO_{\PP(w_1,w_2,w_3,w_4)}(w_2)|_S$ up to permutation of coordinates.
There exists Du Val element $\Omega_Z\in |-K_Z|$ such that $\Omega_Z|_S=\G$. In particular, 
$-(K_S+\Diff_S(0)+\G)$ is an ample divisor and
$(\Omega\ni P)$ is Du Val singularity of type 
$\AAA_{w_2}$, where $\Omega=g(\Omega_Z)$.
\begin{proof} 
The proof follows from Proposition \ref{cantoric} 3).
\end{proof}
\end{proposition}

Now we can apply the construction of Case ${\bf B1)}$.
 
Another construction of same non-toric canonical blow-ups is the following one.
Consider a $\QQ$-factorialization $g\colon\widetilde X\to X$ and $\widetilde T=\Exc g$. 
By $G$ denote the center of $E$ on $\widetilde X$. Applying (if necessary) a flop $\widetilde X\dashrightarrow \widetilde X^+$ we may assume that
$G$ is a point.
Let us apply the above mentioned construction in Case ${\bf B1)}$
for singularity $(\widetilde X\ni G)$. We obtain a non-toric canonical blow-up $f\colon Y\to \widetilde X$. 
Let $Y\dashrightarrow Y^+$ be a log flip for the curve $T_Y$.
Thus we obtain a non-toric canonical blow-up
$f^+\colon (Y^+,E^+)\to (X\ni P)$, where $E^+=\Exc f^+$ and $f^+(E^+)=P$.
\end{example}

\section{\bf Main theorems. Case of point}

\begin{example}
Let $(X\ni P)\cong(\CC_{x_1x_2x_3}^3\ni 0)$.
Let us consider the weighted blow-up $g\colon (Z,S)\to (X\ni P)$ with the weights $(15,10,6)$. Then
$$
\Big(S,\Diff_S(0)\Big)\cong \Big(\PP^2,\frac12L_1+\frac23L_2+\frac45L_3\Big),
$$
where $L_i$ are the straight lines, and the divisor $\sum L_i$ is a complement to open toric orbit of $S$.

Let $\Omega=\{x_1^2+x_2^3+x_3^5=0\}\subset (X\ni P)$ be a divisor with Du Val singularity of type $\EEE_8$.
Then $L=\Omega_Z|_S$ is a straight line. Put $P_i=L_i\cap L$. Then the points $P_i$ are non-toric subvarieties 
of $(S,\Diff_S(0))$.

The main difference of structure of non-toric canonical blow-ups from the structure of non-toric plt blow-ups is shown
in the following statements.

1) We have $P_i\in \mt{CS}(Z,\Omega_Z)$ for every $i$. Thus $P_i$ are the centers of some non-toric canonical blow-ups of $(X\ni P)$, that is,
there exists the canonical blow-up $(Y,E_i)\to (X\ni P)$ such that the center of $E_i$ on $Z$ is the point $P_i$ for every $i$.

2) The points $P_i$ are not the centers of any non-toric plt blow-ups of $(X\ni P)$. The proof of this fact is given in Theorem \ref{dim3plt}. 

The origin of this difference is that $S$ is not (locally) Cartier divisor at the points $P_i$ (cf. Theorem \ref{dim2can}).

The straight line $L\in \mt{CS}(Z,\Omega_Z)$ is a center of some non-toric canonical and plt blow-ups of $(X\ni P)$.
As might appear at first sight the class of non-toric canonical blow-ups is much wider than the class of non-toric plt blow-ups, but it is not true.
To construct the non-toric canonical blow-ups, some necessary conditions used implicitly in this example must be satisfied.
Namely, $g$ is a canonical blow-up, $a(S,\Omega)=0$, the straight line $L$ 
does not contain any center of canonical singularities of $Z$.
\end{example}
 
\begin{theorem}\label{dim3plt}
Let $f\colon (Y,E) \to (X \ni P)$ be a plt blow-up of three-dimensional toric terminal singularity,
where $f(E)=P$. Then, either $f$ is a toric morphism, or $f$ is a non-toric morphism described in Section $\ref{sectionex}$.
\begin{proof}
Let $f$ be a non-toric morphism (up to analytical isomorphism). 
Let $D_Y\in |-n(K_Y+E)|$ be a general element for $n\gg 0$. Put $D_X=f(D_Y)$ and $d=\frac{1}{n}$.
The pair $(X,dD_X)$ is log canonical, $a(E,dD_X)=-1$ and $E$ is a unique exceptional divisor with discrepancy $-1$.

Let $(X\ni P)$ be a $\QQ$--factorial singularity.
According to the construction of partial resolution of $(X,dD_X)$ (see Definition \ref{nondegen}) and Criterion \ref{toricplt} 
there exists a toric divisorial contraction $g\colon Z\to X$ such that
it is dominated by partial resolution 
of $(X,dD_X)$ (up to toric log flips), and one of the following Cases I and II occurs. 
 
\par {\it Case} I. The exceptional set $\Exc g=S$ is an irreducible divisor, 
the divisors $S$ and $E$ define the different discrete valuations of the function field
$\mathcal K(X)$, and $g(S)=P$.
By $\G$ denote the center of $E$ on the surface $S$. Then the center $\G$ is a non-toric subvariety of $Z$.
In the other words $\G$ is a non-toric subvariety of
$(S,\Diff_S(0))$. 
If $\G$ is a point then we assume that it does not lie on any one-dimensional orbit of the surface
$S$
(up to analytical isomorphism $(X\ni P)$ of course).
\par {\it Case} II. The variety $Z$ is $\QQ$-gorenstein, hence it is $\QQ$-factorial.  
The exceptional set $\Exc g=S_1\cup S_2$ is the union of two exceptional irreducible divisors, $S_1$, $S_2$ and $E$ define mutually distinct 
discrete valuations of the function field
$\mathcal K(X)$
and $g(S_1)=g(S_2)=P$. To be definite, let $\rho(S_1)=1$, $\rho(S_2)=2$, and $C=S_1\cap S_2$ is a closure of one-dimensional orbit of
$Z$.
By $\G$ denote the center of $E$ on $Z$. In this case $\G$ is a point and a non-toric subvariety of
$(S_1,\Diff_{S_1}(0))$, $\G\in C$, and the curve $C$ has the coefficient 1 in the divisor $\Diff_{S_1}(S_2+dD_Z)$. 
Mori cone $\NE(Z/X)$ is generated by two extremal rays, denote them by $R_1$ and $R_2$. To be definite, let $R_1$ gives the divisorial contraction
which contracts the divisor $S_1$ to some point $P_1$.
Considering toric blow-ups of $P_1$ we may assume that $\Diff_{S_1}(S_2+dD_Z)$ is a boundary in some analytical neighborhood of the point $\G$.

If $R_2$ gives the divisorial contraction which contracts the divisor $S_2$ (onto curve) then it is {\it Case} IIa.
If $R_2$ gives a small flipping contraction then it is {\it Case} IIb. 

Let us consider {\it Case} IIb in more detail. Let $Z\dashrightarrow Z^+$ be a toric log flip induced by $R_2$.
The corresponding objects on $Z^+$ are denoted by the index $^+$. For the toric divisorial contraction
$g^+\colon Z^+\to X$ we have $\rho(S^+_1)=2$, $\rho(S^+_2)=1$. Note that the point $\G^+\in C^+=S_1^+\cap S_2^+$ of $E$ on $Z^+$ can be
a toric subvariety of $(S^+_2,\Diff_{S^+_2}(0))$. 
The morphism $g^+$ is dominated by partial resolution 
of $(X,dD_X)$ (up to toric log flips), and 
the curve $C^+$ has the coefficient 1 in the divisor $\Diff_{S_2^+}(S_1^++dD_{Z^+})$.

Note that the equality $g(\Exc g)=P$ is proved similarly to Theorem  \ref{dim2plt} in both {\it Cases} I and II. 

Now, according to Section \ref{sectionex} the following lemma implies the proof of theorem (for 
$\QQ$--factorial singularities).
\begin{lemma}
It is possible Case $\mt{I}$ only. Moreover, 
$\G$ is a curve and $K_S+\Diff_S(0)+\G$ is a plt divisor.
\begin{proof}
Let us consider Case I.
Write
$$
K_Z+dD_Z+aS=g^*\big(K_X+dD_X\big),
$$
where $a<1$.
Hence
$$
a\big(E,S+dD_Z\big)< a\big(E,aS+dD_Z\big)=-1.
$$
Therefore $\G\subset \LCS(S,\Diff_S(dD_Z))$ and
$-(K_S+\Diff_S(dD_Z))$ is an ample divisor.

Assume that $\G$ is a (irreducible) curve.
We must prove that $K_S+\Diff_S(0)+\G$ is a plt divisor.
Assume the converse.
By the adjunction formula, $\G$ is a smooth curve, and by connectedness lemma 
the divisor $K_S+\Diff_S(0)+\G$ is not a plt one at unique point denoted by $G$.
The point $G$ is a toric subvariety of $(S,\Diff_S(0))$ by Theorem \ref{point}. Moreover, the curve $\G$ is locally a non-toric subvariety at the point $G$ only.
According to the construction of partial resolution \cite{Var} 
there exists the divisorial toric contraction $\widehat g\colon\widehat Z\to Z$ such that
$\Exc \widehat g=S_2$ is an irreducible divisor, $\widehat g(S_2)=G$ and the following two conditions are satisfied.
\par 1).  Put $S_1=S_{\widehat Z}$ and $C=S_1\cap S_2$. Let $c(\G)$ be the log canonical threshold of $\G$ for the pair $(S,\Diff_S(0))$. 
Then ${\widehat g}|_{S_1}\colon S_1\to S$ is the inductive toric blow-up of $K_S+\Diff_S(0)+c(\G)\G$ (see Theorems \ref{inductive} and \ref{dim2plt}),
and the point $\widehat G=C\cap\G_{S_1}$ is a non-toric subvariety of $(S_2,\Diff_{S_2}(0))$.
\par 2). The divisor $\Diff_{S_2}(dD_{\widehat Z}+S_1)$ 
is a boundary at the point $\widehat G$.

Let $H$ be a general hyperplane section of large degree passing through the point $P$.
Then we have $a(S_i,dD_X+hH)=-1$ and $a(S_j,dD_X+hH)>-1$  for some $h>0$, $i\ne j$.
If $i=1$ and $j=2$ then we have the contradiction with Theorem \ref{point} for the pair
$(S_2,\Diff_{S_2}(dD_{\widehat Z}+S_1))$. Hence, we may assume that $i=2$ and $j=1$.
Mori cone $\NE(\widehat Z/X)$ is generated by two rays, denote them by $\widehat R_1$ and $\widehat R_2$.
To be definite, let $\widehat R_2$ gives the contraction $\widehat g$.

At first assume that $\widehat R_1$ gives the contraction $g_1\colon\widehat Z\to Z_1$ which contracts $S_1$ (onto a curve). 
The contraction $g_1$ is an isomorphism for the surface $S_2$, therefore we denote $g_1(S_2)$ by $S_2$ again for convenience.
If $\Diff_{S_2}(dD_{Z_1})$ is a boundary then we have the contradiction 
with Theorem \ref{point} applied for the pair $(S_2,\Diff_{S_2}(dD_{Z_1}))$. If it is not a boundary then we have the following contradiction
\begin{gather*}
0>(1+a(S_1,dD_X+hH))S_1\cdot C_0=\\=(K_{S_1}+\Diff_{S_1}(dD_{\widehat Z}+S_2+hH_{\widehat Z}))\cdot C_0\geq\\
\geq(K_{S_1}+\Diff_{S_1}(0)'+\G_{S_1}+C+C_0)\cdot C_0\geq (-F_1-F_2+\G_{S_1})\cdot C_0\geq 0,
\end{gather*}
where $C_0$ is the closure of one-dimensional orbit of $S_1$, having zero-intersection with $C$, and $F_1$, $F_2$ are the two toric fibers 
(the closures of corresponding one-dimensional toric orbits) 
of the toric conic bundle
$S_1\to g_1(S_1)$, and the divisor $\Diff_{S_1}(0)'$ is a part of $\Diff_{S_1}(0)$ provided that we equate to zero the coefficients of 
$C$ and $C_0$ in $\Diff_{S_1}(0)$.

At last assume that $\widehat R_1$ gives a flipping contraction. Let $\widehat Z\dashrightarrow \widehat Z^+$ be a corresponding toric log flip.
The corresponding objects on $\widehat Z^+$ are denoted by the index $^+$.
If the point $\widehat G^+$ is a non-toric subvariety of $(S_1^+,\Diff_{S_1^+}(0))$ then we have the contradiction 
with Theorem \ref{point} applied for the pair
$(S_1^+,\Diff_{S_1^+}(S_2^+)+\widehat\G^+)$.
Therefore we can assume that the point $G^+$ is a toric subvariety.
If the curve $\widehat\G^+$ is a non-toric subvariety of $(S_1^+,\Diff_{S_1^+}(0))$, then by the inversion of adjunction 
the pair $(S_1^+,\Diff_{S_1^+}(S_2^+)+\widehat\G^+)$ is plt outside $\widehat G^+$, and we have the contradiction with Proposition
\ref{defS}. Thus we have proved that $\widehat\G^+$ and $G^+$ are the toric subvarieties of $(S_1^+,\Diff_{S_1^+}(0))$.
In particular, $S_1^+\cong\PP(1,r_1,r_2)$, where $\gcd(r_1,r_2)=1$ and $(\widehat\G^+)^2=r_1/r_2$.
Considering the divisor $D(\delta)=(d-\delta)D+h(\delta)H$ for some $\delta\geq 0$ and $h(\delta)>0$ ($h(0)=1$) instead of
the divisor $D(0)=dD$, we may assume that the whole construction is satisfied and $a(E,D(\delta))=-1$.

Let $\Diff_{S_2}(D(\delta)-a(S_1,D(\delta))S_1)\geq 0$ (for example, it holds if $a(S_1,D(\delta))<0$). 
Replacing the divisor $H$ by other general divisor with $\widehat G\in \Supp(H_{\widehat Z})$, we may assume that the three following conditions are satisfied:
1) $\Diff_{S_2}(D(\delta)-a(S_1,D(\delta))S_1)\geq 0$;
2) $\widehat G$ is a center of $\LCS(\widehat Z,D(\delta)_{\widehat Z}-a(S_1,D(\delta))S_1-a(S_2,D(\delta))S_2)$;
3) $a(S_2,D(\delta))>-1$. We obtain the contradiction with Theorem \ref{point} for the pair
$(S_2,\Diff_{S_2}(D(\delta)-a(S_1,D(\delta))S_1))$. 

Let $\Diff_{S_2}(D(\delta)-a(S_1,D(\delta))S_1)$ is not an effective divisor.
The curve $\widehat\G^+$ is locally a toric subvariety in some analytical neighborhood of every point of $\widehat Z^+$, therefore 
there exists a blow-up $\overline g\colon(\overline Z\supset \overline S_3)\to (\widehat Z^+\supset \widehat\G^+)$, where $\Exc \overline g=\overline S_3$ 
is an irreducible divisor such that $\overline g(\overline S_3)=\widehat\G^+$ and the following three conditions are satisfied.
\\
A) The morphism $\overline g$ is locally a toric one at every point of $\widehat \G^+$, in particular, $\overline S_1\cong S_1$. 
\\
B) Let $H$ be a general hyperplane section of $\widehat Z^+$ passing through the general point $\widehat Q\in \widehat\G^+$.
Then $\overline g$ induces a weighted blow-up of $(H\ni \widehat Q)$ with weights $(\beta_1,\beta_2)$, and
$\overline g^*S_1^+=\overline S_1+\beta_2\overline S_3$. 
\\
C) Either the divisors $\overline S_3$ and $E$ define the same discrete valuation of the function field
$\mathcal K(X)$ ({\it Case} C1), or
the curve $\overline \G\subset \overline S_3$ being the center of
$E$ on $\overline Z$ is a non-toric subvariety of
$(\overline S_3,\Diff_{\overline S_3}(0))$ ({\it Case} C2).

By $\overline C_0$ and $\overline F$ denote zero-section and a general fiber of $\overline S_3$ respectively.

Let us consider {\it Case} C1. Then $\overline D(\delta)|_{\overline S_3}\sim_{\QQ} a\overline C_0+b\overline F$ by the generality of $D$,
where $b\geq 0$ and
$a=2+a(S_1,D(\delta))/\beta_1-\frac{\beta_2-1}{\beta_2}-\frac{\beta_1-1}{\beta_1}\geq 1+\frac1{\beta_2}$.  
We obtain the contradiction (the calculations are similar to Lemma \ref{lemm1} and Proposition \ref{pltEDiff}) 
\begin{gather*}
0=(K_{\overline S_3}+\Diff_{\overline S_3}(\overline D(\delta)+\overline S_2^+-a(S_1,D(\delta))\overline S_1^+))\cdot \overline C_0\geq \\ \geq -2+1+\frac{r_2-1}{r_2}+
\overline C_0^2 > (r_1-1)/r_2\geq 0.
\end{gather*}

Let us consider {\it Case} C2. If $a(\overline S_3,D(\delta))\leq-1$ then we require the condition $a(\overline S_3,D(\delta))=-1$ to
be satisfied instead of the condition $a(E,D(\delta))=-1$ in the construction of $D(\delta)$, and we obtain similar contradiction as in {\it Case} C1.
Therefore we may assume that $a(\overline S_3,D(\delta))>-1$. Then $\overline \G\sim a\overline C_0+b\overline F$, where either $a\geq 1$, $b\geq 1$, or
$a\geq 2$, $b\geq 0$, or $a=1$, $b=0$, $\overline \G\ne \overline C_0$, $\beta_2\geq 2$.
Continuing this line of reasoning, we have the same contradictions for any possibility of $\overline \G$. 

Now assume that $\G$ is a point. Theorem \ref{point} implies that $\Diff_S(dD_Z)$ is not a boundary in any analytical
neighborhood of $\G$. Moreover, there is unique curve passing through $\G$ with the coefficient $\geq 1$ in the divisor $\Diff_S(dD_Z)$.
It is clear that it is smooth at the point $\G$, it is a non-toric subvariety of $(S,\Diff_S(0))$
and denote it by $T$. 

Let us prove that
$(S,\Diff_S(0)+T)$ is a plt pair.
Let $H$ be a general hyperplane section of large degree passing through the point
$P$ such that $\G\in H_Z$. As above by Theorem \ref{point},
there exist some rational numbers $0<\delta<d$, $h>0$ and the divisor $D'=(d-\delta)D_X+hH$ such that $(X,D')$ is a log canonical pair, 
$\LCS(Z,D'_Z-a(S,D')S)=T$ and $\G$ is a center of
$(Z,D'_Z-a(S,D')S)$. Moreover, we may assume that there are not another centers differing from $\G$ and $T$ by connectedness lemma.
Now, according to the standard Kawamata's perturbation trick, there exists an effective $\QQ$-divisor $D''$ on
$X$ such that the curve $T$ is unique minimal center
of $(Z,D''_Z-a(S,D'')S)$. So, by the previous statement proved (when $\G$ is a curve) $(S,\Diff_S(0)+T)$ is a plt pair.

Let us consider the blow-up $\overline g\colon(\overline Z\supset \overline S_3)\to (Z\supset T)$ for the pair $(X,D')$ which is similar to 
the blow-up $\overline g\colon(\overline Z\supset \overline S_3)\to (\widehat Z^+\supset \widehat\G^+)$, where $\Exc \overline g=\overline S_3$. 
Let $\overline \G\subset \overline Z$ be a center of $E$.
There are two cases $\overline \G=\overline F$, $\overline \G$ is a point, where $\overline F$ is a fiber over the point $\G$.
Applying Lemma \ref{adjnlc} if $\G$ is a point, 
we obtain the contradiction in same way as above 
$$
0=(K_{\overline S_3}+\Diff_{\overline S_3}(\overline{D'}-a(S,D')\overline S))\cdot \overline C_0>0.
$$

Let us prove that {\it Case} II is impossible. Let $H$ be a general hyperplane section of large degree passing through the point $P$.
Then we have $a(S_i,dD_X+hH)=-1$ and $a(S_j,dD_X+hH)>-1$ for some $h>0$. 

Let us introduce the following notation: let $M=\sum m_iM_i$ be the divisor decomposition on irreducible components, then we put
$M^b=\sum_{i\colon m_i>1}M_i+\sum_{i\colon m_i\leq 1}m_iM_i$.

If $i=2$ and $j=1$ then
we obtain the contradiction with Theorem \ref{point} for the pair $(S_1,\Diff_{S_1}(dD_Z+S_2)^b)$.
Therefore $i=1$ and $j=2$. 

Let us consider {\it Case} IIb. If $\G^+$ is a non-toric subvariety of $(S_2^+,\Diff_{S_2^+}(0))$ then
we obtain the contradiction with Theorem \ref{point} for the pair
$(S_2^+,\Diff_{S_2^+}(dD_{Z^+}+S_1^+)^b)$. 
Therefore we assume that $\G^+$ is a toric subvariety of $(S_2^+,\Diff_{S_2^+}(0))$.
The similar (related) case have been considered, when $\G$ was a curve, therefore we do not repeat its complete description.  
By construction, the curve $C^+\subset S^+_1$ is exceptional and contains at most one singularity of $S^+_1$.
Since the pair $(S^+_1,\Diff_{S_1^+}(dD_{Z^+}+hH_{Z^+})$ is not log canonical at the point $\G^+$, then $(dD_{Z^+}+hH_{Z^+})\cdot C^+=1+\sigma$, where $\sigma>0$.
Since the divisor $-K_{S_1^+}$ is a sum of four one-dimensional orbit closures, then 
\begin{gather*}
a(S^+_2,dD_{Z^+}+hH_{Z^+})S_2^+\cdot C^+=\\=(K_{S_1^+}+\Diff_{S_1^+}(dD_{Z^+}+hH_{Z^+}))\cdot C^+\geq \\ \geq -(C^+)^2_{S_1^+}-1-\frac1r_1+1+\sigma
\geq \sigma >0.
\end{gather*}

Since $S_2^+\cdot C^+<0$ then $a(S^+_2,dD_{Z^+}+hH_{Z^+})<0$.
Now, to obtain the contradiction with Theorem \ref{point} for the pair
$(S_1,\Diff_{S_1}(dD_{Z}+hH_{Z}-a(S_2,dD+hH)S_2)^b)$, it is sufficient to decrease the coefficient $h$ slightly (then $a(S_1,dD+hH)>-1$).

Let us consider {\it Case} IIa. Let $g_1\colon Z\to Z_1$ be a contraction of $R_2$.
The contraction $g_1$ is an isomorphism for the surface $S_1$, therefore we denote $g_1(S_1)$ by $S_1$ again for convenience.
If the divisor $\Diff_{S_1}(dD_{Z_1})$ is a boundary then we have the contradiction with 
Theorem \ref{point} for the pair $(S_1,\Diff_{S_1}(dD_{Z_1}))$, and if it is not a boundary then we have the following contradiction
\begin{gather*}
0>(1+a(S_2,dD_X+hH))S_2\cdot C_0=\\=(K_{S_2}+\Diff_{S_2}(dD_Z+S_1+hH_Z))\cdot C_0\geq\\
\geq(K_{S_2}+\Diff_{S_2}(0)'+F+C+C_0)\cdot C_0\geq 0,
\end{gather*}
where $C_0$ is the closure of one-dimensional orbit of $S_2$ having zero-intersection with $C$, and $F$ is a general fiber of the conic bundle
$S_2\to g_1(S_2)$, and 
the divisor $\Diff_{S_2}(0)'$ is a part of $\Diff_{S_2}(0)$ provided that we equate to zero the coefficients of 
$C$ and $C_0$ in $\Diff_{S_2}(0)$.
Note that
the equality $(D_Z|_{S_2}\cdot C)_{\G}\geq 1$ have been applied here (see Lemma \ref{adjnlc}); it is true since $(S_2,C+D_Z|_{S_2})$ is not a plt pair
at the point $\G$ by the construction.
\end{proof}
\end{lemma}

Let $(X\ni P)$ be a non-$\QQ$--factorial singularity, that is, $(X\ni P)\cong (\{x_1x_2+x_3x_4=0\}\subset (\CC^4_{x_1x_2x_3x_4},0))$.
We repeat the arguments given in Section
\ref{sectionex}.
Let $g\colon \widetilde X\to X$ be a $\QQ$-factorialization and let $C=\Exc g\cong \PP^1$. Note that $\widetilde X$ is a smooth variety.
By $G$ denote the center of $E$ on $\widetilde X$. If $G$ is a point then it is a toric subvariety, and hence
the main theorem is reduced to the case of $\QQ$-factorial singularities. If $G=C$ then we consider the flop 
$\widetilde X\dashrightarrow \widetilde X^+$, and we may assume that $G$ is a point by replacing 
$\widetilde X$ by $\widetilde X^+$. 
\end{proof}
\end{theorem}

\begin{theorem}\label{dim32plt}
Let $f\colon (Y,E) \to (X \ni P)$ be a plt blow-up of three-dimensional toric $\QQ$--factorial singularity, where
$f(E)=P$. Then, either $f$ is a toric morphism, or $f$ is a non-toric morphism described in Section $\ref{sectionex}$.
\begin{proof}
We can repeat the proof of Theorem \ref{dim3plt} without any changes in our case. 
Lemma \ref{termrestr} gives some restrictions, when $(X\ni P)$ is a terminal singularity, but it is not used in what follows.
\end{proof}
\end{theorem}

\begin{theorem}\label{dim3can}
Let $f\colon (Y,E)\to (X\ni P)$ be a canonical blow-up of three-dimensional toric terminal singularity, where
$f(E)=P$. Then, either $f$ is a toric morphism $($see Proposition $\ref{cantoric})$, or $f$ is a non-toric morphism described in Section $\ref{sectionex}$.
\begin{proof}
Let $f$ be a non-toric morphism (up to analytical isomorphism). 
Let $D_Y\in |-nK_Y|$ be a general element for $n\gg 0$. Put $D_X=f(D_Y)$ and $d=\frac{1}{n}$.
The pair $(X,dD_X)$ has canonical singularities and $a(E,dD_X)=0$.

Let $(X\ni P)$ be a $\QQ$--factorial singularity.
There is one of two {\it Cases}
I and II described in the proof of Theorem \ref{dim3plt}.
We will use the notation from the proof of Theorem \ref{dim3plt}.
According to Section \ref{sectionex} the following proposition implies the proof of theorem for $\QQ$--factorial singularities.
\begin{proposition}\label{canred}
There exists a toric blow-up $g$ such that we have Case I always, 
the center $\G$ is a curve, $a(S,dD_X)=0$ and $(X\ni P)$ is a smooth point, in particular, $g$ is a canonical blow-up. 
\begin{proof}
Let us consider {\it Case} II.
We may assume that $C\not\subset \Supp(\Sing Z)$. Actually, by taking toric blow-ups with the center $C$ we obtain either the requirement, or {\it Case} I
(that is, there is some blow-up $g$ such that the center of $E$ is a curve and a non-toric subvariety of corresponding exceptional divisor).
Therefore $S_1$ and $S_2$ are Cartier divisors at the point $\G$. Therefore we have
$$
a(E,S_i+dD_Z)\leq a(E,-a(S_i,dD_X)S_i+dD_Z)-1\leq -1
$$
for $i=1,2$

Let $H$ be a general hyperplane section of large degree passing through the point $P$ and let $\G\in H_Z$.
For any $\delta>0$ there exists a number $h(\delta)>0$ such that $(X,D(\delta)=(d-\delta)D_X+h(\delta)H)$ is a log canonical and not plt pair.
Let $D_Z|_S=\sum d_iD^S_i$ be a decomposition on the irreducible components ($S=S_1+S_2$).
If it is necessary we replace the divisor $D_X$ by $D'_X$ in order to $D'_Z|_S=\sum_{i\colon \G\in D^S_i} d_iD^S_i$.
By the generality of $H$ and connectedness lemma, there exists $\delta>0$ with the following two properties.
 
1) The pair $(X,D(\delta))$ defines a plt blow-up $(Y(\delta),E(\delta))\to (X\ni P)$. 

2) By $T$ denote the center of $E(\delta)$ on $Z$. Then, either $T=\G$, or $T$ is a curve provided that
$T\subset S_2$ and $\G\in T$ (note that case $T\subset S_1$ is impossible, since it was proved in {\it Case} I of 
Theorem \ref{dim3plt}).
 
Let $T=\G$. Then we have {\it Case} II of Theorem \ref{dim3plt}, but it was proved that this case is impossible.

Let $T$ be a curve and let $\psi\colon Z\to Z'$ be a contraction of $R_1$.
The morphism $\psi$ contracts the divisor $S_1$ to the point $P_1$. By construction, $K_{S'_2}+\Diff_{S'_2}(0)+T_{S'_2}$ is not a plt divisor
at the point $P_1$, and it was proved in {\it Case} I of Theorem \ref{dim3plt} that this case is impossible.

Let us consider {\it Case} I. Write $K_Z+dD_Z=g^*(K_X+dD_X)+a(S,dD_X)S$, where $a(S,dD_X)\geq 0$. Since $S$ is Cartier divisor at a general point of 
$\G$ then
$$
a(E,S+dD_Z)\leq a(E,-a(S,dD_X)S+dD_Z)-1=-1.
$$
Hence $\G\subset \LCS(S,\Diff_S(dD_Z))$.

Let $a(S,dD_X)=0$. Then $Z$ has canonical singularities.

Assume that $\G$ is a curve. Then $(X\ni P)$ is a smooth point by 
Lemma \ref{lemm4}, which is of independent interest.
\begin{lemma}\label{lemm4}
Let $g\colon (Z,S)\to (X\ni P)$ be a toric canonical blow-up of three-dimensional $\QQ$-factorial terminal toric singularity. 
Assume that there exists a curve $\G\subset S$ such that it is a non-toric subvariety of 
$(S,\Diff_S(0))$, and it does not contain any center of canonical singularities of $Z$. Let $-(K_S+\Diff_S(0)+\G)$ be an ample divisor.
Assume that there exists a divisor $D'_Z\in |-mK_Z|$ for some $m\in\ZZ_{>0}$ such that $\big(Z,\frac1mD'_Z\big)$ is a canonical pair and 
$\big(\frac1mD'_Z\big)|_S=\G+\sum\gamma_i\G_i$, where $\gamma_i\geq 0$ 
for all $i$. Then $(X\ni P)$ is a smooth point.
\begin{proof}
Assume the converse. We suppose that the reader knows the proof of Proposition
\ref{cantoric} 2), and we use its terminology.
We have $a(S,0)=\frac1r(w_3+rw_2-qw_3+rw_1-w_3)-1$. If $a(S,0)=\frac1r$ then we have a contradiction obviously.
Therefore we suppose that $a(S,0)>\frac1r$. For some $j\in\{1,2,3\}$ we have the inequality $\frac1r\geq a(S,0)/N_j$ and one of the two following
requirements: either $P_j\in\mt{CS}(Z)$, or the singularity at the point $P_j$ is of type $\frac1{N_j}(1,-1,0)$, 
where $N_j\ge 2$, $N_1=w_3$, $N_2=rw_1-w_3$, $N_3=rw_2-qw_3$. 

The non-toric curve $\G$ is conveniently represented as $\Gamma=D_Z\cap S$, where $D=(\psi(x_1,x_2,x_3)=0)/\ZZ_r\subset (\CC^3\ni 0)/\ZZ_r(-1,-q,1)$ and $\psi$ is a
quasihomogeneous polynomial with respect to $(N_1,N_2,N_3)$.            

Then $P_j\in\G$, the singularity is of type $\frac1{N_j}(1,-1,0)$  at the point $P_j$ and $N_j/r\ge 1$. Let us prove it. Let $D'=g(D'_Z)$.
If $P_j\not\in\G$ then we have the contradiction $a(S,\frac1mD')<a(S,0)-N_j/r\leq 0$, since $\G$ is a non-toric subvariety.
Let $P_j\in\G$. Then $P_j\not\in\mt{CS}(Z)$, and if $N_j/r<1$, 
then we have the contradiction $a(S,\frac1mD')\le N_j/r-1<0$ since $\G$ is a non-toric subvariety.

Assume that {\it Case} 2A) of Proposition \ref{cantoric} takes place. Then $j=3$. Since $N_3>\max\{N_1,N_2\}$ then the singularity must be isolated at the point $P_3$.
We obtain the contradiction. It is not hard to prove that {\it Case} 2B) of Proposition \ref{cantoric} is impossible.
\end{proof}
\end{lemma}

Assume that $\G$ is a point. Then $\Diff_S(dD_X)$ is a boundary, and hence we obtain the contradiction with Theorem \ref{point}
for the pair $(S,\Diff_S(dD_X))$ and the point $\G$. 

Let $a(S,dD_X)>0$. We will obtain a contradiction.
Note that the number of exceptional divisors with discrepancy 0 is finite for the pair $(X,dD_X)$.
Now we will carry out the procedure consisting of the two steps: 
i1) replacing $dD_X$ by $D(\delta)$ and i2) replacing $(X,dD_X)$ by other pair with canonical singularities
(the variety $X$ is replaced by other variety also). As the result of finite number of steps of this procedure we will obtain
a contradiction. 
Let $H_1$ be a general hyperplane section of large degree containing the center of $S$ on $X$ (at this first step
the point $P$ is this center, and note that this center can be a curve after replacing $X$ as a result of step i2)).
Also we require that 
$(H_1)_Z|_S \subset S$ is an irreducible reduced subvariety (curve) not containing any zero-dimensional orbit of $S$.
This last condition is necessary to our procedure terminates obviously after a finite number of steps. 

Let us consider the numbers $\delta\geq 0$, $h(\delta)\geq 0$ and the divisor $D(\delta)=(d-\delta)D_X+h(\delta)H_1$ 
such that $(X,D(\delta))$ has canonical singularities,
$\G$ is a center of canonical singularities of
$(Z,D(\delta)_Z-a(S,D(\delta))S)$, and
one of the two following conditions are satisfied: either
a1) $a(S,D(\delta))=0$ or a2) $a(S,D(\delta))>0$ and there exists a center of canonical singularities different from $\G$ for the pair
$(Z,D(\delta)_Z-a(S,D(\delta))S)$. Take the maximal number $\delta$ with such properties.
By $E$ again (for convenience) we denote some exceptional divisor with discrepancy 0 for $(X,D(\delta))$ such that its center is $\G$ on $Z$. It is step i1). 

Let $a(S,D(\delta))=0$ and $\G$ be a curve. 
By the above statement $(X\ni P)$ is a smooth point. 
We claim that $h(\delta)=0$, and thus we have the contradiction.
Let us prove it.
Consider the general point $Q$ of $\G$ and the general (smooth)
hyperplane section $H$ passing through this point. Then $(H\ni Q,(D(\delta)_Z)|_H)$ has canonical non-terminal singularities.
This is equivalent to $\mult_Q(D(\delta)_Z)|_H=1$.
Let us apply the construction of non-toric canonical blow-ups from Section \ref{sectionex} to the curve $\G$ provided that $\beta_1=1$.
As the result we obtain the non-toric canonical non-terminal blow-up $(Y'',E'')\to (X\ni P)$. By the above $a(E'',D(\delta))=0$.
Since $\G\not\subset (H_1)_Z$ then the divisor $(H_1)_{Y''}$ contains the center of canonical singularities of  
$Y''$ (see Section \ref{sectionex}) always. Therefore $h(\delta)=0$.

Let $a(S,D(\delta))=0$ and $\G$ be a point. Then $\Diff_S(D(\delta))$ is a boundary and we have the contradiction with Theorem \ref{point}.

Let $a(S,D(\delta))>0$. Let $\widehat X\to X$ be a log resolution of $(X,D(\delta))$.
Let us consider the set $\mathcal E$ consisting of all exceptional divisors $E'$ on $\widehat X$ with the two conditions:
1) $E'$ can be realized by some toric blow-up of $(X\ni P)$
and 2) $a(E',D(\delta))=0$.

Let $\mathcal E=\emptyset$. Hence, if $T\in\mt{CS}(Z,D(\delta)_Z-a(S,D(\delta))S)$ and $T$ is a curve, then $T$ is a non-toric subvariety of 
$(S,\Diff_S(0))$.
Let us consider the variety $T\in\mt{CS}(Z,D(\delta)_Z-a(S,D(\delta))S)$ which is the maximal obstruction to increase a coefficient $\delta$, that is,
if put $\G=T$ then we can more increase the coefficient $\delta$ as the result of step i1).
If $T$ is a curve then we consider $T$ instead of $\G$ and repeat the first step i1) to increase the coefficient $\delta$
(for the sake to be definite, we denote the curve $T$ by $\G$).
If $T$ is a non-toric point lying on some toric orbit, then we are in {\it Case} II. We have proved that {\it Case} II is reduced to {\it Case} I, besides we can assume that 
we consider 
the pair $(X,D(\delta))$ for some $\delta>0$.
If $T$ is a point not lying on any toric orbit then we can consider the point $T$ instead of $\G$ and increase $\delta$ as the result of step i1).
If $T$ is a toric point then we can consider the point $T$ instead of $\G$ and increase $\delta$ and repeat the procedure from the beginning
with the same notation. 
  
Let $\mathcal E\ne\emptyset$.
Let us consider the toric divisorial contraction $g_1\colon Z_1\to (X\ni P)$ which realizes the set $\mathcal E$ exactly.
In particular, $K_{Z_1}+D(\delta)_{Z_1}=g_1^*(K_X+D(\delta))$. Let $P_1$ be a center of $E$ on $Z_1$.
Let us consider locally the pair $(Z_1\supset P_1,D_1=D(\delta)_{Z_1})$ instead of $(X\ni P,D(\delta))$. It is step i2). 
Let us repeat the whole procedure.
We obtain a new divisor $D_1(\delta)$ on $Z_1$. 
Let $a(S,D_1(\delta))=0$. If the center of $S$ on $Z_1$ is a point then we have the contradiction as above.
If the center of $S$ on $Z_1$ is a closure of one-dimensional toric orbit then we have the similar contradiction, but we must use the results of
Section \ref{chpt} (Example \ref{candim1} and Theorem \ref{dim31_can}) to prove $h(\delta)=0$.
Let $a(S,D_1(\delta))>0$. The case $\mathcal E=\emptyset$ is considered as above (the set $\mathcal E$ will be another one).
In the case $\mathcal E\ne\emptyset$ we obtain a toric divisorial contraction
$g_2\colon Z_2\to (Z_1\supset P_1)$, which is constructed similarly to the construction of $g_1$. After it let us repeat the whole procedure.
By construction of partial resolution of $(X,dD_X)$ we obtain some pair $(Z_k,D_k(\delta))$ in a finite numbers of steps such that
$a(S,D_k(\delta))=0$, and hence we have the contradiction.
\end{proof} 
\end{proposition}
Let $(X\ni P)$ be a non-$\QQ$--factorial singularity, that is, $(X\ni P)\cong (\{x_1x_2+x_3x_4=0\}\subset (\CC^4_{x_1x_2x_3x_4},0))$.
According to Section \ref{sectionex} it is sufficient to prove that the analog of Proposition \ref{canred} is satisfied for this singularity.
Arguing as above in Theorem \ref{dim3plt}, the required statement is reduced to the case of
$\QQ$-factorial singularities, this concludes the proof. 
\end{proof}
\end{theorem}

\begin{corollary}
Under the same assumption as in Theorem
$\ref{dim3can}$ the two following statements are satisfied:
\par $1)$ \cite{Kaw1}, \cite{Kawakita1}, \cite{Corti} if $f$ is a terminal blow-up then $f$ is a toric morphism $($see Proposition $\ref{cantoric})$;
\par $2)$ if $f$ is a non-toric morphism then an index of $(X\ni P)$ is equal to $1$.
\end{corollary}


\begin{thebibliography}{99}

\bibitem{BCHM} 
\emph{Birkar C., Cascini P., Hacon C., McKernan J.} Existence of minimal models
for varieties of log general type // J. Amer. Math. Soc. 2010. V. 23. No. 2. P. 405--468.

\bibitem{Corti}
\emph{Corti A.} Singularities of linear systems and 3-fold birational geometry //
Explicit birational geometry of 3-folds, Cambridge LMS. 2000. V. 281. P. 259--312.

\bibitem{Fur}
\emph{Furushima M.} Singular del Pezzo surfaces and analytic
compactifications of $3$-dimensional complex affine space
$\CC^{3}$ //  Nagoya Math. J. 1986. V. 104. P. 1--28.

\bibitem{IshiiPr}
\emph{Ishii I., Prokhorov Yu.~G.} Hypersurface exceptional singularities //
Internat. J. Math. 2001. V. 12. No. 6. P. 661--687.

\bibitem{KaM}
\emph{Katz S., Morrison D.~R.} Gorenstein threefold singularities with small resolutions
via invariant theory for Weyl groups //
J. of Alg. Geom. 1992. V. 1. No. 3. P. 449--530.

\bibitem{Kawakita1}
\emph{Kawakita M.}
Divisorial contractions in dimension three which contract divisors to smooth points //
Invent. Math. 2001. V. 145. P. 105--119.


\bibitem{Kaw0}
\emph{Kawamata Y.} General hyperplane sections of nonsingular flops in dimension 3 //
Math. Res. Lett.
1994.  V. 1. No. 1. P. 49--52.

\bibitem{Kaw1}
\emph{Kawamata  Y.} Divisorial contractions to 3-dimensional terminal
quotient singularities // Higher-dimensional complex varieties
(Trento 1994), de Gruyter, 1996. P. 241--246.

\bibitem{KMM}
\emph{Kawamata  Y.,Matsuda  K.,Matsuki  K.} Introduction to  the
minimal model program //  Algebraic Geometry, Sendai. Adv. Stud.
Pure Math. 1987. V. 10. P. 283--360.

\bibitem{KeM}
\emph{Keel S., McKernan J.} Rational curves on quasi-projective
surfaces // Memoirs AMS 1999. V. 140. No. 669.

\bibitem{Koetal}
\emph{Kollar J. et al} Flips and abundance for algebraic
threefolds // Ast\'erisque 1992 V. 211.

\bibitem{Kollar}
\emph{Koll\'ar J.} Singularities of pairs // Proc. Symp. Pure
Math. 1997. V. 62 Part 1. P. 221--287.

\bibitem{Kud1}
\emph{Kudryavtsev S.~A.} On plt blow-ups // Math. Notes.
2001. V. 69. No. 6. P. 814--819.


\bibitem{Kud2}
\emph{Kudryavtsev S.~A.}
Classification of three-dimensional exceptional log canonical hypersurface singularities. 
{\small II} // Izvestiya: Mathematics. 2004. V. 68.
No. 2. P. 355--364.

\bibitem{KudLd}
\emph{Kudryavtsev S.~A.}
Classfication of exceptional log del Pezzo surfaces with $\delta=1$ //
Izvestiya: Mathematics. 2003. V. 67.
No. 3. P. 461--497.

\bibitem{Mor}
\emph{Morrison D.} The birational geometry of surfaces with rational double points //
Math. Ann. 1985. V. 271. P. 415--438.

\bibitem{Mor2}
\emph{Morrison D.} Canonical quotient singularities in dimension three //
Proceedings of the American Mathematical Society, Mar., 1985. Vol. 93, No. 3. P. 393--396.

\bibitem{Oda}
\emph{Oda T.} Convex bodies and algebraic geometry. An introduction to
the theory of toric varieties // Springer-Verlag 1988.

\bibitem{Pr2}
\emph{Prokhorov Yu.~G.} Blow-ups of canonical singularities //
Algebra. Proc. Internat. Conf. on the Occasion of the 90th
birthday of A.~G.~Kurosh, Moscow, Russia, May 25-30, 1998 /
Yu.~Bahturin ed., Walter de Gruyter, Berlin, 2000. P. 301--317.

\bibitem{PrSh}
\emph{Prokhorov Yu.~G., Shokurov V.~V.} The first main theorem on complements:
from global to local //
Izvestiya Math. 2001. V. 65. No. 6. P. 1169--1196.

\bibitem{PrEll}
\emph{Prokhorov Yu.~G.} Classification of Mori contractions: the case of an elliptic curve //
Izvestiya Math. 2001. V. 65. No. 1. P. 75--84.

\bibitem{PrLect}
\emph{Prokhorov Yu.~G.} Lectures on complements on log surfaces //
MSJ Memoirs 2001. V. 10.

\bibitem{PrCbf}
\emph{Prokhorov Yu.~G.} An application of the canonical bundle formula //
Proceedings of the Steklov institute of mathematics. 2003. V. 241. P. 210--217.

\bibitem{Reid}
\emph{Reid M.} Decomposition of toric morphisms //
Arithmetic and Geometry II (M. Artin and J. Tate eds.).
Birkhauser. Progress in Math. 1983. V.36. P. 395--418.

\bibitem{Sh1}
\emph{Shokurov V.~V.} 3-fold log flips // Izvestiya: Mathematics. 1993. V. 40.
No. 1. P. 95--202.

\bibitem{Sh2}
\emph{Shokurov V.~V.} Complements on surfaces // J. of Math. Sci.
2000. V. 102. No. 2. P. 3876--3932.

\bibitem{Var}
\emph{Varchenko A.~N.}
Zeta-function of monodromy and Newton's diagram //
Invent. Math. 1976. V. 37. P. 253--262.

\end{thebibliography}
\end{document}